\documentclass[10pt]{amsart}

\usepackage{amssymb,amsmath,amscd,a4,graphicx,latexsym,psfrag,epsfig}

\everymath{\displaystyle}

\newtheorem{teo}{Theorem}[section]
\newtheorem{lem}[teo]{Lemma}
\newtheorem{cor}[teo]{Corollary}

\newtheorem{prop}[teo]{Proposition}
\newtheorem{defi}[teo]{Definition}

\newtheorem{remark}[teo]{Remark}
\newtheorem{remarks}[teo]{Remarks}

\newcommand{\mr}{\mathbb{R}}
\newcommand{\mc}{\mathbb{C}}
\newcommand{\mz}{\mathbb{Z}}
\newcommand{\mh}{\mathbb{H}}

\newcommand{\Bb}{{\mathcal B}}
\newcommand{\Cc}{{\mathcal C}}
\newcommand{\Dd}{{\mathcal D}}

\newcommand{\Ff}{{\mathcal F}}
\newcommand{\Gg}{{\mathcal G}}
\newcommand{\Hh}{{\mathcal H}}
\newcommand{\Ii}{{\mathcal I}}

\newcommand{\Kk}{{\mathcal K}}
\newcommand{\Ll}{{\mathcal L}}

\newcommand{\Rr}{{\mathcal R}}

\newcommand{\Tt}{{\mathcal T}}

\newcommand{\Ww}{{\mathcal W}}

\newcommand{\Zz}{{\mathcal Z}}

\newcommand{\C}{{\mathbb C}}

\newcommand{\Z}{{\mathbb Z}}
\newcommand{\R}{{\mathbb R}}
\newcommand{\N}{{\mathbb N}}

\title{Quantum Hyperbolic Geometry}

\author {St\'ephane Baseilhac$^1$ and Riccardo Benedetti$^2$}

\begin{document}

\maketitle

\vspace{0.5cm}

\noindent $^1$ Universit\'e Grenoble I, Institut Joseph Fourier,
UMR CNRS 5582, 100 rue des Maths, B.P. 74, 38402 Saint Martin
d'H\`eres Cedex, FRANCE

\smallskip

\noindent $^2$ Dipartimento di Matematica, Universit\`a di Pisa,
Largo Bruno Pontecorvo 5, I-56127 Pisa, ITALY\medskip

\centerline{Emails: stephane.baseilhac@ujf-grenoble.fr \ \ and\ \
benedett@dm.unipi.it}

\tableofcontents

\begin{abstract}
We construct a new family, indexed by the odd integers $N\geq 1$, of
 $(2+1)$-dimensional quantum field theories called {\it quantum
 hyperbolic field theories} (QHFT), and we study its main structural
 properties. The QHFT are defined for (marked) $(2+1)$-bordisms
 supported by compact oriented $3$-manifolds $Y$ with a properly
 embedded framed tangle $L_\Ff$ and an {\it arbitrary}
 $PSL(2,\C)$-character $\rho$ of $Y \setminus L_\Ff$ (covering, for
 example, the case of hyperbolic cone manifolds). The marking of QHFT
 bordisms includes a specific set of parameters for the space of
 pleated hyperbolic structures on punctured surfaces.  Each QHFT
 associates in a constructive way to any triple $(Y,L_\Ff,\rho)$ with
 marked boundary components a tensor built on the matrix dilogarithms,
 which is holomorphic in the boundary parameters. We establish {\it
 surgery formulas} for QHFT partitions functions and describe their relations
 with the {\it quantum hyperbolic invariants} of \cite{BB1,BB2}
 (either defined for unframed links in closed manifolds and characters
 trivial at the link meridians, or hyperbolic {\it cusped}
 3-manifolds). For every $PSL(2,\mc)$-character of a punctured
 surface, we produce new families of conjugacy classes of ``moderately
 projective" representations of the mapping class groups.
\end{abstract}

\smallskip


\section{Introduction}\label{INTRO}

In this paper we construct a new family $\{\Hh_N\}$, indexed by the odd
integers $N\geq 1$, of $(2+1)$-dimensional quantum field theories
(QFT) that we call Quantum Hyperbolic Field Theories (QHFT). Here,
following \cite{At,T}, by  QFT we mean a functor from a
$(2+1)$-bordism category, possibly non purely topological, to the
tensorial category of finite dimensional complex linear spaces.
\smallskip

The QHFT bordism category is based on triples $(Y,L_\Ff,\rho)$, where
$Y$ is a compact oriented $3$-manifold, possibly with non-empty
boundary $\partial Y$, $L_\Ff$ is a properly embedded {\it framed
tangle} (ie. a framed $1$-dimensional non oriented submanifold) in
$Y$, and $\rho$ is a flat $sl(2,\mc)$-connection on $Y \setminus
L_\Ff$ up to gauge equivalence (ie. a $PSL(2,\C)$-character of $Y
\setminus L_\Ff$), with {\it arbitrary} holonomy at the meridians of
the tangle components. We require furthermore that $L_\Ff$ is
non-empty when $N>1$, and that it intersects each of the boundary
component, if any.  

We will also consider a variant, denoted QHFT$^0$, such that the
tangles $L$ are {\it unframed}, while the characters $\rho$ are
defined on the whole of $Y$, that is the meridian holonomies are
trivial. Finally we consider a ``fusion'' of QHFT and QHFT$^0$
(still denoted QHFT) that incorporates both, by considering tangles
having a {\it framed} part $L_\Ff$ as well as an {\it unframed} one
$L^0$ (see Section \ref{variant}).

The objects of the bordism category are suitably {\it marked}
surfaces. Every such a  QHFT {\it surface} is a diffeomorphism $f:
(S,\Tt,p(\beta)) \to \Sigma$, where $\Tt$ is a so called
``efficient triangulation" of a fixed base oriented surface $S$
with genus $g$ and $r$ marked framed points $p_i$ ($r>0$ and $r>2$
if $g=0$), and $\beta$ is any $PSL(2,\C)$-character of $S\setminus
\{p_i\}$, represented by points $p(\beta)$ in specific parameter
spaces for ${\rm Hom}(\pi,PSL(2,\C))$, built on $\Tt$ and
particularly suited to the QHFT. In fact we construct several such
parameter spaces with small ``residual gauge groups'' acting on
them, and we point out the relations to each other. One of them,
the so called $(-)$-{\it exponential $\Ii$-parameter space}, is
defined in terms of cross-ratios and incorporates the
Bonahon-Thurston shearbend coordinates for pleated hyperbolic
surfaces with punctures.

Every QHFT bordism has marked boundary and is considered as a
``transition'' from its input QHFT surfaces towards the output
ones. We understand that the characters $\rho$ and $\beta$ are
compatible. Every QHFT functor associates to such a transition a
tensor called the {\it amplitude}, defined up to a sign and
multiplication by $N$th roots of unity.
\smallskip

When $Y=W$ is closed (that is $\partial Y = \emptyset$), the
amplitudes $\Hh_N(W,L_\Ff, \rho)$ are numerical invariants called {\it
  partition functions}. The QHFT$^0$ partition functions
$H_N(W,L,\rho)$ coincide with the ``quantum hyperbolic invariants''
constructed in \cite{BB1, BB2}, while the QHFT ones yield new wide
families of numerical invariants, covering interesting geometric
situations, such as compact hyperbolic cone manifolds. We will analyze
the relations between $\Hh_N$ and $H_N$ partition functions.  In
\cite{BB2} we defined also quantum hyperbolic invariants $H_N(M)$ for
non-compact complete hyperbolic 3-manifolds $M$ of finite volume,
ie. for {\it cusped manifolds}. Although these invariants are not
immediately QHFT partition functions, we will show how they can be
obtained in terms of these last. For that we establish a ``surgery
formula'' for quantum hyperbolic invariants of cusped manifolds and
QHFT partition functions that generalizes the one for
Cheeger-Chern-Simons classes, and makes a crucial use of some of
W. Neumann's arguments in \cite{N}, sections 11 and 14.
\smallskip

By restricting QHFT to the trivial bordisms (the cylinders) we get
a new family of conjugacy classes of ``moderately projective''
representations of the mapping class groups of punctured surfaces,
that is, defined up to a sign and multiplication by $N$th roots of
unity.
\smallskip

We stress that we need that any bordism includes a non-empty link,
intersecting each boundary component (so that the QHFT surfaces
have punctures), in order to build a consistent functor when $N>1$
(see \cite{BB2}, Lemma 6.4). Even when the holonomy is trivial
around the punctures, we cannot forget them, in particular for
what concerns the mapping class groups.
\smallskip

We show that QHFT are in fact restrictions to a {\it geometric bordism
  category} of ``universal functors'' called Quantum Hyperbolic
Geometry (QHG). QHG includes the definition of a specific category of
triangulated 3-dimensional pseudomanifolds equipped with additional
structures, and modeled on the functional properties of the {\it
  matrix dilogarithms} studied in \cite{BB2}.  The QHG functors
associate determined tensors to every such a decorated triangulation,
obtained by tracing the matrix dilogarithms supported by each
tetrahedron. The key point is that such tensors are invariant up to
QHG triangulated pseudomanifold isomorphism.  The main step in order
to construct specializations with a strong geometric content, such as
QHFT, consists in converting each QHFT (marked) bordism to a QHG
triangulated pseudomanifold, unique up to QHG triangulated
pseudomanifold isomorphism.
\smallskip

Hence we view this paper as a kind of achievement of the fundation of
the theory initiated in \cite{BB1, BB2}.
\smallskip

A main interest in the QHFT comes from the fact that they relate
classical $3$-dimensional hyperbolic geometry to the world of quantum
field theories, two main themes of low-dimensional topology that
remained essentially disjoint since their spectacular developments in
the early eighties. In particular, the celebrated Kashaev's Volume
Conjecture for hyperbolic knots in $S^3$ \cite{K2} appears as a
special instance of the challenging general problem of understanding
the relations between the asymptotic behaviour of QHFT partition
functions and fundamental invariants coming from differential
geometry, like the Cheeger-Chern-Simons class (see \cite{BB1}, section
5, \cite{BB2}, section 7, and section \ref{ALL} of the present
paper). We plan to face the asymptotics of QHFT partition functions in
future works.
\smallskip

 In \cite{B} the spectrum of the mapping class group representations
 of Section \ref{MODREP} is studied by using geometric quantization of
 the Bonahon-Thurston complex intersection $2$-form. 

A natural problem left unsettled is to determine the relations between
the QHFT and Turaev's Homotopic QFT \cite{Tu2}. This and formulas
describing the behaviour of the QHFT amplitudes under framing changes
will be treated in a sequel to this paper, as they rely mainly on
$R$-matrix computations.

We refer to \cite{BBo} for a discussion about QHG in the framework of
gravity in dimension $3$.

\begin{remark}\label{OKCCS}
{\rm The results of this paper can be repeated almost verbatim to
define Cheeger-Chern-Simons invariants for QHFT bordisms, by replacing
the matrix dilogarithms for $N>1$ with Neumann's extended Rogers
dilogarithm, corresponding to $N=1$ (see Remark \ref{scissors}). This
is described with all details in \cite{BB2} for the case of cusped
manifolds and triples $(W,L,\rho)$ (ie. for quantum hyperbolic
invariants). Hence in the following we concentrate on the quantum
theory $N>1$, which is technically harder.}
\end{remark}

Here is the content of the paper.
\smallskip

The universal QHG functors are defined in section \ref{QHG}, where
we recall also from \cite{BB1,BB2} and \cite{N0,N} (section
\ref{ToGeoSpec} in particular) the notions and results we need. 
\smallskip

The QHFT bordism category is described in section \ref{BORDCAT}, while
its objects, the QHFT surfaces, are developed starting with
section \ref{PARAM}.
\smallskip

The QHFT functors are defined in section \ref{TENSOR}. This includes
the construction of the distinguished QHG triangulated pseudomanifolds
associated to any triple $(Y,L_\Ff,\rho)$ with marked boundary
components, and of the trace tensors computed on them.
The conjugacy classes of moderately projective
representations of the mapping class groups are treated in section
\ref{MODREP}.
\smallskip

The partition functions $\Hh_N(W,L_\Ff,\rho)$ are considered in
section \ref{CLOSED}.  We show in section \ref{ALL} that when $\rho$
is defined on the whole of $W$, $\Hh_N(W,L_\Ff,\rho)$ coincides with
$H_N(W,L\cup L',\rho)$, where $L'$ is a parallel copy of the unframed
link $L$ given by the framing $\Ff$. In section \ref{DEHNF} we prove
the surgery formula for quantum hyperbolic invariants of cusped
manifolds and QHFT partition functions. In fact the QHG
pseudomanifolds used to compute the trace tensors carry certain
{\it cohomological weights} (see Section \ref{ToGeoSpec}) and the 
partition function values actually depend also on them. These weights
play indeed a subtle role in the surgery formulas.  This eventually
leads to realize $H_N(M)$ as the limit of $H_N(M_n,L_n,\rho_n)$, where
$M_n$ is a sequence of closed hyperbolic manifolds converging
geometrically to $M$, $L_n$ is the link of geodesic cores of the
hyperbolic Dehn fillings of $M$ that produce $M_n$, and $\rho_n$ is
the hyperbolic holonomy of $M_n$.  In section \ref{FIBERED} we discuss
alternative computations of the QHFT partition functions for manifolds
that fiber over $S^1$.  For fibred cusped manifolds, this allows in
particular to identify each $H_N(M)$ with a special instance of QHFT
partition function. We conclude with an example in section
\ref{EXAMPLES}.
\medskip

{\bf Convention.}  Unless otherwise stated all manifolds are oriented,
and the boundary is oriented via the convention: {\it last is the
  ingoing normal}. Often we denote ``$=_N$'' the equality of tensors
up to sign and multiplication by $N$th roots of unity.

\section{Universal QHG}\label{QHG}

\subsection{Building blocks}\label{buildbox} The building blocks of 
QHG are {\it flat/charged $\Ii$-tetrahedra} $(\Delta,b,w,f,c)$, and
{\it matrix dilogarithms} $\Rr_N(\Delta(b,w,f,c)) \in {\rm Aut}(\mc^N
\otimes \mc^N)$ defined for every odd positive integer $N$ \cite{BB2}.

\medskip

\noindent{\bf Flat/charged-$\Ii$ tetrahedra.} Consider the half-space
model of the oriented hyperbolic space $\mh^3$, with the group of
direct isometries identified with $PSL(2,\C)$ by the conformal action
on $\partial \overline {\mh}^3= \mc\mathbb{P}^1 = \C \cup\{ \infty \}$
by Moebius transformations. An $\Ii$-{\it tetrahedron} is an oriented
ideal tetrahedron $\Delta$ in $\partial \overline{\mh}^3$ with distinct
ordered vertices $v_0$, $v_1$, $v_2$ and $v_3$ on $\partial
\overline{\mh}^3$.

In fact we consider $\Delta$ as an abstract oriented simplex
equipped with an additional decoration. The ordering of the
vertices is encoded by a {\it branching} $b$, that is, edge
orientations obtained via the rule: each edge points towards the
biggest end-point. Each $2$-face has an induced branching, and a
{\it $b$-orientation}, which is just compatible with that of two
edges on the boundary. We order the $2$-faces $\delta_0$, $\ldots
\delta_3$ by the opposite vertices, and the edges $e_0$, $e_1$,
$e_2$ of $\delta_3$ by stipulating that for $j=0,1$, $v_j$ is the
first end-point of $e_j$.  For exactly two $2$-faces the
$b$-orientation and the boundary orientation are the same. The
{\it $b$-orientation} of $\Delta$ coincides with the given one if
the $b$-orientation of $\delta_3$ looks anti-clockwise from $v_3$.
We give $\Delta$ and each $2$-face $\delta$ a $b$-{\it sign} $*_b$
and $\sigma (\delta)$ respectively, which is $1$ if the two
orientations agree, and $-1$ otherwise.

The hyperbolic structure is encoded by the cross-ratio moduli that
label the edges of $\Delta$.  Recall that opposite edges share the
same cross-ratio moduli. We set $w=(w_0,w_1,w_2)$ with
$w_j=w(e_j) \in \C\setminus \{0,1\}$.
Hence $w_{j+1} = 1/(1-w_j)$ (indices mod($\mz/3\mz$), and $$w_0 = (
v_2- v_1)(v_3-v_0)/(v_2-v_0)(v_3- v_1).$$ We say that the
$\Ii$-tetrahedron $(\Delta,b,w)$ is {\it non-degenerate} if it is of
non zero volume, that is if the imaginary part of each $w_i$ is not
zero; then they have the same sign $*_w = \pm 1$.
\smallskip

\begin{figure}[ht]
\begin{center}
 \includegraphics[width=8cm]{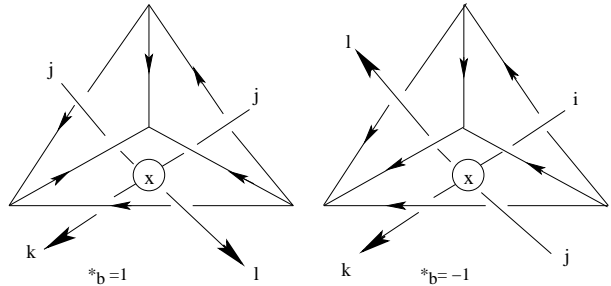}
\caption{\label{CQDidealtensor} $\Ii$-tetrahedra and dual
encoding.}
\end{center}
\end{figure}

It is very convenient to encode $(\Delta,b,w)$ in dual terms. In
Figure \ref{CQDidealtensor} we show the $1$-skeleton of the dual
cell decomposition of Int$(\Delta)$ ($x$ and the indices $i$, $j$,
$k$ and $l$ are considered below). It is understood that (dual) edges
without arrows are incoming at the crossing. Note that an oriented
edge is outgoing exactly when the $b$-sign of the dual $2$-face is
$1$.

A {\it flat/charged} $\Ii$-tetrahedron is an $\Ii$-tetrahedron
equipped with a {\it flattening} $f$ and a {\it charge} $c$, two
notions first introduced in \cite{N0} and \cite{N}. Flattenings
and charges are $\Z$-valued functions defined on the
edges of $\Delta$ that take the same value on opposite edges and
satisfy the following properties, respectively (where $\log$ has
the imaginary part in $]-\pi,\pi]$):
\begin{enumerate}
\item[F.] {\it Flattening condition:} ${\rm l}_0 + {\rm l}_1 +{\rm
    l}_2 = 0$, where
\begin{equation}\label{clb}
{\rm l}_j = {\rm l}_j(b,w,f)= \log(w_j) + \sqrt{-1}\pi f_j,
\end{equation}
\item[C.] {\it Charge condition:} $c_0+c_1+c_2=1$.
\end{enumerate}
We call ${\rm l}_j$ a {\it (classical) log-branch}, and for every
odd $N>1$ we define the \emph{(level $N$) quantum log-branch} as
\begin{equation}\label{qlb}
{\rm l}_{j,N} = \log(w_j) + \sqrt{-1}\pi (N+1)(f_j -*_bc_j).
\end{equation}
The bijective map $$({\rm l}_0,{\rm l}_1,{\rm l}_2) \mapsto
\left(w_0;\frac{{\rm l}_0-\log(w_0)}{\sqrt{-1}\pi},\frac{{\rm
l}_1-\log(w_1)}{\sqrt{-1}\pi}\right)$$ yields an identification of the
set of log-branches on $(\Delta,b)$ with the Riemann surface
$\widehat{\mc}$ of the maps $w_0\mapsto
(\log(w_0)+\varepsilon\pi\sqrt{-1},\log((1-w_0)^{-1})+
\varepsilon'\pi\sqrt{-1})$, with $\varepsilon,\varepsilon'\in
\{0,1\}$. Similarly, the set of triples $(w_0',w_1',w_2')$ with
$w_j'=\exp({\rm l}_{j,N}/N)$ gets identified with the quotient
covering by $N\mz\times N\mz$. These two spaces are the domains of
definition for the matrix dilogarithms to be described next, when
$N=1$ and $N>1$ respectively.
\medskip

\noindent {\bf Matrix dilogarithms.}  Denote by $\log$ the standard
branch of the logarithm, which has the imaginary part in
$]-\pi,\pi]$. Recall that the space of triples of log-branches on a
    branched oriented tetrahedron is identified with the Riemann surface
    $\widehat{\mc}$. 
\smallskip

For $N=1$, we forget the integral charge $c$, so that
$\mathcal{R}_1$ is defined on flattened $\Ii$-tetrahedra:
\begin{equation}\label{symRmatdil}
\mathcal{R}_1(\Delta,b,w,f) =
\exp\biggl(\frac{*_b}{\pi\sqrt{-1}}\mathcal{R}(w_0;f_0,f_1)\biggr) =
\exp\biggl(\frac{*_b}{\pi\sqrt{-1}}\biggl(-\frac{\pi^2}{6} -
\frac{1}{2} \int_0^{w_0} \biggl( \frac{{\rm l}_0(t)}{1-t} - \frac{{\rm
    l}_1(t)}{t} \biggr) \ dt \biggr) \biggr)
\end{equation}
where ${\rm l}_0(t)=\log(t) + \sqrt{-1}\pi f_0$ and ${\rm
  l}_1(t)=\log((1-t)^{-1}) + \sqrt{-1}\pi f_1$. The map
  $\mathcal{R}:\widehat{\mc}\rightarrow \mc/\pi^2\mz$ is holomorphic,
  and takes values in $\mc/2\pi^2\mz$ on the component with even
  valued flattenings (\cite{N}, Proposition 2.5).  \smallskip

For $N=2m+1>1$ and every complex number $x$ set $x^{1/N} =
\exp(\log(x)/N)$ ($0^{1/N} = 0$ by convention). Put
\begin{equation}\label{funcg}
g(x) := \prod_{j=1}^{N-1}(1 - x\zeta^{-j})^{j/N}.
\end{equation}
The function $g$ is defined over $\mc$, and analytic over the
complement of the rays from $x =\zeta^k$ to infinity, $k=1,\ldots,\
N-1$. Set $h(x) := g(x)/g(1)$ (we have $\vert g(1) \vert = N^{1/2}$).
For any $u',v'\in \C$ satisfying $(u')^N + (v')^N = 1$ and any $ n \in
\mathbb{N}$, let
$$\omega(u',v'\vert n) = \prod_{j=1}^n \frac{v'}{1-u'\zeta^j}$$ with
$\omega(u',v'\vert 0) = 1$ by convention. The functions $\omega$ are
periodic in their integer argument, with period $N$. Write
$[x]=N^{-1}(1-x^N)/(1-x)$. Given a flat/charged $\Ii$-tetrahedron
$(\Delta,b,w,f,c)$, set
\begin{equation}\label{qvar}
w_j'=\exp({\rm l}_{j,N}/N)
\end{equation}
with ${\rm l}_{j,N}$ as in (\ref{qlb}). Put the standard tensor
product basis on $\C^N\otimes \C^N$. The matrix dilogarithm of level
$N>1$ is the tensor valued function of flat/charged $\Ii$ tetrahedra
defined by
\begin{equation}\label{symqmatdil}
\mathcal{R}_N(\Delta,b,w,f,c) =
\bigl((w_0')^{-c_1}(w_1')^{c_0}\bigr)^{\frac{N-1}{2}}\
(\Ll_N)^{*_b}(w_0',(w_1')^{-1})\in {\rm
Aut}(\C^N\otimes\C^N)
\end{equation}
where (recall that $N=2m+1$)
$$\begin{array}{lll} \Ll_N(u',v')_{k,l}^{i,j}& = & h(u')\
\zeta^{kj+(m+1)k^2}\ \omega(u',v'\vert i-k) \ \delta(i + j - l) \\
\bigl( \Ll_N(u',v')^{-1}\bigr)_{k,l}^{i,j}& = & \frac{[u']}{h(u')}\
\zeta^{-il-(m+1)i^2}\ \frac{\delta(k + l - j)}{\omega(u'/\zeta,v'\vert
k-i)}
\end{array}$$
with $\delta$ the Kronecker symbol with period $N$, that is $\delta(n)
= 1$ if $n \equiv 0$ mod($N$), and $\delta(n) = 0$ otherwise. Note
that we use the branching in order to associate an index among $i$,
$j$, $k$ and $l$ to each $2$-face of $\Delta$. The rule is shown in
Figure \ref{CQDidealtensor}.

Up to multiplication by $N$th roots of unity, the map $\mathcal{R}_N$
is holomorphic on $\widehat{\mc}$, because $w_j'= \exp((1/N) ({\rm
  l}_j - *_bc_j)) \exp((f_j - *_bc_j)\pi i)$ and the parity of $f_j$
is unaltered when we move $l_j$ continuously. The ambiguity of
$\mathcal{R}_N$ is a consequence of the jumps of $g$ along the cuts
from the $\zeta^k$ to infinity.
\medskip

\noindent {\bf On matrix dilogarithms via geometric quantization.}  As
explained in Remark \ref{OKCCS}, in this paper we concentrate on the
{\it quantum} matrix dilogarithms (ie.  with $N>1$). These were
derived in \cite{BB2} from the Kashaev's 6j-symbols for the cyclic
representation theory of a Borel quantum subalgebra $\Bb_\zeta$ of
$U_\zeta(sl(2,\mc))$, where $\zeta=\exp(2i\pi/N)$. Let us outline here
very briefly an alternative construction based on {\it geometric
quantization of $\Rr_1$}, thus clarifying their geometric origin (for
details see \cite{B}).
\smallskip

\begin{figure}[ht]
\begin{center}
\includegraphics[width=6cm]{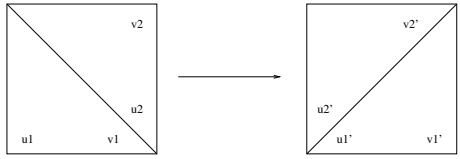}
\caption{\label{diagex} A diagonal exchange.}
\end{center}
\end{figure}

\begin{figure}[ht]
\begin{center}
 \includegraphics[width=7cm]{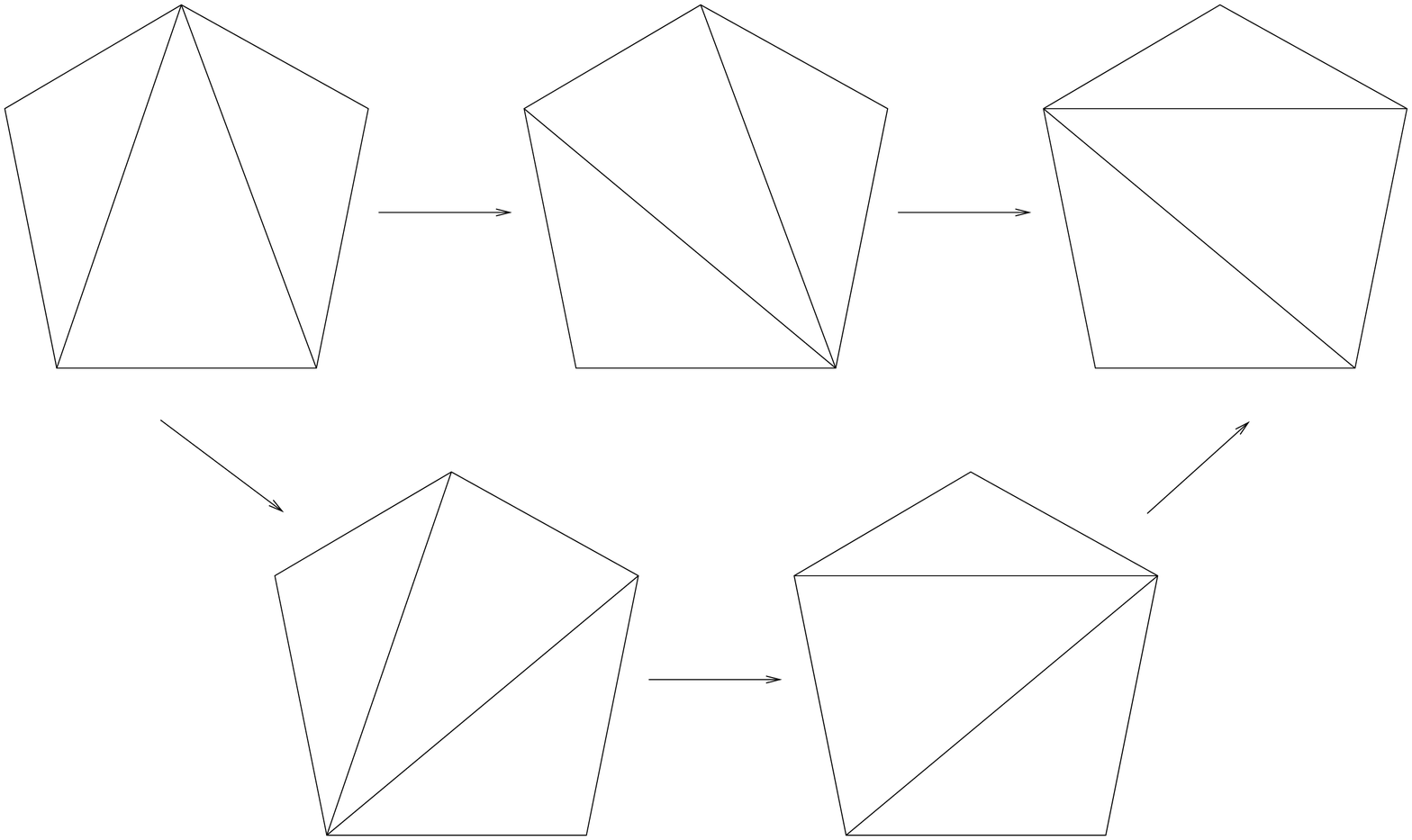}
\caption{\label{figpent} The pentagon relation.}
\end{center}
\end{figure}

Consider an abstract oriented quadrilateral $Q$, triangulated by two
triangles.  Order the triangles of $Q$, and associate to each a copy
$X_i$, $i=1,2$, of $X=\{(u,v,w)\in(\mc^*)^3\mid uvw=-1\}$, where $u$,
$v$ and $w$ correspond to the corners, ordered cyclically by using the
orientation. Let us regard $\mc$ as the Riemann surface of $\log$,
ie. 
$$\mc=\{(z;p) \in \mc^*\times 2\mz \}/((z+i0;p)\sim (z-i0;p+2),
\forall z\in (-\infty;0)) \ . $$ As in (\ref{clb}), we set ${\rm
l}(z;p)= \log(z)+\sqrt{-1}\pi p $.  Let
$$\widehat{X}=\{((u;p),(v;q),(w;r)) \in \tilde{X}\mid {\rm l}(u;p)+{\rm
l}(v;q)+{\rm l}(w;r)=0\}$$ a subspace of the universal covering
$\tilde{X}$ of $X$, and put
$$\widehat{\mc}(z;p,q)=\{({\rm l}(z;p),{\rm
l}((1-z)^{-1};q))\}=\{((z;p),((1-z)^{-1};q),(1-z^{-1};r))\in
\widehat{X}\} \ . $$ Denote $b_Q=d\log u_1 \wedge d\log v_1 + d\log
u_2 \wedge d\log v_2$ the canonical complex symplectic form on
$X_1\times X_2$, and
$$\theta_Q=(1/2)\left(({\rm l}(u_1;p_1)d\log v_1 - {\rm
  l}(v_1;q_1)d\log u_1) + ({\rm l}(u_2;p_2)d\log v_2 - {\rm
  l}(v_2;q_2)d\log u_2))\right)$$
the symplectic potential for the lift of $b_Q$ to $\widehat{X}_1
\times \widehat{X}_2$. Consider the four punctured sphere $S_Q^2$
obtained by gluing along the boundary in the natural way $Q$ and the
quadrilateral $Q'$ obtained from it by exchanging the diagonal. The
$S_4$-action on the vertices of $S_Q^2$ reorders the copies $X_i$
($i=1,\ldots,4$) attached to the triangles, and induces the usual action
of $PSL(2,\mz)$ on each $\widehat{X}_i\cong \mc^2$, with standard
basis ${\rm l}(u_i;p_i)$ and ${\rm l}(v_i;q_i)$, by
symplectomorphisms. 
\begin{prop}\label{flipform} 
There is a canonical family $\{\tilde{\phi}_Q^n\}_{n\in \mz}:
(\widehat{X}_1 \times \widehat{X}_2,b_Q)\rightarrow (\widehat{X}_3
\times \widehat{X}_4,b_{Q'})$ of analytic symplectomorphisms
equivariant for the $S_4$-action on $S_Q^2$, given by (see
Figure \ref{diagex})
$$ \tilde{\phi}_Q^n(({\rm l}(u_1;p_1),{\rm l}(v_1;q_1)),({\rm
    l}(u_2;p_2),{\rm l}(v_2;q_2))) = (({\rm l}(u_1';p_1'),{\rm
    l}(v_1';q_1')),({\rm l}(u_2';q_2'),{\rm l}(v_2';q_2'))),$$ where we
    identify $\widehat{X}$ with $\mc^2$ (first two coordinates) and
$$\left\lbrace\begin{array}{lll}{\rm l}(u_1';p_1') & = & {\rm
l}(u_1;p_1) -{\rm l}(u_2';p_2') \\{\rm l}(v_1';q_1')& = & {\rm
l}(v_1;q_1) + {\rm l}(u_2;p_2)\end{array}\right.\quad
\left\lbrace\begin{array}{lll} {\rm l}(u_2';p_2') & = & {\rm
l}(u_2;p_2)-{\rm l}((1-u_1v_2)^{-1};n)\\{\rm l}(v_2';q_2')& = & {\rm
l}(v_2;q_2) +{\rm l}(u_1';p_1')+{\rm l}(v_1';q_1')
\end{array}\right. .$$
The maps $\tilde{\phi}_Q^n$ satisfy the pentagon relation, that is,
$\tilde{\phi}_{Q_1}^{n_1}\circ\tilde{\phi}_{Q_2}^{n_2}
\circ\tilde{\phi}_{Q_3}^{n_3}\circ(\tilde{\phi}_{Q_4}^{n_4})^{-1}\circ
(\tilde{\phi}_{Q_5}^{n_5})^{-1}={\rm Id}_{\widehat{X}^3}$, where $Q_i$
has diagonal the $i$th edge exchanged in Figure \ref{figpent}, for the
positive cyclic ordering starting from the top left pentagon. Moreover
$\theta_Q -\tilde{\phi}_n^*\theta_{Q'}= -d\mathcal{R}$, where
$\mathcal{R}:\widehat{\mc}(u_1v_2;m,n)\rightarrow \mc/2\pi\mz$ is the
extended dilogarithm of (\ref{symRmatdil}), and
$\widehat{\mc}(u_1v_2;m,n)=\{({\rm l}(u_1;p_1)+{\rm l}(v_2;q_2),{\rm
l}((1-u_1v_2)^{-1};n))\}$ is attached to the diagonal of $Q$.
\end{prop}
By ``canonical'' we mean that $\{\tilde{\phi}_Q^n\}$ is the unique
satisfying some natural properties related to configuration spaces of
points in $\mc P^1$. In fact, $\widehat{\mc}(z;p,q)$ is isomorphic to
the moduli space of similarity classes of triangles in the complex
plane endowed with lifts to $\mr$ of the angles, or, equivalently, to
the moduli space of isometry classes of hyperbolic ideal tetrahedra
with even valued flattenings. The complex symplectic form $d\log u
\wedge d\log v$ restricts to the {\it real} symplectic form
$w=-d\log(z) \wedge d\log(1-z)$ on $\widehat{\mc}(z;p,q)$. It is the
differential version of the extended complex Dehn invariant
$\widehat{\delta}:\widehat{\mc}(z;p,q) \rightarrow \mc \wedge \mc$,
$(z;p,q)\mapsto {\rm l}(z;p)\wedge_\mz {\rm l}(1-z;q)$ of \cite{N}.

When Im$(z)\ne 0$ the form $w$ is K\"ahler for the usual complex
structure on $\mr^2$. Replace $\widehat{\mc}$ with the ramified
covering $\widehat{\mc}_N$ of $\mc^*$ obtained by taking the quotient
with $N\mz\times N\mz$. Put the forms $Nw$ and $N\theta$ on
$\widehat{\mc}_N$. Standard half-form quantization of
$\widehat{\mc}_N$ produces a $N$-dimensional vector space $\Gamma^N$
of sections of a line bundle over $\widehat{\mc}_N$ (the {\it coherent
  states}), and the maps $\tilde{\Phi}_Q^n:
\Gamma_1^N\otimes\Gamma_2^N\rightarrow \Gamma_3^N\otimes \Gamma_4^N$
induced by pull back via $\{\tilde{\phi}_Q^n\}$, that is such that
$\tilde{\Phi}_Q^n(s)=s\circ \tilde{\phi}_Q^n$ for any $s \in
\Gamma_1^N\otimes\Gamma_2^N$, coincide with the matrix dilogarithms
$\Rr_N$ for some suitable basis. In particular, the pentagon relation for
$\tilde{\phi}_Q^n$ lifts to the five term identities mentionned after
Proposition \ref{fund_func} below.

\subsection{QHG triangulated pseudomanifolds}\label{QHGPM} We
restrict the discussion to pseudomanifolds for simplicity, but all
what follows makes sense for arbitrary singular $3$-cycles whose non
manifold locus is of codimension $\geq 2$. By a pseudomanifold $Z$,
possibly with non-empty boundary $\partial Z$, we mean a compact
oriented polyhedron with at most a finite set of non manifold
points. The boundary is a pseudo-surface.

A {\it QHG triangulated pseudomanifold} $(Z,\Tt)$ is a pseudomanifold
$Z$ obtained as the quotient of a finite family
$\mathcal{Z}=\{(\Delta^i,b^i,w^i,f^i,c^i)\}$ of flat/charged
$\Ii$-tetrahedra, via a system of orientation reversing simplicial
identifications of pairs of $2$-faces such that the branchings
match. The resulting triangulation $T$ of $Z$ is endowed with a
{\it global branching} $b$, and is possibly singular (multiply
adjacent as well as self adjacent tetrahedra are allowed). The set of
non manifold points of $Z$ is contained in the set of vertices of
$T$. We do not impose for the moment any global constraint on the
moduli, flattening and charges. Hence $Z$ is equipped with a {\it
rough flat/charged $\Ii$-triangulation} $\ \Tt=(T,b,w,f,c)$, where
$w=\{w^i\}$ and so on.

\begin{figure}[ht]
\begin{center}
\includegraphics[width=6cm]{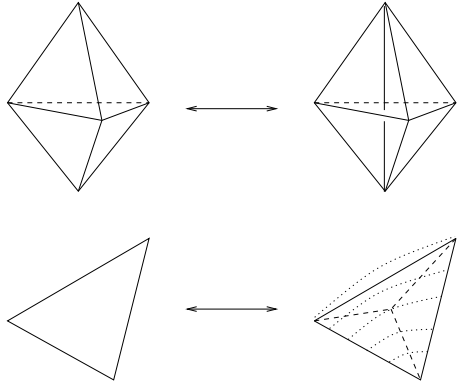}
\caption{\label{CQDfigmove1} The moves on naked singular
triangulations.}
\end{center}
\end{figure}

Next we define the {\it QHG triangulated pseudomanifold
isomorphisms}. Fix a QHG triangulated pseudomanifold $(Z,\Tt)$,
and let $\epsilon_T : E(\Zz)\to E(T)$ be the identification map of
edges. We define the {\it total modulus, log-branch} and {\it
charge} of an edge $e$, respectively, by:
\begin{eqnarray}\label{totw} W_\Tt(e)=\prod_{h \in
\epsilon_T^{-1}(e)} w(h)^{*_{b}}\nonumber\\ L_\Tt(e)= \sum_{h\in
\epsilon_T^{-1}(e)} *_{b}\ {\rm l}(h)\\ C_\Tt(e)= \sum_{h\in
\epsilon_T^{-1}(e)} c(h)\ \ \ \nonumber
\end{eqnarray}
where $*_{b}=\pm 1$ according to the $b$-orientation of the
tetrahedron in $Z$ that contains $h$, $w(h)$ is the cross-ratio
modulus at $h$, $l(h)$ the log-branch at $h$, and $c(h)$ the
charge at $h$.
\begin{remark}\label{W_mean} {\rm It is easily seen that $W_\Tt(e)$ is
a cross ratio for the four `extremal' points on
$\partial \overline{\mathbb{H}}^3$ determined by gluing oriented
hyperbolic ideal tetrahedra with moduli $w(h)^{*_{b}}$ along a
common edge (by continuity $W_\Tt(e)=1$ for a degenerate
quadrilateral with
    three distinct vertices).}
\end{remark}
It is well-known that any two arbitrary {\it naked} triangulations $T$,
$T'$ of $Z$ with the same boundary triangulation can be connected
(keeping the boundary triangulation fixed) by a finite sequence of the
{\it local moves} shown in Figure \ref{CQDfigmove1}, the {\it $2
\leftrightarrow 3$ move} (top) and the {\it bubble move} (bottom).

For any such a local move $T \leftrightarrow T'$ we have two
triangulations of a same portion of a polyhedron $Q$. Assume that
both $T$ and $T'$ extend to portions $(Q,\Tt)$ and $(Q,\Tt')$ of
QHG triangulated pseudomanifolds. We have to specify the
admissible {\it QHG transits} $\Tt \leftrightarrow \Tt'$. In any
case we require that they are local (that is, the portions
complements remain unchanged), and that the branchings coincide at
every common edge of $T$ and $T'$.

For the $2 \leftrightarrow 3$ move we also require that at every
common edge $e$ as above the total modulus, log-branch and charge
coincide
\begin{equation}\label{preserveq}
W_\Tt(e)=W_{\Tt'}(e),\quad L_\Tt(e)=L_{\Tt'}(e),\quad C_\Tt(e)=C_{\Tt'}(e).
\end{equation}
The same rule restricted to the total modulus and total log-branch
holds also for the bubble move; however, the total charge behaves
in a different way: any bubble transit $\Tt \leftrightarrow \Tt'$
includes a marked edge $e$ common to $T$ and $T'$. Referring to
the bottom of Figure \ref{CQDfigmove1} we require that
$C_\Tt(e)=C_{\Tt'}(e)-2$, while for the other two common edges the
total charges are unchanged.

\begin{remark}\label{forced}{\rm
For every QHG transit supported by a $2 \leftrightarrow 3$ move
(top of Figure \ref{CQDfigmove1}), if $E_0$ denotes the new edge
in $T'$ then \begin{equation}\label{2-3comp} W_{\Tt'}(E_0)=1, \
L_{\Tt'}(E_0)=0, \ C_{\Tt'}(E_0)=2.\end{equation} For every QHG
bubble transit, let $f'$ be the unique new $2$-simplex of $T'$
that contains the marked edge $e$. Denote by $E_1$ and $E_2$ the
other edges of $f'$, and by $E_3$ the further new edge of $T'$.
Then we have} \begin{eqnarray}\label{bubblecomp}
 W_{\Tt'}(E_j)=1, \quad  L_{\Tt'}(E_j)=0, \quad j=1,2,3 \\
C_{\Tt'}(E_1)=C_{\Tt'}(E_2)= 0, \quad C_{\Tt'}(E_3)=2.\nonumber
\end{eqnarray}
\end{remark}

\begin{defi}\label{QHGisom}{\rm A {\it QHG isomorphism} between
    QHG triangulated pseudomanifolds is any finite composition of QHG
    transit configurations and oriented simplicial homeomorphisms
that preserve the whole decoration.}
\end{defi}

\subsection{QHG universal functor}\label{QHGFUNC} For every odd
$N\geq 1$, we associate to every QHG triangulated pseudomanifold
$(Z,\Tt)$ a {\it trace tensor} $\Hh_N(\Tt)$, as follows.  Define
an {\it $N$-state} of $T$ as a function that gives every
$2$-simplex an index, with values in $\{0,\dots,N-1\}$. Every
$N$-state determines an entry for each matrix dilogarithm
$\Rr_N(\Delta, b, w, f,c)$. As two tetrahedra induce opposite
orientations on a common $2$-face, an index is down for the
$\Rr_N$ of one tetrahedron while it is up for the other (see
Figure \ref{CQDidealtensor}). By summing over repeated indices we
get the {\it total contraction} of the tensors $\{\Rr_N(\Delta, b,
w, f,c)\}$, that we denote $\textstyle \prod_{\Delta \subset T}
\Rr_N(\Delta,b,w,f,c)$. Let $v_I$ and $v_\delta$ be the number of
vertices of $T\setminus \partial T$ and $\partial T$,
respectively, that correspond to manifold points. We set
\begin{equation}\label{deftt}
\Hh_N(\Tt) = N^{-(v_\delta/2 + v_I)}\prod_{\Delta \subset T}
\Rr_N(\Delta,b,w,f,c)
\end{equation}
The type of a trace tensor $\Hh_N(\Tt)$ depends on the $b$-signs
of the boundary triangles of $(T,b)$. The matrix dilogarithms
themselves are special instances of trace tensors.

\begin{figure}[ht]
\begin{center}
 \includegraphics[width=7cm]{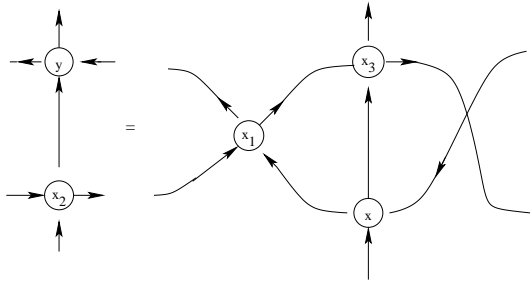}
\caption{\label{CQDSchaeffer} A $2 \leftrightarrow 3$ QHG move
($x_1=y/x$, $x_2=y(1-x)/x(1-y)$, and $x_3=(1-x)/(1-y)$).}
\end{center}
\end{figure}

We have
\begin{prop}\label{fund_func}
For every odd $N\geq 1$, up to sign and multiplication by $N$th
roots of unity the trace tensor $\Hh_N(\Tt)$ is invariant up to
QHG isomorphism.
\end{prop}
This result is a restatement of Theorem 2.1 (2) and Lemma 6.7 of
\cite{BB2}, and summarizes the fundamental functional relations
satisfied by the matrix dilogarithms. In particular those
corresponding to $2 \leftrightarrow 3$ QHG transits are usually called
{\it five terms identities}. In Figure \ref{CQDSchaeffer} we show one
instance in dual terms. The normalization factor $ N^{-v_I}$ in
(\ref{deftt}) is due to the bubble move, that changes by $1$ the
number of internal vertices.
\medskip

The notions of ``boundary'' (also allowing only portions of the
standard boundary - see eg. \cite{T}), bordism, and bordism gluing are
well defined for the category of QHG triangulated pseudomanifolds
considered up to QHG isomorphisms . Hence, for every odd $N\geq 1$,
the association of the trace tensor to each bordism defines a
functor. Note that the identities between trace tensors hold up to the
phase ambiguity of Proposition \ref{fund_func}. Also, the
normalization factor $N^{-v_\delta/2}$ in (\ref{deftt}) compensates
the change in the number of internal vertices when gluing along
complete {\it connected} components of the (usual) boundary (this is
easily adapted to more general gluings).

Any bordism is the gluing of flat/charged $\Ii$-tetrahedra,
considered as elementary bordisms between the couples of $2$-faces
having $b$-sign equal to $-1$ or $+1$. The matrix dilogarithms can
be interpreted as the amplitudes of a diagonal exchange that
relates two quadrilateral triangulations (recall the above discussion
about the geometric quantization derivation of matrix dilogarithms.)

\subsection{Towards geometric specializations}\label{ToGeoSpec} The
category of QHG triangulated pseudomanifolds is built {\it ad hoc}
on the functional properties of the matrix dilogarithms. In order
to eventually get specializations with a geometric content, such
as QHFT, it is necessary to refine more and more our {\it rough}
flat/charged $\Ii$-triangulations. In fact we have to impose some
{\it global constraints} that would be preserved by (possibly
refined) QHG isomorphisms.  It is useful to recall at once some of
these refinements, and related matter. We adopt the notations
introduced above.
\begin{defi}\label{non-rough}{\rm A triple $\Tt=(T,b,w)$ is an {\it
      $\Ii$-triangulation} if at each edge $e$ not contained in
    $\partial T$ we have the {\it edge compatibility relation}
    $W_\Tt(e)=1$.  We say that $\Tt=(T,b,w,f)$ is a {\it flattened}
    $\Ii$-triangulation (and $f$ a {\it global  flattening}) if moreover
    $L_\Tt(e)=0$.}
\end{defi}
It is easily seen that the collection of cross-ratio moduli of an
$\Ii$-triangulation $\Tt$ of a pseudomanifold $Z$ defines a PL {\it
pseudo developing map} $d: \tilde{Z} \to \overline{\mh}^3$, unique up to
post composition with the action of $PSL(2,\C)$, and a {\it holonomy
representation} $h: \pi_1(Z) \to PSL(2,\C)$ such that $d(\gamma(y)) =
h(\gamma)(d(y))$ for every $\gamma \in \pi_1(Z)$, $y\in \tilde{Z}$.

Our main tool for producing $\Ii$-triangulations is the following
  {\it idealization procedure}:
\begin{defi}\label{3dim-id}{\rm
Let $\Tt=(T,b,z)$ be a branched triangulation of a pseudomanifold
$Z$ equipped with a $PSL(2,\C)$-valued $1$-cocycle $z$, where the
cocycle relation is $z(e_0)z(e_1)z(e_2)^{-1} = 1$ on a branched
triangle with edges $e_0$, $e_1$ and $e_2$. A $3$-simplex $\Delta$
of $(T,b,z)$ with vertices $x_0,x_1,x_2,x_3$ is {\it idealizable}
if
\begin{equation}\label{vertid}
u_0=0,\ u_1= z(x_0x_1)(0),\ u_2= z(x_0x_2)(0),\ u_3= z(x_0x_3)(0)
\end{equation}
are four distinct points in $\mc \subset \mc\mathbb{P}^1= \partial
\overline{\mh}^3$. The convex hull is a hyperbolic ideal tetrahedron
with ordered vertices, denoted $(u_0,u_1,u_2,u_3)$. Let us give
each edge $e$ of $\Delta$ the cross-ratio modulus $w(e)\in \C
\setminus \{0,1\}$ of the corresponding edge in
$(u_0,u_1,u_2,u_3)$, and put $w=(w_0,w_1,w_2) \in (\C \setminus
\{0,\ 1 \})^3$. We call $(\Delta,b,w)$ the {\it idealization} of
$(\Delta,b,z)$. If all $3$-simplices of $(T,b,z)$ are idealizable
we say that $\Tt_{\Ii}=(T,b,w)$ is the {\it idealization} of
$\Tt$.}
\end{defi}
Note that the idealization includes the choice of a base point
(ie. $0$ in Definition \ref{3dim-id}; we make this choice once for
ever).  The idealization of each single $3$-simplex strongly
depends on this choice, but taking a different base point is
equivalent to conjugating the cocycle $z$. Not all branched
triangulations support idealizable cocycles (for example an ideal
triangulation of a knot in $S^3$ only supports the trivial
constant cocycle). However, when it exists we have:
\begin{lem}\label{idealization}{\rm (\cite{BB1}, section 2.4)}
The idealization $\Tt_\Ii$ is an $\Ii$-triangulation of $Z$.
\end{lem}
The following fact will be important below (a version of it was
implicitely used in \cite{BB1}; see (3), Lemma 2.12 and \S 3.2 in
that paper):
\begin{lem}\label{canonicalflat} The idealization $\Tt_\Ii$ is
canonically flattened by taking for each $3$-simplex the
log-branches (with the notations of Definition \ref{3dim-id}):
$$\begin{array}{l}{\rm
l}_0:=\log(u_2-u_1)+\log(u_3)-\log(u_2)-\log(u_3-u_1)\\ {\rm
l}_1:=\log(u_2)+\log(u_3-u_1)-\log(u_0)-\log(u_3-u_2)+\sqrt{-1}\pi\\{\rm
l}_2:=\log(u_3-u_2)+\log(u_0)-\log(u_3)-\log(u_2-u_1)-\sqrt{-1}\pi.
\end{array}$$
Moreover, replacing at each edge of $\Tt_\Ii$ the standard $\log$ with
any other log determination still makes it flattened.
\end{lem}
\noindent {\it Proof.} These expressions are just (corrected)
signed sums of the standard logs of the edge vectors in the
cross-ratio moduli $w_0=(u_2-u_1)u_3/u_2(u_3-u_1)$,
$w_1=u_2(u_3-u_1)/u_0(u_3-u_2)$ and
$w_2=-(u_3-u_2)u_0/(u_2-u_1)u_3$. They clearly define triples of
log-branches. The idealizations for two distinct branchings with
the same $b$-orientation are related by an element in
$PSL(2,\mc)$, that is, a conformal transformation of
$\mathbb{C}P^1$. We deduce that triples $({\rm l}_0,{\rm l}_1,{\rm
l}_2)$ are invariant under a change of branching, because the
angle formed by any pair of vectors is preserved (for instance
$u_2-u_1$ and $u_2$, or $u_3$ and $u_3-u_1$). Then, the edge
compatibility relations $L_\Tt(e)=0$ follow easily by developing
the tetrahedra around $e$ with branchings such that $e=e_0$,
similarly as for Lemma \ref{idealization}. For the last claim we
note that any $\log$ correction appears for each $3$-simplex in
two distinct ${\rm l}_j$, with opposite signs.\hfill$\Box$
\medskip

\noindent{\bf Distinguished flat/charged $\Ii$-triangulations.}
Given an arbitrary non empty and properly embedded tangle $L$ in
$Z$, we say that $(T,H)$ and $\Tt=(T,H,b,w,f,c)$ are {\it
distinguished} if $H$ is a subcomplex of the $1$-skeleton of $T$
isotopic to $L$ and passing through all the vertices that are
manifold points, and containing no singular vertices (we say that
$H$ is Hamiltonian).
\begin{defi}\label{DISFLATCHARGE}{\rm A distinguished $\Ii$-triangulation
$\Tt=(T,H,b,w,f,c)$ is {\it flat/charged} if $f$ is a global flattening,
and
\begin{equation}\label{condc}C_{\Tt}(e)=\left
\lbrace\begin{array}{l} 0\\ 2\end{array}\right.\begin{array}{l}
{\rm if}\ e \subset H\\ {\rm if}\ e \subset T\setminus(H \cup
\partial T)\end{array}
\end{equation} In such a case $c$ is said a {\it global charge}.}
\end{defi}
The existence of global flattenings of an arbitrary $\Ii$-triangulation with
empty boundary, and of global charges (such that $C_\Tt(e)=2$ at all
edges $e$) on any topological ideal triangulation of an oriented
$3$-manifold whose boundary consists of tori, was proved by Neumann in
\cite{N}, section 9, and \cite{N0}, section 6, respectively. This last
result is easily adapted to the existence of global charges on
distinguished triangulations $(T,H)$ (see \cite{BB1}, Theorem 4.7).
\smallskip

We have to refine the QHG isomorphisms in order to deal with
distinguished flat/charged  $\ \Ii$-triangulations. First we have to
incorporate the Hamiltonian tangles into the bare moves. Any
positive $2 \to 3$ move $T\to T'$ naturally specializes to a move
$(T,H)\to (T',H')$; in fact $H'=H$ is still Hamiltonian. For
positive bubble moves, we assume that an edge $e$ of $H$ lies in
the boundary of the involved $2$-simplex $f$; then $e$ lies in the
boundary of a unique $2$-simplex $f'$ of $T'$ containing the new
vertex of $T'$. We define the Hamiltonian subcomplex $H'$ of $T'$
just by replacing $e$ with the other two edges of $f'$. The
inverse moves are defined in the same way; in particular, for
negative $3\to 2$ moves we require that the edge disappearing in
$T$ belongs to $T' \setminus H'$. The $2 \leftrightarrow 3$ QHG
transit specializes verbatim. For the bubble transit we just impose
that the marked edge $e$ (see section \ref{QHGPM}) coincides with
the above edge of $H$. Thanks to (\ref{2-3comp}) and
(\ref{bubblecomp}) we see that distinguished flat/charged
$\Ii$-triangulations are closed under such refined QHG transits.

\begin{remark}\label{univfc}{\rm The residue mod$(2\sqrt{-1}\pi N)$ of
the classical log-branches of a flattened $\Ii$-triangulation are
equivalently given by $N$th roots $w'(h)$ of the cross-ratio moduli
$w(h)$ such that $w'(e_1)w'(e_2)w'(e_3)=1$ at each $3$-simplex and
$\textstyle \prod_{h \in \epsilon_T^{-1}(e)} w'(h)^{*_{b}}=1$ at each
edge, with the notations of (\ref{totw}).  For level $N$ quantum
log-branches of distinguished flat/charged $\Ii$-triangulations this
is replaced with} $$w'(e_1)w'(e_2)w'(e_3)=\exp(-*_b\sqrt{-1}\pi/N)$$
{\rm and} $$\prod_{h \in \epsilon_T^{-1}(e)} w'(h)^{*_{b}}=\left
\lbrace\begin{array}{l} 1 \\ \exp(-2\sqrt{-1}\pi/N)
\end{array}\right.\begin{array}{l} {\rm if}\ e \subset H\\
{\rm otherwise.}\end{array}$$
\end{remark}

\noindent{\bf Cohomological weights and structural facts.} Given
an $\Ii$-triangulation $\Tt=(T,b,w)$ without boundary, let $T_0$
be the complement of an open cone neighborhood of each
$0$-simplex. This is a disjoint union of triangulated closed
oriented surfaces, the links of the vertices in $T$. To each
flattening $f$ of $\Tt$ a class $\gamma(f) \in H^1(\partial
T_0;\mc)$ is associated as follows.  Represent any non zero
integral $1$-homology class $a$ of $\partial T_0$ by ``normal
paths", that is, a disjoint union of oriented essential simple
closed curves transverse to the triangulation and such that no
component enters and exits from the same face of a $2$-simplex.
Such a curve selects a vertex for each $2$-simplex. The value
$\gamma(f)(a)$ is defined as the signed sum of the log-branches of
the edges ending at the vertices selected by $a$. For each vertex
$v$ the sign is $*_b$ if the path goes in the direction given by
the orientation of $\partial T_0$ as viewed from $v$, and $-*_b$
otherwise. Using the edge compatibility relations of Definition
\ref{non-rough} it is easily checked that $\gamma(f)$ does not
depend on the choice of normal path representative. If the
holonomy of $\Tt$ is trivial or parabolic about each (non
manifold) point this class takes values in $2\mz$ (\cite{N},
Proposition 5.2). Similarly we can define $\gamma_2(f) \in
H^1(T_0;\mz/2\mz)$ by using normal paths in $T_0$ and taking
modulo $2$ sum of the flattenings we meet along the paths. We call
$(\gamma(f),\gamma_2(f))$ the {\it (cohomological) weight} of $f$.
These definitions extend immediately to global charge, replacing
log-branches with charges; then branchings are not needed to
compute the signed sums (ie. put $*_b=1$).

In \cite{N0}, section 4 (see also \cite{N}, section 9), Neumann
defines an integral chain complex $\mathcal{J}$ such that global
charges are defined from cycles at level $3$. The maps $\gamma$
and $\gamma_2$ above are well-defined on the third homology group
$H_3(\mathcal{J})$ and satisfy:
\begin{teo}\label{theoweight}{\rm (\cite{N0}, Theorem 5.1)} The sequence
$$0\rightarrow H_3(\mathcal{J})
\stackrel{(\gamma_2,\gamma)}{\longrightarrow} H^1(T_0;\mz/2\mz)
\oplus H^1(\partial T_0;\mz) \stackrel{i^*-r}{\longrightarrow}
H^1(\partial T_0;\mz/2\mz) \rightarrow 0$$ is exact, where $r:
H^1(\partial T_0;\mz) \rightarrow H^1(\partial T_0;\mz/2\mz)$ is
the coefficient map and $i^*: H^1(T_0;\mz/2\mz) \rightarrow
H^1(\partial T_0;\mz/2\mz)$ is induced by the inclusion $\partial
T_0 \rightarrow T_0$.
\end{teo}
In particular, any pair $(h,k) \in H^1(T_0;\mz/2\mz) \times
H^1(\partial T_0;\mz)$ with $r(k)=i^*(h)$ is a weight for some global
charges. This fact extends to log-branches as follows. Assume that the
holonomy of $\Tt$ restricted to each component of $\partial T_0$ takes
values in a Borel subgroup of $PSL(2,\mc)$, so that it fixes some
point in $\partial \overline{\mathbb{H}}^3$. Let $\gamma'$ be defined as
above, except that for each $3$-simplex log-branches are taken with
flattenings $f_0=f_1=0$. It can be checked that $\gamma'(a)$ is the
logarithm of the derivative of the holonomy of $a$ (a similarity), up
to multiples of $2\sqrt{-1}\pi$. Then, any pair $(h,k) \in
H^1(T_0;\mz/2\mz) \times H^1(\partial T_0;\mc)$ such that
\begin{equation}\label{weightcondition}\left\lbrace
\begin{array}{l}(k-\gamma')/\sqrt{-1}\pi \in H^1(\partial T_0;\mz)\\
r((k-\gamma')/\sqrt{-1}\pi)=i^*(h)\end{array}\right.
\end{equation}
 is a weight for some flattening. For instance, when
$\gamma'(a)\in 2\sqrt{-1}\pi\mz$ for all $a$ the first condition
means that $k$ is integral.

Finally, the structure of the spaces of flattenings and integral
charges is given by Theorem 2.4 in \cite{N0}. For each fixed
weight $(h,k)$ they form an affine space over an integral lattice.
Generators have the following combinatorial realization: for each
$3$-simplex in the star of an edge $e$, add $+1$ to the
flat/charges of one of the two other pairs of opposite edges, and
$-1$ for the other pair, so that the total log-branches or charges
stay equal everywhere. In particular, for flattenings of an
idealization $\Tt_\Ii$ any generator is obtained by adding $+1$ to
the log determination at some edge. Hence any flattening of
$\Tt_\Ii$ inducing the weight of the canonical flattening of Lemma
\ref{canonicalflat} differs from it as described in the statement.

We note that the above refined QHG isomorphisms preserve the weights
(see eg. Lemma 4.12 in \cite{BB1}). The difference in considering
global flattenings or charges with different mod($2$) weights in
$H^1(T_0;\mz/2\mz)$ seems to carry not so essential information (see
Theorem \ref{main} and Remark \ref{scissors}). This contrasts with
boundary weights in $H^1(\partial T_0;\mz)$, which play a key role in
surgery formulae (see \cite{N}, theorems 14.5 and 14.7, and section
\ref{DEHNF} below). We note that a process involving $2$-handle
surgery allows to define explicit isomorphisms between lattices of
flat/charges with different boundary weights (see \cite{N}, section
11, p. 457).

\section{ Parameters for $PSL(2,\C)$-characters of surfaces}
\label{PARAM} Fix a compact closed oriented surface $S$ of genus
$g$ with a non empty set $V= \{v_1,\dots \ , v_r\}$ of marked
points, and negative Euler characteristic $\chi(S\setminus V)<0$.
Denote by $\pi$ the fundamental group of $S\setminus V$, and by
$$\Rr(g,r) = {\rm Hom}(\pi, PSL(2,\C))/PSL(2,\C)$$ the set of {\it
all} conjugacy classes of $PSL(2,\C)$-valued representations of
$\pi$. The group $\pi$ is free of rank $\kappa = 2g + r -1$. Any
choice of free generators of $\pi$ identifies the set ${\rm
Hom}(\pi, PSL(2,\C))$ with $PSL(2,\C)^\kappa$. Different such
identifications are related by algebraic automorphisms of
$PSL(2,\C)^\kappa$.  Moreover, the isomorphism $PSL(2,\C)\cong
SO(3,\C)$ induced by the adjoint action ${\rm Ad}: PSL(2,\mc)
\rightarrow {\rm Aut}(sl(2,\mc))$ implies that ${\rm Hom}(\pi,
PSL(2,\C))$ is an {\it affine} complex algebraic set, with the
complex algebraic action of $PSL(2,\C)$ by conjugation.

As the rough topological quotient space $\Rr(g,r)$ is not
Hausdorff, it is more convenient to consider the algebro-geometric
quotient $$X(\pi)={\rm
  Hom}(\pi, PSL(2,\C))//PSL(2,\C)$$ of invariant theory, called the
{\it variety of $PSL(2,\C)$-characters} (see eg. \cite{HP}).
Recall that $X(\pi)$ is a complex affine algebraic set such that
there exists a surjective regular map $$t:{\rm Hom}(\pi,
PSL(2,\C))\to X(\pi)$$ inducing an isomorphism $t^*$ between the
regular functions on $X(\pi)$ and the regular functions on ${\rm
Hom}(\pi, PSL(2,\C))$ invariant by conjugation. In general
$t(\gamma)=t(\sigma)$ does not imply that $\gamma$ and $\sigma$
are conjugate, but this is true if we restrict to {\it
irreducible} representations: we have ${\rm
  Hom}^{irr}(\pi, PSL(2,\C)) = t^{-1}(X^{irr}(\pi))$, where
$X^{irr}(\pi)= t({\rm Hom}^{irr}(\pi, PSL(2,\C))$, so that the
(restricted) rough quotient $\Rr(g,r)^{irr}$ and the algebraic
quotient $X^{irr}(\pi)$ coincide.

 We can deal with the whole of $\Rr(g,r)$ to construct
the QHFT. Hence below we treat the complex dimension of subsets of
$\Rr(g,r)$ somewhat formally, as everything can be substantiated
in terms of $X(\pi)$, or by restriction to $X^{irr}(\pi)$. For
instance, as $PSL(2,\C)$ has trivial centre and complex dimension
$3$ we say that the complex dimension of $\Rr(g,r)$ is $3\kappa -
3 = -3\chi(F)$.
\subsection {Efficient triangulations}\label{e-triang}
Fix a surface $F$ with $r$ boundary components, obtained by
removing from $S$ the interior of small $2$-disks $D_i$ such that
$v_i \in \partial D_i$.

\begin{figure}[ht]
\begin{center}
\includegraphics[width=6cm]{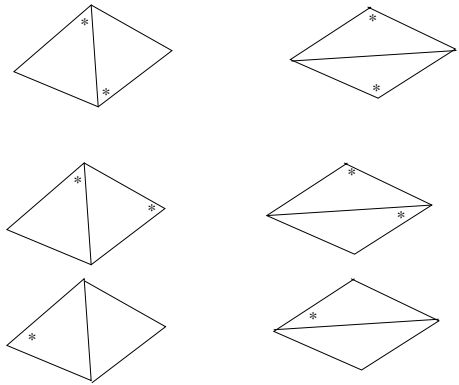}
\caption{\label{eTflip} The flips with marked corners.}
\end{center}
\end{figure}

A triangulation $T'$ of $S$ with the set of vertices equal to $V$ is
called a {\it topological ideal triangulation} of $S\setminus
V$. Given such a $T'$, we need a marking of corners of the
$2$-simplices. The best suited to $3$-dimensional extension are
induced by global branchings $b'$ of $T'$, and it is known that pairs
$(T',b')$ always exist. Then, as in section \ref{QHG}, we have a sign
function $\sigma = \sigma_{(T',b')}$. A {\it corner map} $v\mapsto
c_v$ associates to each vertex $v$ of $T'$ the corner at $v$ of a
triangle, say $t_v$, in its star. We say that $v\mapsto c_v$ is {\it
$t$-injective} if $v\mapsto t_v$ is injective.
\begin{lem}\label{injcorner}
For every $(g,r)\ne (0,3)$, every triangulation $T'$ of $S$ with $r$
vertices admits $t$-injective corner maps.
\end{lem}
\noindent {\it Proof.} For $(g,r) =(0,4)$ or $g>0$ and $r=1$, it
is immediate to construct such a triangulation. Subdividing a
triangle by taking the cone from an interior point preserves the
existence of $t$-injective corner maps. By induction on $r$ we
deduce that for every $(g,r)\ne (0,3)$ there exist triangulations
of $S$ as in the statement.\hfill\break In Figure \ref{eTflip} the
corner selection is specified by a $*$, and the rows show all the
possible flip moves on triangulations of $S$, up to obvious
symmetries, that preserve the injectivity of corner maps. Consider
triangulations $T'$, $T''$ of $S$ with $r$ vertices, such that
$T''$ supports a $t$-injective corner map. It is well known that
$T''$ is connected to $T'$ via a finite sequence of flips. The
$t$-injective corner map for $T''$ yields a $t$-injective corner
map for $T'$ by decorating these flips as in Figure
\ref{eTflip}.\hfill$\Box$\medskip

There are no obstructions to use arbitrary corners maps in what
follows. Specializing to injective ones just simplifies the
exposition. For $(g,r)= (0,3)$, triangulations have two
$2$-simplices, one with two selected corners.  From now on, we
assume that $(g,r) \neq (0,3)$, the extension to the $(0,3)$ case
being straightforward.

\smallskip

\begin{figure}[ht]
\begin{center}
\includegraphics[width=6cm]{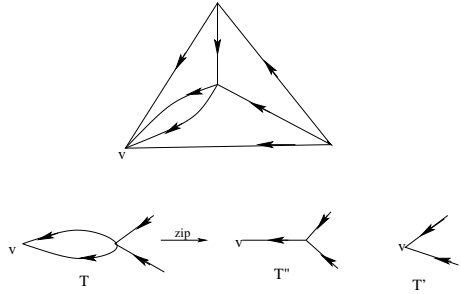}
\caption{\label{newb} The branched triangulated cell $s_v$.}
\end{center}
\end{figure}

Given a pair $(g,r) \ne (0,3)$ and $(T',b')$ as in Lemma
\ref{injcorner}, fix a $t$-injective corner map $v\mapsto c_v$. In the
interior of each triangle $t_v$ consider a bigon $D_v$ with one vertex
at $v$, and call $v'$ the other vertex. Remove from $t_v$ the interior
of $D_v$, and triangulate the resulting cell $s_v = t_v \setminus {\rm
  Int}(D_v)$ by taking the cone with base $v'$. Repeating this
procedure for each $t_v$, we get a triangulation $T$ of $F$ with $2r$
vertices and $p+2r$ triangles, where $p$ denotes the number of
triangles of $T'$.  The set of edges of $T$, $E(T)$, contains $E(T')$
in a natural way, and $|E(T)|= |E(T')|+4r$. We extend $b'$ to a
branching $b$ on $T$ as in in Figure \ref{newb}, and the sign function
$\sigma_{(T',b')}$ to the triangles of $(T,b)$ in the natural way.
Note that the figure shows only one of the possible branching
configuration.  In general we extend $b'$ to $b$ so that we can
recover $(T',b')$ from $(T,b)$ by ``zipping'' and ``collapsing'', as
suggested at the bottom of Figure \ref{newb}. In what follows, for
simplicity we will refer to this configuration, as the treatment of
the others is similar.
\begin{defi}\label{defet}
{\rm We call the pair $(T,b)$ an efficient triangulation (for short:
    $e$-{\it triangulation}) of $F$. For each vertex $v$, the {\it
    preferred innner triangle at} $v$ is the triangle $\tau_v$ in
    $s_v$ with an edge on $\partial F$ whose $b$-orientation coincides
    with the boundary orientation of $F$.}
\end{defi}

\subsection{Cocycle parameters}\label{Dpar}
The inclusions of Int$(F)$ into $F$ and $S\setminus V$ induce
identifications of the respective fundamental groups.

Fix an $e$-triangulation $(T,b)$ of $F$. Denote by $Z(T,b)$ the space
of $PSL(2,\C)$-valued $1$-cocycles on $(T,b)$; we stipulate
that on a triangle with ordered $b$-oriented edges $e_0,e_1,e_2$ the
cocycle relation is $z(e_0)z(e_1)z(e_2)^{-1} = 1$. In this section we
construct a parametrization of $\Rr(g,r)$ based on cocycle
coefficients by specifying subsets of $Z(T,b)$ with small `residual
gauge groups', that make principal algebraic bundles over the
``strata'' of a suitable partition of $\Rr(g,r)$. These strata are
determined by the holonomies around the boundary components of
$F$.

Write $C(T,b)$ for the space of $PSL(2,\C)$-valued $0$-cochains on
$(T,b)$, that is the $PSL(2,\C)$-valued functions defined on the
set of vertices of $T$.  Two $1$-cocycles $z$ and $z'$ are said
equivalent {\it up to gauge transformation} if there is a
$0$-cochain $\lambda$ such that, for every oriented edge
$e=[x_0,x_1]$, we have $$z'(e) = \lambda(x_0)^{-1}z(e)\lambda(x_1)
\ . $$ It is well known that the quotient set $H(T,b) =
Z(T,b)/C(T,b)$ is in one-one correspondence with $\Rr(g,r)$.
Indeed, for any fixed $x_0 \in T$ we have a natural surjective map
$$f_{x_0}: Z(T,b) \to {\rm Hom}(\pi, PSL(2,\C)),$$ where
$\pi=\pi_1(F,x_0)$, and two representations $f_{x_0}(z)$ and
$f_{x_0}(z')$ define the same point in $\Rr(g,r)$ if and only if
$z$ and $z'$ are equivalent up to gauge transformation. Note that
the complex dimension of $H(T,b)$ is $3(|E(T)| - (p+2r) - 2r) =
-3\chi(F)$, which is the dimension of $\Rr(g,r)$. This is because
the complex dimension of $C(T,b)$ is equal to $6r$, the set
$Z(T,b)$ is defined by $3(p+2r)$ polynomial relations on $3|E(T)|$
variables, and $PSL(2,\C)$ has trivial centre. (Recall that $p$ is
the number of triangles of the initial triangulation $T'$ of $S$
with $r$ vertices.)

Denote by $B^+(2,\C)$ (respectively $B^-(2,\C)$) the Borel subgroup of
$SL(2,\C)$ of upper (respectively lower) triangular matrices. Let
$PB^\pm(2,\C)= B^\pm(2,\C)/\pm I$ and put
$$ P=\left( \begin{array}{cc} 0 & 1 \\ 1 & 0  \end{array}\right).$$
Define a map
\begin{equation}\label{convert}
\Psi : SL(2,\C)\to SL(2,\C),\quad \Psi(A)=PAP
\end{equation}
and denote by $\Psi$ the induced automorphism of $PSL(2,\C)$. We
have $\Psi(B^\pm(2,\C)) = B^\mp(2,\C)$. 
For any $g\in PSL(2,\C)$ we distinguish the {\it type} of $g$ as
{\it
  trivial}, {\it parabolic} or {\it generic} (the latter for elliptic or
loxodromic) with the obvious meaning. We denote
$$\mathcal{C}^+(g) \in PB^+(2,\C)$$ the canonical
upper triangular matrix (up to sign) representative of the conjugacy
class of $g$ (for generic $g$ we normalize $\mathcal{C}^+(g)$ by stipulating
that the top diagonal entry has absolute value $>1$), and we set
$$\mathcal{C}^-(g) = \Psi(\mathcal{C}^+(g)). $$
We define the type of $\rho \in \Rr(g,r)$ as the $n$-uple of types of
the $\rho$-holonomies of the oriented boundary components
$\gamma_1,\ldots,\gamma_r$ of $F$. Put:
$$\begin{array}{l} \Rr(g,r,t) = \{ \rho \in \Rr(g,r) \mid \rho {\rm\
has \ type\ } t\} \\ \Rr(g,r,\mathcal{C}^\pm) = \{ \rho \in \Rr(g,r,t)
\mid \mathcal{C}^\pm(\rho(\gamma_i)) = \mathcal{C}_i^\pm\}
\end{array}$$
where $\mathcal{C}^\pm=(\mathcal{C}_1^\pm,\ldots,\mathcal{C}_r^\pm)$,
$\mathcal{C}_i^\pm$ is an arbitrary diagonal or unipotent element in
$PB^\pm(2,\C)$, and $\mathcal{C}^\pm$ has type $t$.
\smallskip
For any $g \in PB^\pm(2,\C)$ let us write $g=[a,b]^\pm$, where $a$ is
the top diagonal entry of $g$ and $b$ is the non diagonal one (we do
the abuse of confusing $g$ with its projective class). For every
vertex $v_i$ of $(T,b)$, let $e_1^i$ be the boundary edge of the
preferred inner triangle $\tau_{v_i}$ (see Definition \ref{defet}),
and $\gamma_i$ the oriented boundary loop of $F$ based at
$v_i$. Recall the projections $f_{v_i}: Z(T,b) \to {\rm Hom}(\pi,
PSL(2,\C))$.
\begin{defi} \label{defcoc}
The set of ($\pm$)-{\it cocycle parameters} for
$\Rr(g,r,\mathcal{C}^\pm)$ is $$\Zz(T,b,\mathcal{C}^\pm) = \{ z
\in Z(T,b) \mid \forall i=1,\ldots,r,\ f_{v_i}(z)(\gamma_i) =
\mathcal{C}^\pm_i,\ z(e_1^i)= [1,1/2]^\pm \}.$$
\end{defi}
Clearly $\Zz(T,b,\mathcal{C}^\pm)$ is non empty. In fact, any
$z\in Z(T,b)$ with ${\rm conj} \circ f(z) \in
\Rr(g,r,\mathcal{C}^\pm)$ is equivalent to one in
$\Zz(T,b,\mathcal{C}^\pm)$ via some gauge transformation.
Moreover, given a vertex $v_i$ and a $0$-cochain $s$ with support
at $v_i$ and $v_i'$, $s$ maps $\Zz(T,b,\mathcal{C}^\pm)$ onto
itself if and only if $$s(v_i)\in {\rm Stab}(\mathcal{C}^\pm_i)$$
and $$s(v_i')^{-1}[1,1/2]^\pm s(v_i) = [1,1/2]^\pm. $$ As the
triangles $\tau_{v_i}$ are in one-one correspondence with the
$v_i$, we deduce that there is a projection
\begin{equation}\label{projc}
p^\pm_{\mathcal{C}}: \Zz(T,b,\mathcal{C}^\pm)\to
\Rr(g,r,\mathcal{C}^\pm) \end{equation} with fiber isomorphic to
the group $$\Gg (T,b,\mathcal{C}^\pm):={\rm
Stab}(\mathcal{C}^\pm_1)\times \dots \times {\rm
Stab}(\mathcal{C}^\pm_r)$$ that we call the {\it residual gauge
transformations}.  Denote by $\Zz(T,b,t)^\pm$ (respectively
$p^\pm_t$) the union of the $\Zz(T,b,\mathcal{C}^\pm)$
(respectively $p^\pm _{\mathcal{C}}$) over all $\mathcal{C}^\pm$
with type $t$, and put $$\textstyle \Zz(T,b)^\pm = \coprod_t
\Zz(T,b,t)^\pm, \quad p^\pm =  \coprod_t p^\pm_t).$$ For futur
reference, we summarize the above constructions in the following
proposition.
\begin{prop}\label{topcomp} For every $\Cc^\pm$, the projection
  $p^\pm : \Zz(T,b)^\pm \rightarrow \Rr(g,r)$ restricts
to a complex affine algebraic principal $\Gg
(T,b,\mathcal{C}^\pm)$-bundle $$\Zz(T,b,\mathcal{C}^\pm) \to
\Rr(g,r,\mathcal{C}).$$ The map $\Psi$ in (\ref{convert}) yields
an isomorphism $\Zz(T,b,\mathcal{C}^\pm) \cong
\Zz(T,b,\mathcal{C}^\mp)$. Moreover we have: $$\begin{array}{lll}
{\rm dim}(\Zz(T,b,\mathcal{C}^\pm)= -3\chi(F)\\ {\rm
dim}(\Rr(g,r,\mathcal{C})) = -3\chi(F) - {\rm dim}(\Gg
(T,b,t_{\mathcal{C}}))\\ {\rm dim}(\Rr(g,r,t_{\mathcal{C}})) =
{\rm dim}(\Rr(g,r,\mathcal{C})) + \alpha(t_{\mathcal{C}})
\end{array}$$
where $t_{\mathcal{C}}$ is the type of $\mathcal{C}$ and
$\alpha(t_{\mathcal{C}})$ is the number of
generic entries of $t_{\mathcal{C}}$.
\end{prop}
Observe that: if $t_{{\rm gen}}$ is the purely generic type, then
${\rm dim}(\Rr(g,r,t_{{\rm gen}}))= -3\chi(F)$; if $t_{{\rm par}}$ is
purely parabolic, then ${\rm dim}(\Rr(g,r,t_{{\rm par}}))= -3\chi(F)-
r$; if $t_I$ is purely trivial, then ${\rm dim}(\Rr(g,r,t_I))=
-3\chi(F) -3r = -3\chi(S)= 6g-6$. Hence $\Rr(g,r,t_{{\rm gen}})$ is a
dense open subset of $\Rr(g,r)$. Moreover, via the inclusion of
closures, we have a filtration of $\Rr(g,r)$ for which $\Rr(g,r,t_I)$
is the `deepest' part. It would be interesting to study the
singularities of the closure of each $\Rr(g,r,t)$ in $\Rr(g,r)$, in
order to check if this filtration induces a ``stratification'' of
$\Rr(g,r)$.

\begin{remark}\label{real} {\rm For parabolic elements $g \in
    PSL(2,\R)$ we have two conjugacy classes, that can be
distinguished by a sign. So, replacing $PSL(2,\C)$ with
$PSL(2,\R)$, the constructions of this section still work by
associating a {\it sign} to each parabolic end of the surface $S
\setminus V$. Note that for $PSL(2,\R)$-valued cocycles the
idealization procedure described in Section \ref{Ipar} below gives
only degenerate triangles and tetrahedra (with real shapes or
cross-ratio moduli).}
\end{remark}

\noindent {\bf Example: the Fricke space.} Suppose that $S$ has a
unique marked point $v$. Choose a standard curve system
$\mathcal{S}=\{a_i,b_i\}_{i=1}^{g}$ based at $v$, so that
$$\pi_1(S,v) =<a_1,b_1,\ldots,a_g,b_g\mid [a_1,b_1]\ldots
[a_g,b_g] =1>.$$ Cutting open $S$ along $\mathcal{S}$ we get a
$4g$-gon $P$ with oriented boundary edges. Taking the cone to a
vertex it is easy to construct a branched triangulation of $P$,
which induces one, say $(T',b')$, for $(S,v)$.  Denote by $(T,b)$
any $e$-triangulation obtained from $(T',b')$. \hfill\break Recall
that the Teichm\"uller space $\mathcal{T}(S)$ can be identified
with the set of conjugacy classes of $PSL(2,\mr)$-valued discrete
faithful representations of $\pi_1(S)$, and that we have the
well-known (real-analytic) Fricke parametrization $\mathcal{T}(S)
\cong \mr^{6g-6}$ (see e.g. \cite{Ab}). For each $z \in
\mathcal{T}(S)$ the Fricke coordinates of $z$ are matrix entries
of the $\hat{z}(\gamma)$ for all $\gamma \in \mathcal{S}$, where
$\hat{z}$ is a representative of $z$ specified by fixing once and
for all three of the fixed points of $z(a_g)$ and $z(b_g)$. So
$\mathcal{T}(S)$ embeds in the space of real cocycle parameters
for $F$ with trivial type $t_I$. This embedding is generalized
easily to the Teichm\"uller space of arbitrary bordered Riemann
surfaces, by considering the spaces $\Zz(T,b,t)$ for all types $t$
and $e$-triangulations of $(S,V)$ with arbitrary $V$ (see Remark
\ref{real}).

\subsection {Cross-ratio parameters}\label{Ipar}
In this section we derive from the cocycle parameters
$\Zz(T,b)^\pm$ other parameters for $\Rr(g,r)$ which are related
to the shear-bend coordinates for pleated hyperbolic surfaces.
These parameters are obtained via an {\it idealization procedure}
that includes the choice of a base point on the Riemann sphere. We
fix this base point as $0$. Note that the $B^+(2,\C)$-orbit of $0$
is the whole of $\C$, while $B^-(2,\C)$ fixes $0$. Hence the
symmetry between $\Zz(T,b)^+$ and $\Zz(T,b)^-$ given by the map
$\Psi$ in (\ref{convert}) shall be broken.
\begin{defi} \label{triideal-pre} {\rm
(Compare with Definition \ref{3dim-id}.) Let $(K,b)$ be any oriented
surface branched triangulation. Let $z$ be any $PSL(2,\C)$-valued
cocycle on $(K,b)$.  We say that $z$ is {\it idealizable} if for any
triangle $t$ of $(T,b)$ with $b$-ordered edges $e_0$, $e_1$ and $e_2$,
the points $u_0=0$, $u_1= z(e_0)(0)$ and $u_2= z(e_2)(0)$ are distinct
in $\mc$. We say that the complex triangle with vertices $u_0$, $u_1$
and $u_2$ is the \emph{idealization} of $t$.}
\end{defi}
\begin{defi} \label{triideal}
{\rm Let $(T,b)$ be an $e$-triangulation of $F$ obtained from a
branched ideal triangulation $(T',b')$ of $S \setminus V$. A
cocycle $z\in \Zz(T,b)^+$ is {\it (strongly) idealizable} if:
\smallskip

(a) $z$ is idealizable;
\smallskip

(b) $\Psi(z) \in \Zz(T,b)^-$, which {\it cannot} be idealizable at the
triangles having an edge on $\partial F$, is nevertheless idealizable
at every other triangle.
\smallskip

We denote by $\Zz_I(T,b)^+$ the set of (strongly) idealizable
cocycles, and we put $\Zz_I(T,b)^- = \Psi(\Zz_I(T,b)^+)$ and
$\Rr_I(T,b) = p^+(\Zz_I(T,b)^+) = p^-(\Zz_I(T,b)^- )$.}
\end{defi}

Clearly, $\Zz_I(T,b)^+$ is a non-empty dense open subset of
$\Zz(T,b)^+$. If every edge of $T$ has distinct endpoints (in case
we say that $T$ {\it quasi-regular}), then $\Rr_I(T,b) =
\Rr(g,r)$. In general, characters of representations with a free
action on a non-empty domain of $\mc\mathbb{P}^1$ (such as
quasi-Fuchsian representations) always belong to $\Rr_I(T,b)$, for
any $(T,b)$. By using the arguments of \cite{Kmod}, Theorem 1, it
can be shown that for any character $\rho\in \Rr(g,r)$ of
irreducible representations, there exists an $e$-triangulation
$(T,b)$ with idealizable cocycles representing $\rho$. In fact,
the union of a finite number of spaces $\Zz_I(T,b)^+$ cover the
whole of $\Rr(g,r)^{irr}$.

\smallskip

\noindent{\bf Exponential $\Ii$-parameters.}  For any cocycle $z
\in \Zz_I(T,b)^+$, we associate a non zero complex weight
$W^+(z)(e)$ to each edge $e$ of $T$ that is not contained in
$\partial F$, as follows. Let $p_e$ be the initial endpoint of
$e$, and $t_l$ and $t_r$ the left and right adjacent triangles (as
viewed from $e$). Locally modify the branching on $t_l \cup t_r$
by cyclically reordering the vertices on each triangle, so that
$p_e$ is eventually the source of the new branching on both $t_l$
and $t_r$. The {\it ($+$)-exponential $\Ii$-parameter} $W^+(z)(e)$
is the cross-ratio modulus at $e$ of the (possibly degenerate)
branched oriented hyperbolic ideal tetrahedron spanned by the
idealization of $t_r \cup t_l$, where the branching completes the
one of $t_l \cup t_r$ so that $*_b=1$. 
\smallskip

Let us assume now that $\Psi(z)\in \Zz_I(T,b)^-$. For each $v \in V$,
denote by $e_v^1$ and $e_v^2$ the edges of the triangle $t_v$ of $T'$
having $c_v$ as corner. At an edge $e$ of $T'$ distinct from any of
the $e_v^i$, we define the {\it ($-$)-exponential $\Ii$-parameter}
$W^-(\Psi(z))(e)$ in the same way as $W^+(z)(e)$, but taking
the idealization for $\Psi(z)$ instead of $z$. If $e$ is one of the
$e_v^i$, the formula works as well, except that each left/right
triangle with an edge on $\partial F$ is replaced with the innermost
triangle in the corresponding triangle $t_v$ of $T'$. Note that $t_r
\cup t_l$ is again a quadrilateral, because $0$ is a fixed point of
the cocycle values at the boundary edges of $F$.
\smallskip

By varying the edge in $T$ or $T'$ for every $\Cc^\pm$, we have two
{\it ($\pm$)-exponential $\Ii$ parameter maps}
(recall that the number of edges of $T'$ is $ -3\chi(F)$, and $p$
denotes the number of triangles of $T'$):
$$W^+: \Zz_I(T,b,\mathcal{C}^+)
    \rightarrow (\mc \setminus \{ 0\})^{-3\chi(F)+2p}$$
$$W^-: \Zz_I(T,b,\mathcal{C}^-) \rightarrow (\mc \setminus \{ 0
    \})^{-3\chi(F)}.$$

    \begin{defi}\label{phasespace} {\rm We call
        $W^\pm(T,b,\mathcal{C}) = W^\pm(\Zz_I(T,b,\mathcal{C}^\pm))$
        the {\it ($\pm$)-exponential $\Ii$-parameter space} of
        $\Rr_I(T,b,\mathcal{C})$. }
\end{defi}
We stress again that the result of the idealization strongly
depends on the choice of the base point, here $0$. In particular,
the exponential $\Ii$-parameters are {\it not} invariant under
arbitrary gauge transformations of cocycles.  A remarkable
exception is for gauge transformations associated to $0$-cochains
$\lambda$ with values in $PB^-(2,\C)$. Indeed, these act on the
idealization of quadrilaterals as conformal transformations of the
{\it four} vertices (this is because $0$ is fixed by every
$\lambda(v)$, $v \in V$), and cross-ratios are conformal
invariants. This makes a big difference between $W^+$ and $W^-$.
In fact the whole of the residual gauge group
$\Gg(T,b,\mathcal{C}^+)$ acts on $W^+(T,b,\mathcal{C})$ via the
map $W^+$. On the other hand, consider the subgroup of
$\Gg(T,b,\mathcal{C}^-)$ defined as $$B\Gg(T,b,\mathcal{C}^-):=
\prod_{i=1,\dots, r}{\rm Stab}(\mathcal{C}^-_i)\cap PB^-(2,\C) \
.$$ Any fiber $(W^-)^{-1}(W^-(z))$ is given by the
$B\Gg(T,b,\mathcal{C}^-)$-orbit of $z$ in
$\Zz_I(T,b,\mathcal{C}^-)$. Hence the actual residual gauge
transformations of $W^-(T,b,\mathcal{C})$ are in one-one
correspondence with the quotient set $\Gg(T,b,\mathcal{C}^-)/
B\Gg(T,b,\mathcal{C}^-)$, where the equivalence class of $\lambda$
is $B\Gg(T,b,\mathcal{C}^-)\lambda$. For the matter of notational
convenience, formally put $\textstyle B\Gg(T,b,\mathcal{C}^+):=
\prod_{i=1,\dots, r} {\rm Id}$. The map $$\Theta^\pm:
W^\pm(T,b,\mathcal{C}) \to \Rr(g,r,\mathcal{C}^\pm)$$ given by $
\Theta^\pm(W^\pm(z))= p^\pm_{\mathcal{C}}(z)$, with
$p^\pm_{\mathcal{C}}$ defined in (\ref{projc}), is a principal
$\Gg (T,b,\mathcal{C}^\pm)/ B\Gg(T,b,\mathcal{C}^\pm)$-bundle. The
situation is particularly clean when $\mathcal{C}$ has no trivial
entries :
\begin{prop}\label{no-gauge} If $\mathcal{C}$ has no trivial entries
the ($-$)-exponential $\Ii$-parameter map $W^-:
  \Zz_I(T,b,\mathcal{C}^-) \rightarrow W^-(T,b,\mathcal{C})$ is
  invariant under gauge transformations.
\end{prop}

\noindent{\bf Geometry and computation of $\Theta^\pm$.} The maps
$\Theta^\pm$ can be defined directly in terms of $\Ii$-parameters. A
way to see this is to consider $\Ii$ triangulations of cylinders $C =
F \times [0,1]$, as discussed in Section \ref{FIBERED}, and the
corresponding pseudo-developing maps. We give here another
description. Let us just consider $\Theta^-$. Recall that $PSL(2,\mc)$
is isomorphic to Isom$^+(\mh^3)$, with the natural conformal action on
$\mc\mathbb{P}^1=\partial\bar{\mh}^3$ via linear fractional
transformations. For any $\rho\in \Rr(g,r,\mathcal{C}^-)$, take a
representative $\tilde{\rho}: \pi_1(F,q) \rightarrow PSL(2,\mc)$ in
the conjugacy class. Consider the associated flat principal
$PSL(2,\mc)$-bundle $F_{\tilde{\rho}}$. A trivializing atlas defines a
cocycle $z \in \Zz(T,b,\mathcal{C}^-)$. If $\rho\in
\Rr_I(T,b,\mathcal{C}^-)$, we can take $z \in
\Zz_I(T,b,\mathcal{C}^-)$, so that the trivializing atlas of
$F_{\tilde{\rho}}$ associated to the cellulation $T'^*$ {\it dual} to
$T'$, with edges oriented by using the orientation of $F$ and the
branching orientation of the edges of $T'$, has non trivial transition
functions. These can be viewed as transition functions for the fiber
bundle associated to $F_{\tilde{\rho}}$ and with fibre $\mh^3$. For
each $2$-simplex $t$ of $T'$, there is a unique $g_t \in {\rm
  Isom}(\mh^3)$ (possibly reversing the orientation) mapping the
vertices $u_0$, $u_1$ and $u_2$ of the idealization of $t$ to $0$,
$\infty$ and $-1$ respectively. Then, the transition function along
the edge of $T'^*$ {\it positively transverse} to a given edge $e$ of
$T'$ is of the form $(g_{t_l})^{-1}\circ \varphi(z)(e)\circ g_{t_l}
\in {\rm Isom}^+(\mh^3)$, where $t_l$ is the triangle on the left of
$e$, and $\varphi(z)(e)$ is the isometry of $\mh^3$ of hyperbolic type
fixing $0$ and $\infty$ and mapping $1$ to $W^-(z)(e)$. Analytic
continuation defines the parallel transport of $F_{\tilde{\rho}}$
along paths transverse to $T'$, whence a representation into
$PSL(2,\mc)$ of the groupoid of such paths, well-defined up to
homotopy rel($\partial$)). In particular, it gives a practical recipe
to compute $\tilde{\rho}$, that we describe now.
\begin{figure}[ht]
\begin{center}
\includegraphics[width=10cm]{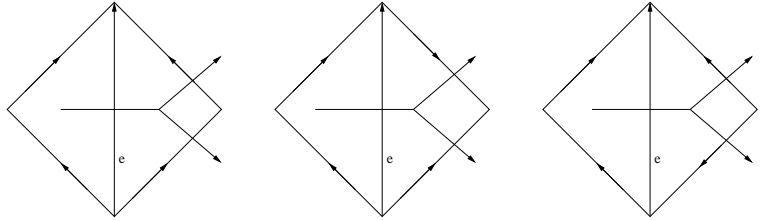}
\caption{\label{matrixrep} The recipe for reading off holonomies
from exponential $\Ii$-parameters.}
\end{center}
\end{figure}

Let the base point $q$ be not in the $1$-skeleton of $T'$.  Given
an element of $\pi_1(F,q)$, represent it by a closed curve
$\gamma$ in $F$ transverse to $T'$, and which do not departs from
an edge it just entered. Assume that $\gamma$ intersects an edge
$e$ of $T'$ positively with respect to the orientation of $F$.
Figure \ref{matrixrep} shows three possible branching
configurations for the two triangles glued along $e$. Fix
arbitrarily a square root $W^-(z)(e)^{1/2}$ of $W^-(z)(e)$.
Consider the elements of $PSL(2,\mc)$ given by $$\gamma(e)=\left(
\begin{array}{cc} W^-(z)(e)^{1/2} & 0 \\ 0 & W^-(z)(e)^{-1/2}
  \end{array} \right), \quad p=\left(
  \begin{array}{cc} 0 & 1 \\ 1 & 0  \end{array}\right),
\quad l=\left(\begin{array}{cc} -1 & 1 \\ -1 & 0
\end{array}\right)$$
and $r=l^{-1}$. The matrix $\gamma(e)$ represents the isometry
with fixed points $0,\infty \in \mc\mathbb{P}^1$ and mapping $1$
to $W^-(z)(e)$.  The elliptic elements $p$ and $l$ send
$(0,1,\infty)$ to $(\infty,1,0)$ and $(\infty, 0,1)$ respectively.
For the portion of $\gamma$ on the left of Figure \ref{matrixrep},
if $\gamma$ turns to the left after crossing $e$ the parallel
transport operator is $\gamma(e)\cdot p\cdot l$, while it is
$\gamma(e)\cdot p\cdot r$ if $\gamma$ turns to the right.  (The
matrix multiplication is on the right, as is the action of
$PSL(2,\mc)$ on the total space of $F_{\tilde{\rho}}$). In the
middle and right pictures the parallel transport operators along
the portion of $\gamma$ are given by $\gamma(e)\cdot l$ or
$\gamma(e)\cdot p\cdot l$, and $\gamma(e)\cdot p\cdot r$ or
$\gamma(e)\cdot r$ respectively. We see that the action of $p$,
$l$ and $r$ depends on the reordering of the vertices after the
mapping $\gamma(e)$. If $\gamma$ intersects $e$ negatively, we
replace $\gamma(e)$ with $\gamma(e)^{-1}$ in the above
expressions. A similar recipe applies for any other branching of
the two triangles glued along $e$.

Continuing this way each time $\gamma$ crosses an edge of $T'$
until it comes back to $q$, we get an element $[W(z)]_q(\gamma)
\in PSL(2,\mc)$ that depends only on the homotopy class of
$\gamma$ based at $q$ and coincides with $\tilde{\rho}(\gamma)$,
because of the identity $D_z(\gamma \cdot x) = [W(z)]_q(\gamma)
\cdot D_z(x)$ for any $x \in \tilde{F}$, where $D_z: \tilde{F}
\rightarrow \mc\mathbb{P}^1$ is a $\tilde{\rho}$-equivariant
pseudo-developing map from the universal cover $\tilde{F}$ of $F$.
Varying $\gamma$, we eventually get the representation $[W(z)]_q
=\tilde{\rho}: \pi_1(F,q) \rightarrow PSL(2,\mc)$.
\smallskip

\noindent{\bf ($-$)-exponential $\Ii$-parameters and pleated
surfaces.} When $\Cc^-$ has no trivial entries, there is a nice
interpretation of the parameter space $W^-(T,b,\mathcal{C})$ in
terms of pleated hyperbolic surface structures on $S\setminus V$.
(see eg. \cite{CEG}; compare also with \cite{BoL}, Section 8.)
\smallskip

As before, let $T'$ be a triangulation of $S$ with vertices $V$,
viewed as an ideal triangulation of $S^{\circ} = S \setminus V$. A
{\it pleated surface} (with pleating locus $T'$) is a pair
$(\tilde{\mathfrak{f}},r)$, where $r: \pi_1(S^{\circ}) \rightarrow
PSL(2,\mc)$ is a group homomorphism (not up to conjugacy), and
$\tilde{\mathfrak{f}} : \widetilde{S^{\circ}} \rightarrow \mh^3$
is a map from the universal cover $\widetilde{S^{\circ}}$ of
$S^{\circ}$ such that:
\begin{itemize}
\item $\tilde{\mathfrak{f}}$ sends homeomorphically each component of
  the preimage $\tilde{T'}$ of $T'$ in $\widetilde{S^{\circ}}$ to a
  complete geodesic in $\mh^3$;

\item $\tilde{\mathfrak{f}}$ sends homeomorphically the closure of
  each component of $\tilde{S^{\circ}} \setminus \tilde{T'}$ to an
  ideal triangle in $\mh^3$;

\item $\tilde{\mathfrak{f}}$ is $r$-equivariant, that is, for all $x
  \in \tilde{S^{\circ}}$, $\gamma \in \pi_1(S^{\circ})$ we have
  $\tilde{\mathfrak{f}}(\gamma x) = r(\gamma)
  \tilde{\mathfrak{f}}(x)$.
\end{itemize}
Two pleated surfaces $(\tilde{\mathfrak{f}},r)$ and
$(\tilde{\mathfrak{f}}',r')$ are said {\it isometric} if there
exists an isometry $A \in PSL(2,\mc)$ and a lift $\tilde{\phi} :
\tilde{S^{\circ}} \rightarrow \tilde{S^{\circ}}$ of an isotopy of
$S^{\circ}$ such that $\tilde{\mathfrak{f}}' = A \circ
\tilde{\mathfrak{f}} \circ \tilde{\phi}$ and $r'(\gamma) = A
r(\gamma) A^{-1}$ for all $\gamma \in \pi_1(S^{\circ})$. From
\cite{Bo} it is known that isometry classes of pleated surfaces
$(\tilde{\mathfrak{f}},r)$ are in one-one correspondence with
arrays $\{x_e\}_e$ of non zero complex numbers $x_e$ associated to
the edges $e$ of $T'$, the {\it exponential shear-bend
parameters}.

To a pleated surface $(\tilde{\mathfrak{f}},r)$ we can associate
the type $\mathcal{C}^-=(\mathcal{C}_1,\ldots,\mathcal{C}_r)$ of
the conjugacy class of $r$, defined as in Section \ref{Dpar}. For
instance, it is shown in \cite{Th} that the set of isometry
classes of pleated surfaces given by real positive shear-bend
parameters is real-analytic diffeomorphic to the Teichm\"uller
space of hyperbolic metrics on $S^{\circ}$, with totally geodesic
boundary completion or non compact finite area completion at
$v_i$, according to the type of $\mathcal{C}_i$, loxodromic or
parabolic.

In general, we can also associate a {\it sign} to each loxodromic
puncture. Namely, if $A$ is a small annulus neignborhood of $v_i$,
each connected component $\tilde{A}$ of the preimage of $A$ in
$\tilde{S^{\circ}}$ is the fixed point set of a subgroup $\pi$ of
$\pi_1(S^{\circ})$. This subgroup is the image of the fundamental
group of $A$ for some choice of base points and paths between
these base points. All the edges of $\tilde{T}'$ that meet
$\tilde{A}$ are sent by $\tilde{\mathfrak{f}}$ to geodesics lines
that meet at one of the two fixed points of $r(\pi)$. We specify
this fixed point by a sign, as it determines an orientation
(whence a generator) for the axis of the group $r(\pi)$. Since any
two subgroups $\pi$ as above are conjugated, for each puncture the
fixed point assignment is $r$-equivariant, so that the sign is
canonically associated to the puncture.

\smallskip

Recall that in Definition \ref{defcoc} the types $\mathcal{C}_i$
where assigned to the components of $\partial F$ endowed with the
boundary orientation. Let us remove this constraint, and associate
to each boundary component with loxodromic type an orientation
that we specify by a sign, positive for the boundary orientation,
and negative otherwise. Then, for each type $\mathcal{C}^-$ with
$l$ loxodromic entries and each $l$-uple of signs $s$, we get a
space $\Zz_I(T,b,\mathcal{C}^-,s)$. We define
$W^-(T,b,\mathcal{C},s)=W(\Zz_I(T,b,\mathcal{C}^-,s))$.
\begin{prop}\label{psI}{\rm For each type
    $\mathcal{C}$ with non trivial entries $\mathcal{C}_i$, $l$ being
    loxodromic, and for each $l$-uple of signs $s$, the space
    $W^-(T,b,\mathcal{C},s)$ coincides with the exponential shear-bend
    parameter space $\mathcal{P}\mathcal{S}_{\mathcal{C}}^s$ of
    isometry classes of pleated surfaces of type $\mathcal{C}$ and
    signs $s$:
\begin{equation}\label{polPS}
\mathcal{P}\mathcal{S}_{\mathcal{C}}^s = \lbrace \{x_e\} \in (\mc
  \setminus \{ 0 \})^{-3\chi(S^{\circ})} \mid \forall i=1,\ldots,r,\
  \prod x_e = \mu_i \rbrace.
\end{equation}
Here the product is over all edges $e$ with $v_i \in\partial e$
(counted with multiplicities), $\mu_i$ is $1$ if $\mathcal{C}_i$
is parabolic, and, if $\mathcal{C}_i$ is loxodromic, $\mu_i$ is
the dilation factor of the generator of $<\mathcal{C}_i>$
specified by the sign of $v_i$. (Hence, in either case this is an
eigenvalue of $\mathcal{C}_i$.)}
\end{prop}

\noindent{\it Proof.} By the results of \cite{Bo} recalled above,
each point of $W^-(T,b,\mathcal{C},s)$ is identified with the
family of exponential shear-bend parameters of an isometry class
of pleated hyperbolic surface on $S^{\circ}$ with pleated locus
$T'$. The definition of the map $\Theta: W^-(T,b,\mathcal{C},s)
\rightarrow \Rr_I(T,b,\mathcal{C})$ shows that this isometry class
has type $\mathcal{C}$. Furthermore, by using the recipe given
above for computing $\Theta$, we check that the holonomy of any
positively oriented (with respect to the boundary orientation)
small loop about the puncture $v_i$ is exactly $\mathcal{C}_i$.
The upper left diagonal entry is just the product of a square root
of the exponential $\Ii$-parameters at the edges with endpoint
$v_i$. From (\ref{polPS}), which is an easy consequence of results
in \cite{Bo} (see sections 12.2-12.3 in that paper), we deduce
that $W^-(T,b,\mathcal{C},s) \subset
\mathcal{P}\mathcal{S}_{\mathcal{C}}^s$.

Conversely, for any $(\tilde{\mathfrak{f}},r) \in
\mathcal{P}\mathcal{S}_{\mathcal{C}}^s$ we have $ {\rm conj}(r) \in
\Rr_I(T,b,\mathcal{C})$, because $r$ is injective. Also, by
\cite{BoL}, Proposition 33, it is known that the isometry class of
$(\tilde{\mathfrak{f}},r)$ is determined by $r$ and the signs $s$.
As the map $\Theta$ is onto, let us take $z \in
\Zz_I(T,b,\mathcal{C}^-,s)$ with ${\rm conj} \circ f(z) = {\rm conj}(r)$
and consider the isometry class of pleated hyperbolic surfaces
associated to $W^-(z)$. There is a representative
$(\tilde{\mathfrak{f}}',r)$ with the same holonomy $r$. As the signs
of $(\tilde{\mathfrak{f}}',r)$ are $s$, the same as for
$(\tilde{\mathfrak{f}},r)$, the two pleated surfaces coincide. So
$\mathcal{P}\mathcal{S}_{\mathcal{C}}^s \subset
W^-(T,b,\mathcal{C},s)$.\hfill$\Box$

\begin{figure}[ht]
\begin{center}
\includegraphics[width=3.5cm]{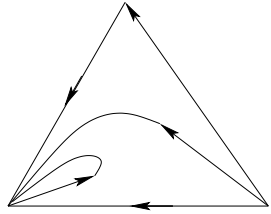}
\caption{\label{more-e} More economic e-triangulation.}
\end{center}
\end{figure}

\begin{remark}\label{branchwhy}{\rm
Exponential shear-bend parameters do not depend on branchings but the
orientation of $F$ (taking the opposite branching edge orientation
simultaneously exchanges the left and right triangles). However,
branchings govern all choices in QHFT tensors. Also, they allow us to
interpret exponential shear-bend parameters as exponential
$\Ii$-parameters, thus coming from $1$-cocycles representing {\it
  arbitrary} $PSL(2,\mc)$-characters on the triangulated boundary of
{\it arbitrary} compact orientable $3$-manifolds.  When $\Cc$ has no
trivial entries, the maps $W^\pm$ define {\it decorated} shear-bend
parameter spaces similar to those occurring in \cite{Kmod, Pe}.}
\end{remark}
\begin{remark}\label{varpar}
{\rm For types $\mathcal{C}$ with no trivial entries, we can use
simpler $e$-triangulations of $F$ (see Figure \ref{more-e}): for
each of the triangles $t_v$ of the base surface $S$ we remove the
interior of a {\it monogon} inside $t_v$, and triangulate the
resulting quadrilateral by adding an edge $e_v$ with endpoints $v$
and the $b$-output vertex of the opposite edge. We extend $b$ by
orienting $e_v$ from that vertex to $v$. Proposition \ref{topcomp}
applies to the cocycle parameters based on such
$e$-triangulations, for points of $\Rr(g,r)$, and also the
treatment of exponential $\Ii$-parameters works as well.}
\end{remark}

\section{The QHFT bordism category}\label{BORDCAT}
We define first a {\it topological} $(2+1)$-bordism category. Then we
will give it more structure, including the parameter spaces of the
previous sections.

\subsection{Marked topological bordisms}\label{top-bord}
Like in Section \ref{PARAM}, for every $(g,r)\in \N \times \N$
such that $g\geq 0$, $r>0$, and $r>2$ if $g=0$, fix a compact
closed oriented base surface $S$ of genus $g$ with a set
$V=\{v_1,\ldots,v_r\}$ of $r$ marked points. Denote by $-S$ the
same surface with the opposite orientation, and write $*S$ for
$S=+S$ or $-S$. Moreover, fix a set of disjoint embedded closed
segments $a_{v_i}$ in $S$ such that $a_{v_i}$ has one end point at
$v_i$. We say that {\it $v_i$ is framed by $a_{v_i}$}.

We say that two orientation preserving diffeomorphisms $\phi_1: *
S \to \Sigma_1$ and $\phi_2: * S \to \Sigma_2$ are equivalent if
there is an orientation-preserving diffeomorphism $h: \Sigma_1 \to
\Sigma_2$ such that $(\phi_2)^{-1}\circ h \circ \phi_1$ pointwise
fixes the segments $a_{v_i}$ and is isotopic to the identity
automorphism of $S$ relatively to $\{a_{v_i}\}$. We write $[*
S,\phi]$ for such an equivalence class.

Let $Y$ be an oriented compact $3$-manifold with (possibly empty)
boundary $\partial Y$, with an input/output bipartition $\partial Y =
\partial_-Y \cup \partial_+ Y$ of the boundary components (we say that
$\partial_- Y$ is ``at the bottom'' of $Y$, while $\partial_+Y$ is
``on the top'').  Each boundary component inherits the boundary
orientation, via the usual convention {\it last is the ingoing
normal}. Let $L_\Ff$ be a properly embedded {\it non-empty} framed
tangle in $Y$. This means that $L_\Ff$ is a disjoint union of properly
embedded orientable {\it ribbons}. We split $L_\Ff = (L_\Ff)_i\cup
(L_\Ff)_b$, where $(L_\Ff)_i$ is the {\it internal} part of $L_\Ff$
made by its closed connected components, homeomorphic to the annulus
$S^1 \times [0,1]$, while $(L_\Ff)_b$ is the union of the components
homeomorphic to the quadrilateral $I\times [0,1]$.  We consider
$L_\Ff$ up to proper ambient isotopy. For every boundary component
$\Sigma$ of $Y$, we assume that $(L_\Ff)_b\cap \Sigma \ne \emptyset$,
and that $(L_\Ff)_b\cap \Sigma$ consists of at least two segments if
$g(\Sigma)=0$ (we do not require that every component of $(L_\Ff)_b$
goes from $\partial_-$ towards $\partial_+$). Hence, on each
boundary component $\Sigma$ of $Y$ we have a set of marked points
framed by $L_\Ff \cap \Sigma$. Also, associated to $\partial_\pm$ we
have a finite disjoint union $\textstyle \alpha_\pm = \coprod_{\phi(*
S) \in \partial_\pm} [*S,\phi]$ of equivalence classes of
diffeomorphisms as above.

\smallskip

Consider the topological $(2+1)$-bordism category with objects the
empty set and any finite union of the $[* S,\phi]$, and morphisms the
triples $(Y,L_\Ff,\alpha_\pm)$ as above. We say that
$(Y,L_\Ff,\alpha_\pm)$ is a {\it bordism from $\alpha_-$ to $\alpha_+$
with support $(Y,L_\Ff)$}. We allow the case when $Y$ is a closed
manifold, so that $\partial Y = \emptyset$ and $(Y,L_\Ff)$ is a
morphism from the empty set to itself. We stress that $L_\Ff$ is
non-empty in any case.

We can reformulate this category in a setup closer to that of the
phase space parameters of Section \ref{PARAM}, as follows. If we cut
open each $a_{v_i}$ in $S$ we get an oriented surface $F$ with $r$
boundary bigon components. This is the domain of an elementary object
$[* F,\phi]$, where the target surfaces $\Sigma$ have now $r$ boundary
components, and the diffeomorphisms $\phi$ are considered up to
isotopy rel($\partial$).

\smallskip

Consider a bordism $(Y,L_\Ff,\alpha_\pm)$. On the boundary of each
ribbon component of $L_\Ff$ we keep track of a {\it tangle} line
$X\times \{0\}$ ($X=S^1$ or $I$) for the corresponding component of
the unframed tangle $L$, and there is a {\it longitudinal} line
$X\times \{1\}$ that specifies the framing of the normal bundle of the
parallel tangle line. These make a pair
$\bar{\lambda}=(\lambda,\lambda')$ of parallel {\it unframed} tangles
in $Y$. Cutting open each ribbon we get a $3$-manifold with corners
$\tilde{Y}$. The boundary $\partial \tilde{Y}$ has two {\it
  ``horizontal''} parts $\partial_\pm \tilde{Y}$ contained in
$\partial_\pm Y$, and a ``vertical'' {\it tunnel} part $\widetilde
{L}_\Ff$. The horizontal parts intersect the tunnel part at the corner
locus; this is a union of bigons contained in $\partial Y$. Each
boundary component $\Sigma$ of $Y$ corresponds to a horizontal
boundary component of $\tilde{Y}$, still denoted by $\Sigma$. Each
tunnel boundary component is made by the union of two copies of $X
\times (0,1)$, glued each to the other at $\lambda \cup \lambda'$. The
horizontal boundary components are the targets of elementary objects
$[* F,\phi]$, and each triple $(\tilde{Y},\tilde{L}_\Ff)$ supports a
morphism between such objects. Clearly, we can recover $(Y,L_\Ff)$
from $(\tilde{Y},\tilde{L}_\Ff)$, so that we have two equivalent
settings to describe the same topological bordism category.

\subsection{Boundary structures}\label{QHB-def}
\begin{figure}[ht]
\begin{center}
\includegraphics[width=4cm]{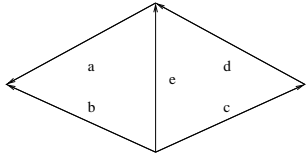}
\caption{\label{LogIfig} Notations for Log-$\Ii$-parameters.}
\end{center}
\end{figure}

We equip a marked topological bordism with additional boundary
structures by using the notions of section \ref{PARAM}. Fix an
$e$-triangulation $(T,b)$ of $F$. Following Remark \ref{W_mean}
and the definition of log-branches, we put:

\begin{defi}\label{log-I}{\rm Let $z \in \Zz_I(T,b)^\pm$, and $e$ be a
non boundary edge of $T$. Denote ${\rm l}_z(e)$ the canonical
log-branch of $W^\pm(z)(e)$, computed from the idealization of Star$(e)$ as in Lemma \ref{canonicalflat} and before Definition \ref{phasespace}. For any
collection $m=\{m_\lambda\}_\lambda$ of integers, one for each edge
$\lambda$ in $T$, the {\it classical Log-$\Ii$-parameter} of $(z,m)$
at $e$ is $${\rm l}_{(z,m)}(e) = {\rm l}_z(e)+
\sqrt{-1}\pi(m_{a}+m_{c}-m_{b}-m_{d}),$$ where $a$, $\ldots, d$ make
Link$(e)$ with $a$, $c$ opposite and $e$, $a$ have coherent branching
orientations (see Figure \ref{LogIfig}). Similarly, for every $N>1$
and any other collection $n=\{n_\lambda\}_\lambda$ of integers, the
{\it quantum Log-$\Ii$-parameter} at $e$ is $${\rm
  l}_{(z,m)}(e)+\sqrt{-1}\pi
N(m_{a}+m_{c}-m_{b}-m_{d})-*_b\sqrt{-1}\pi(N+1)(n_{a}+n_{c}-n_{b}-n_{d}).$$
We call the collections $$f=\{f_e\}_e=\{({\rm
  l}_{(z,m)}(e)-\log(W^\pm(z)(e)))/\pi\sqrt{-1}\}_e$$ and
$$c=\{c_e\}_e=\{n_{a}+n_{c}-n_{b}-n_{d}\}_e$$ the flattenings of
$(z,m)$ and the charge of $n$ respectively, and we denote
generically by $\Ll$ any such a system of classical or quantum
($\pm$)-Log-$\Ii$-parameters.}
\end{defi}
Remark that if $z\in \Zz_I(T,b,\mathcal{C}^-)$ and $\mathcal{C}$
has no trivial entries, then the ($-$)-Log-$\Ii$-parameters depend
on $W^-(z)$, not on $z$, because of Proposition \ref{no-gauge}.

\begin{defi} \label{QHBdef} {\rm The QHFT category is the
    $(2+1)$-bordism category with objects the empty set
    and any finite union of the $[*F,(T,b,\Ll),\phi]$, where $\Ll$ is a
    system of ($+$)-Log-$\Ii$-parameters, and morphisms the 4-uples
    $(\tilde{Y},\tilde{L}_\Ff,\rho,\alpha_\pm)$, where: $\rho$ is a
    conjugacy class of $PSL(2,\C)$-valued representations of
    $\pi_1(\tilde{Y} \setminus \tilde{L}_\Ff)$; $\alpha_\pm$ are QHFT
    objects with targets $\partial_\pm \tilde{Y}$, such that for every
    $[* F,(T,b,\Ll),\phi]$, the character $\phi^*(\rho)$ coincides
    with $\Theta(W^+(z))$. We say that $(\tilde{Y},\tilde{L}_\Ff,
    \rho)$ is the support of a QHFT bordism from $\alpha_-$ to
    $\alpha_+$, and that $\alpha_\pm$ is a QHFT surface.}
\end{defi}
Given bordisms $\Bb$ and $\Bb'$ from $\alpha_-$ to $\alpha_+$ and
$\alpha_-'$ to $\alpha_+'$, respectively, assume that $\beta_+$
and $\beta'_-$ are subobjects of $\alpha_+$ and $\alpha'_-$ that
coincide up to the change of orientation.
\begin{defi}\label{composition}{\rm
The bordism $\Bb''$ from $\alpha''_- = \alpha_- \cup (\alpha'_-
\setminus \beta'_-)$ to $\alpha''_+ = \alpha'_+ \cup (\alpha_+
\setminus \beta_+)$ obtained by gluing $\Bb$ and $\Bb'$ along
$\beta_+$ is called the {\it composition} of $\Bb$ followed by
$\Bb'$. We write $\Bb'' = \Bb'*\Bb$. }
\end{defi}

\noindent {\bf Examples from hyperbolic geometry.} Any
topologically tame hyperbolic $3$-manifold $Y$ with hyperbolic
holonomy $\rho$ and a tangle $L_\Ff$ of singularities makes a QHFT
bordism. More specifically, any geometrically finite non compact
complete hyperbolic $3$-manifold $Y$ defines a triple
$(Y',L_\Ff,\rho)$ with a {\it non empty} link $L_\Ff$, as follows.
The manifold $Y$ has a natural compactification $\bar{Y}$, with $Y
\cong {\rm Int}(\bar{Y})$, which is a ``pared'' manifold
$(\bar{Y},P)$. Here $P$ is a union of disjoint tori or annuli
embedded in the boundary of $\bar{Y}$. The tori correspond to the
cusps of $Y$. Each annulus $A$ of $P$ comes from a couple of cusps
on some boundary component of (a small neighborhood of) the convex
core of $Y$; $A$ is fibered by geodesic arcs. If $A$ is separating
the cusps belong to different components. Define $Y''$ as the
result of attaching a $2$-handle to $\bar{Y}$ at each annulus $A$,
so that $P$ is contained in the interior of $Y''$ and is
transverse to $\partial Y''$. Equivalently, $Y''$ contains a
properly embedded framed $1$-tangle $L_\Ff'$ made by the cocores
of the $2$-handles, the framing being determined by the fibration
by intervals of the annuli of $P$. Let us choose a framing at each
torus of $P$. By Dehn filling we get a manifold $Y'$, and $Y''$ is
the exterior in $Y'$ of the union $L_\Ff''$ of the framed cores of
the filling solid tori. Hence, if every boundary component of the
convex core of $Y$ contains at least one cusp, associated to $Y$
and the cusp framings we have $(Y',L_\Ff,\rho)$, where $L_\Ff
=L_\Ff' \cup L_\Ff''$ intersects all the boundary components of
$Y'$, and $\rho$ is a $PSL(2,\mc)$-character of $Y'\setminus L_\Ff
\cong Y$. If furthermore $\rho$ is the holonomy of a complete
hyperbolic metric on $Y$ and $Y$ has infinite volume ends, then
$\partial Y'$ is non empty. We can give the triple
$(Y',L_\Ff,\rho)$ a natural boundary structure $\alpha_- \cup
\alpha_+$, induced by exponential $\Ii$-parameters of the pleated
surfaces in the boundary of the convex core (see eg. \cite{Ep}).

\section {The QHFT functor}\label{TENSOR} Consider a QHFT bordism
$\Bb =(\tilde{Y},\tilde{L}_\Ff,\rho,\alpha_\pm)$. For every odd
integer $N\geq 1$, we associate a finite dimensional complex
linear space $V(\alpha_\pm)$ to $\alpha_{\pm}$, and a linear map
$\Hh_N(\Bb) : V(\alpha_-) \to V(\alpha_+)$ to $\Bb$, well-defined
up to sign and multiplication by $N$th roots of unity. This
defines a (moderately projective) functor $\Hh_N: {\rm QHB}
\rightarrow {\rm Vect}$, where Vect is the tensor category of
complex linear spaces. The construction immediately implies that
$\Hh_N$ is a modular functor, in the sense of \cite{T}, III.1.2.

\subsection {From QHFT bordisms to QHG-triangulated pseudomanifolds}
\label{Z(B)}
First we associate to $\Bb$ a pseudomanifold $Z(\Bb)$.
Fill each tunnel boundary component of $\tilde{Y}$
with a solid tube, thus recovering a copy of the manifold $Y$. The
cores of the solid tubes make a parallel unframed copy $\lambda''$
of $L$. We define $Z(\Bb)$ as the result of collapsing to one
point each component of $\lambda''$. In other words, we glue to
each tunnel component of $\tilde{Y}$ the {\it oriented topological
cusp} $\hat{\bf C} = B \times [0,+\infty]/(B \times \{\infty\})$
with base equal to either $B=S^1\times [-1,1]$ or $B=S^1\times
S^1$.
\smallskip

Next we describe a procedure to convert $Z(\Bb)$ to a
distinguished QHG-triangulated pseudomanifold. We refer to the
notions introduced in Subsection \ref{ToGeoSpec}.

\smallskip

We say that a branched triangulation $(T,b)$ of $B$ as above is {\it
  admissible} if $B \cap \partial_\pm \tilde{Y}$ and the tangles
  $\lambda$, $\lambda'$ are covered by the $1$-skeleton. We denote
  $(\hat{T},\hat{b})$ the branched triangulation of $\hat{\bf C}$,
  where $\hat{T}$ is the cone over $T$ from the non manifold point,
  say $\infty$, and $\hat{b}$ extends $b$ so that $\infty$ is a pit
  for every branched tetrahedron of $\hat{T}$.  Assume we are given an
  idealizable $PB^+(2,\mc)$-valued cocycle $z$ on an admissible
  triangulation of $B$. The idealization of $z$ determines for each
  $2$-simplex of $B$ a face of an ideal hyperbolic tetrahedron with
  further vertex at $\infty$ (see Figure \ref{tcusp}, where opposite
  vertical triangles are identified).  Since the fundamental group of
  $B$ is Abelian, the resulting family of oriented ideal hyperbolic
  tetrahedra actually makes an $\Ii$-triangulation
  $(\hat{T},\hat{b},w)$, which we call an {\it $\Ii$-cusp}. By
  conjugating if necessary, we see that $\Ii$-cusps make sense also
  when $z$ takes values more generally in $PSL(2,\mc)$. We get
  flattenings similarly as in Lemma \ref{canonicalflat}: at a corner
  of a $2$-simplex formed by edges $e_l$ and $e_r$ we put the
  difference of the logarithms of the vectors in $\mc$ associated by
  the idealization to $e_l$ and $e_r$, respectively.
\begin{figure}[ht]
\begin{center}
\includegraphics[width=4cm]{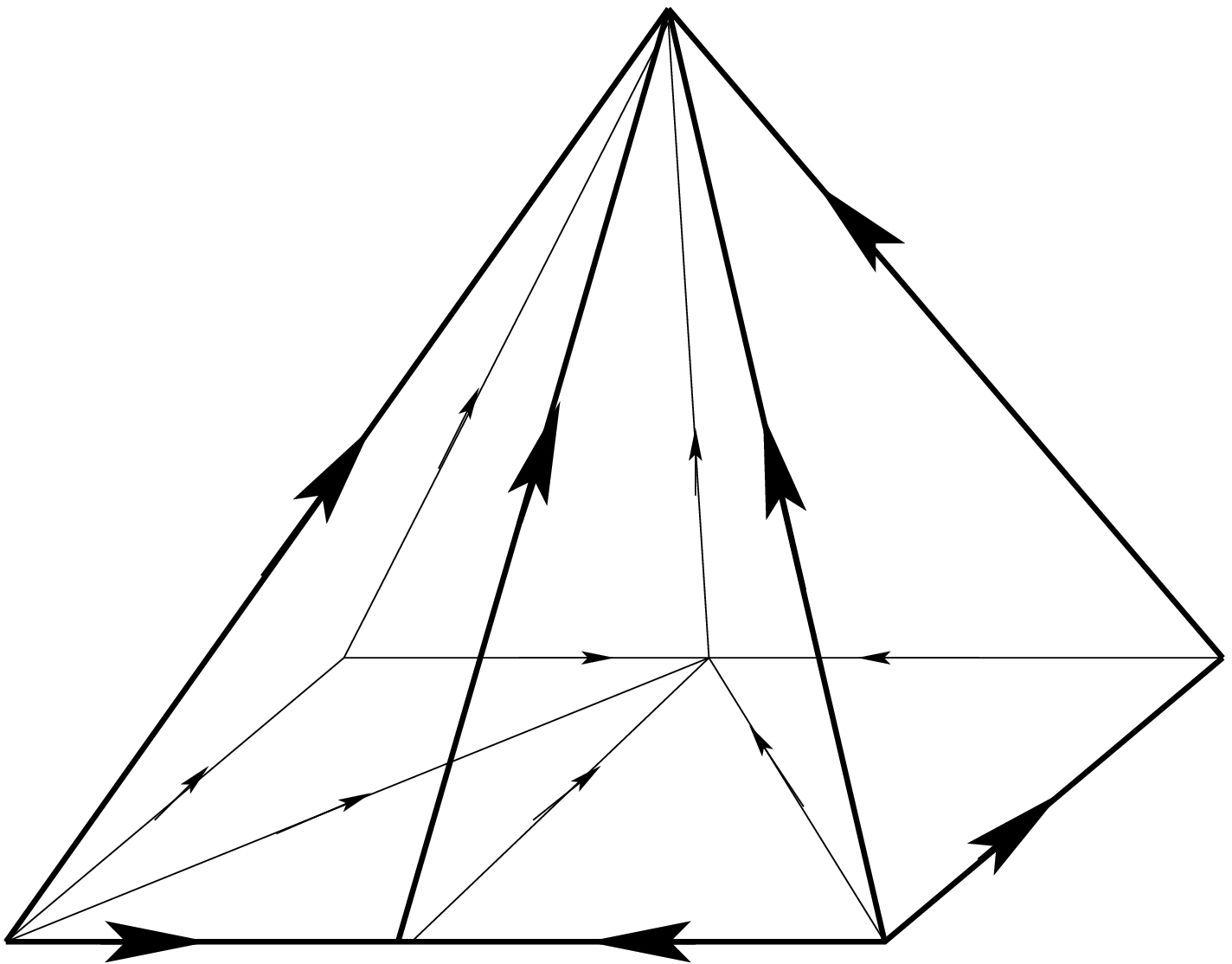}
\caption{\label{tcusp} An $\Ii$-cusp.}
\end{center}
\end{figure}

\begin{defi}\label{Dtriang}{\rm A {\it $\Dd$-triangulation} of
$\Bb =(\tilde{Y},\tilde{L}_\Ff,\rho,\alpha_\pm)$ consists of a
$4$-uple $\Kk =(K,\bar{H},b,z)$ where:

\smallskip

(a) $(K,b)$ is a branched triangulation of $\tilde{Y}$ extending
that on $\partial_\pm \tilde{Y}$, and inducing an admissible cusp
base triangulation at each tunnel component of $\tilde{L}_{\Ff}$.

\smallskip

(b) The $1$-dimensional subcomplex $\bar{H}= H\cup H'$ of $K$ is
ambiently isotopic to the tangle $\bar{\lambda}=\lambda \cup
\lambda'$, and $\bar{H}$ contains
all the vertices of $K$.

\smallskip

(c) $z$ is an idealizable $PSL(2,\C)$-valued $1$-cocycle on $(K,b)$
such that:

(i) the conjugacy class of $PSL(2,\mc)$-representations of $\pi_1(\tilde{Y})$
associated to $z$ coincides with $\rho$;

(ii) the ($+$)-exponential $\Ii$-parameters given by the
restriction of $z$ to $\partial_{\pm} \tilde{Y}$ coincide with
that of the objects $\alpha_- \cup \alpha_+$ (see Definition \ref{log-I}).
\smallskip

(d) the restriction of $z$ to each vertical tunnel component of
$\tilde{L}_{\Ff}$ takes values in the Borel subgroup $PB^+(2,\mc)$
of $PSL(2,\mc)$.}
\end{defi}
For any $\Dd$-triangulation $\Kk =(K,\bar{H},b,z)$ of $\Bb$ we get
a distinguished $\Ii$-triangulation $\Kk_\Ii=(K,\bar{H},b,w)$ of
$(Z(\Bb), \bar{\lambda})$ by gluing the idealization of $\Kk$ with
the $\Ii$-cusp given by the cocycle at each tunnel component. Note
that $\bar{H}$ contains all the vertices of $K_I$ that are
manifold points.
\begin{defi}\label{flat/charge-def} {\rm We say that $\mathcal{T}(\mathcal{B})
=(\Kk_\Ii,f,c) = (K,\bar{H},b,w,f,c)$ is a {\it distinguished
  flat/charged $\Ii$-triangulation} of $(Z(\Bb),\bar{\lambda})$ if it
satisfies Definition \ref{DISFLATCHARGE}, and at every boundary edge
of $Z(\Bb)$ the total (classical or quantum) log-branch of
$(\Kk_\Ii,f,c)$ of (\ref{totw}) coincides with the (classical or
quantum) Log-$\Ii$-parameter of the boundary object $\alpha_- \cup
\alpha_+$. }
\end{defi}
Recall the cohomological weights from section \ref{ToGeoSpec}.  These
notions still make sense for distinguished flat/charged
$\Ii$-triangulations, where the homology of $\partial T_0$ is replaced
with that of the tunnel components $\tilde{L}_\Ff$. Since we have
weights $(h_f,k_f)$ for log-branches and $(h_c,k_c)$ for charges at
the same time, we will denote them $(h,k)=((h_f,h_c),(k_f,k_c))$. We
can also define, in the very same way, boundary weights $k_f\in
H^1(\partial \tilde{Y}; \mc)$, but these are completely encoded by
$\alpha_\pm$. We have:
\begin{teo} \label{teoTOP} For every bordism
$\Bb=(\tilde{Y},\tilde{L}_\Ff,\rho,\alpha_\pm)$ and every $(h,k)
\in H^1(\tilde{Y};\Z/2\Z) \times H^1(\tilde{L}_\Ff; \mc)$
satisfying the properties (\ref{weightcondition}), there are
distinguished flat/charged $\Ii$-triangulations $\Tt(\Bb)$ of
$Z(\Bb)$ with weight $(h,k)$, and any two are QHG-isomorphic.
\end{teo}

\noindent {\it Proof.} The existence of distinguished
$\Ii$-triangulations $\Kk_\Ii$ of $Z(\Bb)$ follows from a tedious
but straightforward generalization of Theorem 4.13 in \cite{BB1}.
Global flattenings and integral charges with arbitrary weight
exist on the double $D\Kk_\Ii$ of $\Kk_\Ii$ by the results
recalled with Theorem \ref{theoweight}.

Consider the ($+$)-Log-$\Ii$-parameters $\{W^\pm(z)(e)\}_e$ at
$\alpha_\pm$. They are in one-one correspondence with the interior
edges of the corresponding $e$-triangulations, which is less than the
cardinality of the families $m$ and $n$ used to defined flattenings
and charges in Definition \ref{log-I}. Hence any family of
determinations of the logarithms of the $W^\pm(z)(e)$ is a system of
($+$)-Log-$\Ii$-parameters. Also, Lemma \ref{canonicalflat} implies
that any system of ($+$)-Log-$\Ii$-parameters at $\alpha_\pm$ extends
to a distinguished flat/charge $\Ii$-triangulation of the pseudo
manifold obtained from the trivial cylinders over
$\partial_{\pm}\tilde{Y}$ by collapsing to a point each annulus of
$\partial( \partial_{\pm}\tilde{Y}) \times [-1,1]$.

This means that any QHFT surface bounds a QHFT bordism, and that
for the bordism $\Bb$ there are flat/charges on $D\Kk_\Ii$ whose
restriction to $\Kk_\Ii$ induce the Log-$\Ii$-parameters of
$\alpha_\pm$. Hence we get global flat/charges as in Definition
\ref{flat/charge-def}. In fact the affine spaces of flat/charges
on $D\Kk_\Ii$ project onto that on $\Kk_\Ii$ compatible with
$\alpha_\pm$ (see the end of section \ref{QHG}). Then, the
Mayer-Vietoris exact sequence in cohomology for the triad
$(DK_\Ii,\Kk_\Ii,-\Kk_\Ii)$ shows that $(h,k)$ is induced by some
weight on $D\Kk_\Ii$. As we can choose the latter arbitrarily,
this concludes the proof of the first claim.

The second is harder, but follows strictly from the arguments in
the proof of Theorem 6.8 (2) in \cite{BB2}. The only new
ingredient is the presence of $\Ii$-cusps, which mimic the ends of
cusped manifolds treated in that paper. \hfill$\Box$

An alternative characterization of classical/quantum log-branches
of distinguished flat/charged $\Ii$-triangulations $\Tt(\Bb)$
follows from Remark \ref{univfc}.
\subsection{Amplitudes}\label{QHFTTENS} Fix an odd positive integer $N$.
Write $V=\C^N$, with the canonical basis $\{e_i\}$, and $V^{-1}$
for the dual space. Both are endowed with the hermitian inner
product with orthonormal basis the vectors $e_0$ and
$(e_i+e_{N-i})/\sqrt{2}$, $i=1,\ldots,N-1$.

Recall the notations of section \ref{PARAM}. For each base surface
$F$ fix an $e$-triangulation $(T,b)$, an ordering of the set
$T^{(2)}$ of $2$-simplices, and let $V(T,b) = \otimes_{t\in
T^{(2)}} V^{\sigma_\pm (t)}$. Given a QHFT bordism
$\Bb=(\tilde{Y},\tilde{L}_\Ff,\rho,\alpha_\pm)$ with a
distinguished flat/charged $\Ii$-triangulation $\Tt(\Bb)$, the
trace tensor in (\ref{deftt}) is a morphism $\Hh_N(\Tt(\Bb)) \in
{\rm Hom}(V(\alpha_-),V(\alpha_+))$, where $V(\alpha_\pm)$ is the
tensor product of isomorphic copies of the spaces $V(T,b)$ over
the (ordered) components of $\alpha_\pm$.
\begin{teo}\label{main} The morphism $\Hh_N(\Bb,h,k)=\Hh_N(\Tt(\Bb))$ does not
depend on the choice of $\Tt(\Bb)$ up to sign and multiplication by
$N$th roots of unity, and $\Hh_N(\Bb,h,k)=\Hh_N(\Bb,h',k')$ if the
mod($N\mz$) reductions of $h$ and $h'$ (resp. $k$ and $k'$) are the
same, that is, if we have $k-k' \in H^1(\tilde{L}_\Ff;\mz)$ and $k-k'=0 \in
H^1(\tilde{L}_\Ff;\mz/N\mz)$, and similarly for $h$ and $h'$. Moreover,
there is no sign ambiguity if we restrict to even valued flattenings
as in Lemma \ref{canonicalflat} (hence with $h=0$). We call
$\Hh_N(\Bb,h,k)$ the {\rm amplitude} of $(\Bb,h,k)$.
\end{teo}
\noindent {\it Proof.} The result up to sign is an immediate
consequence of Proposition \ref{fund_func} and the last claim in
Theorem \ref{teoTOP}. For the dependance with respect to the
mod($N\mz$) reductions of weights, we note that the associated systems
of $N$th roots of moduli (see Remark \ref{univfc}) are connected by
QHG isomorphisms. Indeed, the difference $k-k' \in
H^1(\tilde{L}_\Ff;N\mz)$ coincides with $\gamma(f) -
\gamma(f')/\sqrt{-1}\pi$ and $\gamma(c) - \gamma(c')/\sqrt{-1}\pi$ for
some $f$, $f'$ and $c$, $c'$, and similarly for $h-h'$ and the
$\gamma_2$ maps. By first considering $(k-k')/N$ we can eventually
take the collections of values of $f-f'$ and $c-c'$ in $N\mz$, and
equal. Hence the conclusion follows from (\ref{qvar}) and
(\ref{symqmatdil}). For even valued flattenings, the
claim follows from the fact that Proposition \ref{fund_func} has no
sign ambiguity. For clarity let us state the result with some details.

Consider all the possible branching configurations of $2
\leftrightarrow 3$ QHG transits, up to obvious symmetries. They are
obtained from a single one by composing with the transpositions
$(01)$, $(12)$, $(23)$ and $(34)$ of the vertices (ordered in
accordance with the branching). Any such a transposition changes the
matrix dilogarithm of each tetrahedron by matrix conjugation and
multiplication by a determined scalar factor (see Corollary 5.6 of
\cite{BB2}). Now, there is a prefered ``basic'' $2 \leftrightarrow 3$
QHG transit, for which {\it even} flattenings give no sign ambiguity
in Proposition \ref{fund_func} (see the proof of Theorem 5.7 in
\cite{BB2}). The corresponding branching configuration is defined by
ordering the vertices as $1$ and $3$ on the bottom and top of the
central edge, and $0$, $2$, $4$ counterclockwise as viewed from the
$3$rd vertex. The following table describes for each tetrahedron
$\Delta^i$, opposite to the $i$-th vertex, the scalar factors induced
by the above transpositions (we put $v=\exp(\sqrt{-1}\pi(1-N)/2N)$,
which is $-\zeta^{-(m+1)(1-N)/2}$ in the notations of \cite{BB2}, and
the $c_k^l$ are charge values):
$$\begin{tabular}{||c||c|c||c|c|c||} \hline &
  \rule{0cm}{0.4cm}{$\Delta^1$}& $\Delta^3$ & $\Delta^0$ & $\Delta^2$
  & $\Delta^4$ \\ \hline

\rule{0cm}{0.8cm} \raisebox{0.2cm}{$(01)$} & \raisebox{0.2cm}{$0$} &
\raisebox{0.2cm}{$v^{c_0^3}$} & \raisebox{0.2cm}{$0$} &
\raisebox{0.2cm}{$v^{c_0^2}$} & \raisebox{0.2cm}{$v^{c_0^4}$}\\ \hline

\rule{0cm}{0.8cm} \raisebox{0.2cm}{$(12)$} & \raisebox{0.2cm}{$0$} &
\raisebox{0.2cm}{$v^{c_1^3}$} & \raisebox{0.2cm}{$v^{c_0^0}$} &
\raisebox{0.2cm}{$0$} & \raisebox{0.2cm}{$v^{c_1^4}$} \\ \hline

\rule{0cm}{0.8cm} \raisebox{0.2cm}{$(23)$} &
\raisebox{0.2cm}{$v^{c_1^1}$} & \raisebox{0.2cm}{$0$} &
\raisebox{0.2cm}{$v^{c_1^0}$} & \raisebox{0.2cm}{$0$} &
\raisebox{0.2cm}{$v^{c_0^4}$} \\ \hline

\rule{0cm}{0.8cm} \raisebox{0.2cm}{$(34)$} &
\raisebox{0.2cm}{$v^{c_0^1}$} & \raisebox{0.2cm}{$0$} &
\raisebox{0.2cm}{$v^{c_0^0}$} & \raisebox{0.2cm}{$v^{c_0^2}$} &
\raisebox{0.2cm}{$0$} \\ \hline
\end{tabular}$$
Because of (\ref{preserveq}), we see that the scalars at both sides
are equal. Hence, for any $2 \leftrightarrow 3$ log-branch transit
with even flattenings, the only ambiguity in Proposition
(\ref{fund_func}) is by multiplication by $N$th roots of unity, which
is due to the definition of the function $g$ in (\ref{funcg}). The
conclusion follows as before by using QHG-isomorphisms preserving the
parity of flattenings (whence based on even multiples of the
generators of flat/charge lattices).\hfill$\Box$
\begin{remark}\label{scissors}(Cheeger-Chern-Simons invariants and $\Hh_1$.)
{\rm By the results of \cite{N}, there is an injective
homomorphism from $H_3(PSL(2,\mc);\mz)$ (discrete homology) to a
scissors congruence group $\hat{\mathcal{P}}(\mc)$, such that the
dilogarithm (\ref{symRmatdil}), defined on
$\hat{\mathcal{P}}(\mc)$, restricts to the universal
Cheeger-Chern-Simons class $\hat{C}_2: H_3(PSL(2,\mc);\mz)
\rightarrow \mc/\pi^2\mz$. Hence $\Hh_1$ is a natural extension of
$\exp(\hat{C}_2/i\pi)$, the exponential of a constant times ${\rm
Vol} +iCS$, to classes representing QHFT bordisms. Recently J.
Dupont and C. Zickert produced dilogarithmic formulas for the lift
$\hat{C}_2'$ of $\hat{C}_2$ to $H_3(SL(2,\mc);\mz)$, such that
$\exp(\hat{C}_2'/i\pi)$ coincides with the lift of $\Hh_1$
determined by even flattenings in Theorem \ref{main} \cite{DZ}. We are
indebted to their work for pointing out the existence of such
flattenings.}
\end{remark}
Recall Definition \ref{composition}. Assume that $\Bb'' =
\Bb'*\Bb$ exists, and let $\Tt(\Bb)$ and $\Tt(\Bb')$ be given
weights $(h,k)$ and $(h',k')$, respectively. Then
$\Tt(\Bb')*\Tt(\Bb)$ is a distinguished flat/charged
$\Ii$-triangulation $\Tt(\Bb'')$, with some weight $(h'', k'')$.
(It follows from the Mayer-Vietoris exact sequence
    for $(\Bb'',\Bb',\Bb)$ that even if $h=h'=0$, it can happen
    that $h''\ne 0$. However, if the glued part of the boundary is
    connected, or is a boundary in $\Bb * \Bb'$, then $h=h'=0$ implies
    $h''=0$.)
\begin{prop} \label{functQHFT}{\rm (Functoriality)} For any composition
$\Bb'' = \Bb'*\Bb$ of bordisms, $\Hh_N(\Bb'',h'',k'')$ coincides
with $\Hh_N(\Bb',h',k')\circ \Hh_N(\Bb,h,k)$ up to sign and
multiplication by $N$th roots of unity.
\end{prop}
This is a direct consequence of Theorem \ref{main}. Also, we prove
as in Proposition 4.29 of \cite{BB2}:
\begin{prop}\label{unitarity} {\rm (Polarity)}
Write $\bar{\Bb}$ for the QHFT bordism with opposite orientation and
  complex conjugate holonomy $\bar{\rho}$. Then
  $\Hh_N(\bar{\Bb},-h,-k)$ and $\Hh_N(\Bb,h,k)^*$, the adjoint for the
  hermitian structure of $V(\alpha_\pm)$, coincide up to sign and
  multiplication by $N$th roots of unity.
\end{prop}
In the proof of Theorem \ref{teoTOP} we have seen that the space
$\tilde{W}^+(T,b)$ of ($+$)-Log-$\Ii$-parameters over $W^+(T,b)$
(the disjoint union of spaces $W^+(T,b,\mathcal{C})$) is
isomorphic to $\tilde{\mc}^{-3\chi(F)+2p}$, where $\tilde{\mc}$ is
the universal cover of $\mc \setminus \{0\}$. Similarly, for any
admissible triangulation $\tau$ of a topological cusp with $n$
$2$-simplices, we have the analytic subspace $Def(\tau)$ of
$\widehat{\mc}^n$ made of the $n$-uples of log-branches for the
tetrahedra of $\Ii$-cusps with base triangulation $\tau$, where
$\widehat{\mc}$ is defined in Section \ref{buildbox}. Such
log-branches satisfy the compatibility relations $L_\Tt(e)=0$ at
interior edges.
\begin{defi}\label{defps}{\rm Let $X=(Y,L_\Ff,\alpha_\pm)$ be a marked
topological bordism with $e$-triangulated or admissibly
triangulated boundary components. The {\it phase space} of $X$ is
the (analytic) subset $Def(X)$ in the product of the spaces
$\tilde{W}^+(T,b)$ and $Def(\tau)$ over the components of
$\partial \tilde{Y}$ and $L$, determined by the family of
distinguished flat/charged $\Ii$-triangulations of QHFT bordisms
supported by $(Y,L_\Ff,\alpha_\pm)$.}
\end{defi}
When $Y$ has empty boundary, $Def(X)$ is a generalization of the well
known deformation space of hyperbolic structures supported by ideal
triangulations of $Y\setminus L$, introduced in \cite{Th}, and
recently studied in \cite{Choi} and \cite{Cha}. When $\tilde{Y}$ is
the mapping cylinder of a diffeomorphism $\phi$ of $F$, the amplitudes
of QHFT bordisms supported by $\tilde{Y}$ define a morphism of the
trivial vector bundle
\begin{equation}\label{statebundle}
E(F)^+=E(T,b)^+: \tilde{W}^+(T,b) \times V(T,b) \rightarrow
\tilde{W}^+(T,b). \end{equation}(Note that any two choices of
$e$-triangulations $(T_1,b_1)$ and $(T_2,b_2)$ of $F$ yield
isomorphic bundles $E(T_1,b_1)$ and $E(T_2,b_2)$, with
birationally equivalent bases.) The resulting mapping of sections
of $E(F)$, the {\it states of $F$}, are studied in \cite{B}.
\begin{prop}\label{phase-space}{\rm (Analyticity)} For every $N\geq 1$,
the amplitudes of QHFT bordisms supported by
$X=(Y,L_\Ff,\alpha_\pm)$ vary analytically with the boundary
structure in $Def(X)$, up to sign and multiplication by $N$th
roots of unity.
\end{prop}
This follows immediately from the fact that the matrix
dilogarithms are analytic, together with the fact that any
path in $Def(X)$ lifts to a path of log-branches via the relations
induced by Definition \ref{flat/charge-def} (4).

\subsection{QHFT variants}\label{variant} By varying the
bordism category we can vary the corresponding QFT.
\smallskip

{\bf QHFT$^0$:} Consider the bordism category supported by triples
$(Y,L,\rho)$, where $L$ is an non-empty {\it unframed} tangle in $Y$
and $\rho$ is a $PSL(2,\mc)$-character on the whole of $Y$ (ie. $\rho$
is trivial at the meridians of $L$). In fact, we restrict to
holonomies $\rho$ such that $(Y,L,\rho)$ admits $\Dd$-triangulations
that extend a topological branched ideal triangulation $(T',b')$ of
each boundary component, say $(S,V)$, and for which the link $L$ is
realized as a Hamiltonian subcomplex (hence with no $\Ii$-cusp). In
particular the objects of this bordism category incorporate the
idealization of (necessarily idealizable) cocycles on $(T',b')$, that
represent the restriction of $\rho$ to $S$.  The arguments of Theorem
\ref{main} can be easily adapted to produce tensors $H_N(\Bb,h,k)$
associated to such a bordism $\Bb$, and eventually the so called
QHFT$^0$ variant of quantum hyperbolic field theory.
\medskip

{\bf Fusion of QHFT and QHFT$^0$:} We can consider triples
$(Y,L_\Ff,L^0,\rho)$, where $L=L_\Ff\cup L^0$ is a tangle with a {\it
framed part} $L_\Ff$ and an {\it unframed} one $L^0$.  We also
stipulate that $\rho$ is trivial at each meridian of $L^0$.  For every
object support $(S,V)$, we have a partition $V= V_\Ff \cup V_0$ and we
use ``mixed'' triangulations that looks like an efficient one at
$p\in V_\Ff$ and like an ideal one at $p\in V_0$. A similar mixed
behaviour holds for the adapted $\Dd$-triangulations of such
bordisms. We eventually get tensors still denoted
$\Hh_N(\Bb,h,k)$ giving variants, still denoted QHFT, that extend both
the previous one ($L^0 = \emptyset$), and QHFT$^0$ ($L_\Ff = \emptyset$).
\medskip

{\bf QHFT$^e$:} Let $(Y,L_\Ff,L^0,\rho)$ be as above, and let us
specialize to $\rho$ that, as usual, are trivial at the meridians of
$L^0$, but are {\it not} trivial at the meridians of $L_\Ff$. Now we
use mixed triangulations of each object support $(S,V_\Ff\cup V_0)$
that look like an {\it economic} triangulation (see Remark
\ref{varpar}) at each $p\in V_\Ff$. Concerning the adapted
$\Dd$-triangulations, each component of $L_\Ff$ contributes the
hamiltonian subcomplex with just a copy of the parallel curve
specifying the framing (recall that by using ordinary efficient
triangulations, it contributed with {\it two} parallel curves). We get
tensors now denoted $\Hh^e_N(\Bb,h,k)$, and a variant denoted
QHFT$^e$.
\medskip

Of course, there are no deep structural differences between these
variants; nevertheless each one has its own interest (see also Section
\ref{CLOSED}).

\subsection{Mapping class group representations}\label{MODREP}
Fix $F$. Set $\tilde{Y}_F= F\times [-1,1]$, $\tilde{L}=\partial
F\times [-1,1]$ with trivial vertical framing, and let ${\rm
Mod}(g,r)$ be the mapping class group of $F$, that is, the group of
homotopy classes of orientation-preserving diffeomorphisms of $F$
fixing pointwise each boundary component. Given $\psi_{\pm}:(\pm F,
(T,b,\Ll)) \rightarrow F \times \{\pm 1\}$, put $\psi=
\psi_+^{-1}\psi_-$ and $[\psi]$ for the corresponding element in ${\rm
Mod}(g,r)$. Denote by $\tilde{W}_{[\psi]}$ the {\it mapping torus}
$(F\times [-1,1])/ (x,-1) \sim (\psi(x),1)$ of $\psi$, with tunnel
boundary $\tilde{L}_{[\psi]}$. Let $\rho$ be the conjugacy class of
$PSL(2,\C)$-valued representations of $\pi_1(F\times [-1,1])$
(identified with $\pi_1(F))$) associated to $\Ll$.  To simplify
notations, in all statements of this section we do not mention the
weights (we understand they are fixed).

\begin{lem}\label{CYLtensor} Up to sign and multiplication by
  $N$th roots of unity (denoted "$=_N$") we have :

\smallskip

(1) For any fixed $(T,b,\Ll)$ the amplitudes
$\Hh_N(\tilde{Y},\tilde{L}_F,\rho,(\pm F,(T,b,\Ll)), [\psi_\pm])$
depend only on $[\psi]$. We denote them $\Hh_N([\psi])$.
\smallskip

  (2) $\Hh_N([id])$ is the identity map from $E_N(\alpha_-)$ to
  $E_N(\alpha_+)$, and $\Hh_N([h_2])\circ \Hh_N([h_1]) =_N
  \Hh_N([h_2h_1])$. In particular, for any $[\psi] \in {\rm Mod}(g,r)$
  the QHFT tensor $\Hh_N([\psi])$ is invertible, with inverse
  $\Hh_N([\psi^{-1}])$, and for a homotopically $d$-periodic $[\psi]$
  the QHFT tensor $\Hh_N([\psi])$ is of finite order less than or
  equal to $d$.  \smallskip

  (3) If $\psi(\rho)$ coincides with $\rho$, then ${\rm
    Trace}\bigl(\Hh_N([\psi])\bigr)=_N\Hh_N(\tilde{W}_{[\psi]},
  \tilde{L}_{[\psi]},\rho).$
\end{lem}
\noindent{\it Proof.} Point (1) follows from Theorem \ref{main},
because $$\Hh_N(\tilde{Y},\tilde{L}_F,\rho,(\pm
F,(T,b,\Ll)),[\psi_\pm]) =_N \Hh_N(\tilde{Y},\tilde{L}_F,\rho,(\pm
F,(T,b,\Ll)),[id],[\psi])$$ (the homeomorphism $\psi_-^{-1} \times
{\rm id}$ sends the first mapping cylinder to the second). By
Proposition \ref{functQHFT} we have $\Hh_N^2([id]) =_N
\Hh_N([id])$, so $\Hh_N([id])$ is an idempotent. It is invertible
because the matrix dilogarithms are. Both facts imply the first
claim in (2). The rest is a direct consequence of Proposition
\ref{functQHFT} and formula (\ref{deftt}). \hfill $\Box$
\medskip

The arguments in the proof of Lemma \ref{CYLtensor} imply also
that, letting $[\psi] = [id]$ and $\rho$ fixed, the amplitudes of
any marking variation $(T,b,\Ll) \rightarrow (T_1,b_1,\Ll_1)$ are
invertible. Hence $\Hh_N([\psi])$ is conjugated to the tensor
$\Hh_N(\tilde{Y},\tilde{L}_F,\rho,(\pm F,(T_1,b_1,\Ll_1)),
[\psi_\pm])$. Moreover $\Hh_N(-\tilde{Y},\tilde{L}_F,\rho,(\mp
F,(T,b,\Ll)), [\psi_\pm])$, the amplitude with reversed
orientation, clearly coincides with $\Hh_N([\psi]^{-1})$. Using
Proposition \ref{unitarity} we deduce:
\begin{cor}\label{mcgrep} For any fixed $\rho\in \Rr(g,r)$, the homomorphisms
$\psi \mapsto \Hh_N([\psi])$ induce a conjugacy class of
  linear representations of ${\rm Mod}(g,r)$, well-defined up to sign and
  multiplication by $N$th roots of unity. For $SL(2,\mr)$-valued
  characters $\rho$ these representations are unitary.
\end{cor}

\subsection{Tunneling the ($+$)/($-$) states}\label{+-versus} We use
($+$)-Log-$\Ii$-parameters to define the QHFT because of the existence
of strongly idealizable cocycles on QHFT bordism triangulations, which
makes functoriality easy to check. Here we exhibit a family of tensors
correlating the ($\pm$)-Log-$\Ii$-parameters, thus recovering, in
particular, the direct and nice interpretation of boundary structures
having non trivial holonomy at the punctures in terms of pleated
hyperbolic surfaces (see Section \ref{Ipar}). These tensors are also
used in Section \ref{FIBERED}.
\smallskip

For any base surface $F$ with an $e$-triangulation $(T,b)$, let $Z(F)$
be the pseudo-manifold obtained by collapsing to a point each boundary
annulus of the cylinder $C(F)=F\times [-1,1]$.  Recall the bundle
$E(F)^{+}$ in (\ref{statebundle}), and consider similarly
$E(F)^{-}:\tilde{W}^-(T,b) \times V(T,b) \rightarrow
\tilde{W}^-(T,b)$. We have:
\begin{prop}\label{tunnel+-} There exists a canonical family $\mathcal{F}$
of flat/charged $\Ii$-triangulations covering a portion of $Z(F)$,
with invertible trace tensor $\Hh_N(\mathcal{F}): E(F)^{+}
\rightarrow E(F)^-$.
\end{prop}
\noindent {\it Proof.} Orient $C(F)$ so that $\pm F$ is identified
with $F\times \{\pm 1\}$. Let $P(T,b)$ be the cell decomposition
of $C(F)$ made by the prisms with base the $2$-simplices of $T$.
Orient all the ``vertical'' (ie. parallel to $[-1,1]$) edges of
$P(T,b)$ towards $+F$. For every abstract prism $P$, every
vertical boundary quadrilateral $R$ has both the two horizontal
and the two vertical edges endowed with parallel orientations. So
exactly one vertex of $R$ is a source (that belongs to $-F$), and
exactly one is a pit (that belongs to $+F$). Triangulate each $R$
by the oriented diagonal going from the source to the pit. Finally
extend the resulting triangulation of $\partial P$ to a
triangulation of $P$ made of $3$ tetrahedra, by taking the cone
from the $b$-first vertex of the bottom base triangle of $P$ (note
that no further vertices nor further edges have been introduced).
Repeating this for every prism, we get a branched distinguished
triangulation $(C(T,b),H)$ of $C(F)$, where the vertical edges
make the Hamiltonian tangle $H$. As in the proof of Theorem
\ref{teoTOP} there exists integral charges on $(C(T,b),H)$.

Let $F\times [-1,3]$ be triangulated by two adjacent copies of
$C(T,b)$, glued each to the other at $F\times \{+1\}$. For any $z \in
\Zz_I(T,b)^+$, consider the unique cocycle $C_0(z)$ on the composition
$C(T,b)*C(T,b)$ that: extends $z\cup \Psi(z)$, given on $(F\times
\{-1\})\cup (F\times \{3\})$; takes the value $P$ of (\ref{convert})
on each vertical edge contained in $F\times [1,2]$; takes the value
$1$ on each vertical edge contained in $F\times [-1,1]$. Perturb
$C_0(z)$ with a $0$-cochain $s$ that: takes the value $1$ on $(F\times
\{-1\})\cup (F\times \{3\})$; takes values in $PB^+(2,\C)$ at each
vertical boundary annulus; restricts to an idealizable cocycle on
$F\times [-1,1]$ and to a {\it maximally idealizable} cocycle on
$F\times [1,3]$ (see Definition \ref{triideal}). Finally, glue
$\Ii$-cusps to the idealization.  Note that the only non idealizable
tetrahedra are those in the star of a boundary edge of the
triangulation of $F\times \{3\}$.  Lemma \ref{canonicalflat} gives
flattenings for the idealizable tetrahedra.

Look at the ideal triangulation $(T',b')$ of $S\setminus V$
corresponding to the copy of the triangulation $(T,b)$ for the
boundary component $F\times \{3\}$. For every cochain $s$ as above and
every edge $e$ of $T'$ we have two complex numbers: the
($-$)-exponential $\Ii$-parameter $W^-(\Psi(z))(e)$, and, as in
(\ref{totw}), the total product $\Ww(e)$ of the cross-ratio moduli at
the edges of $T$ that enter the definition of
$W^-(\Psi(z))(e)$. Recall that there are two distinct such edges only when $e$
contributes to make a marked corner. It is possible to normalize $s$
so that for every edge $e$ of $T'$ we have $W^-(\Psi(z))(e) =
\Ww(e)^{-1}$. 

Varying the cocycle $z \in \Zz_I(T,b)^+$, this choice determines the
family $\mathcal{F}$ in the statement. By perturbing the initial
cocycle $z$ with $0$-cochains $t$ with values in $PSL(2,\mc)\setminus
PB^+(2,\mc)$, the same construction leads to families $\mathcal{F}_t$
of flat/charged $\Ii$-triangulations covering the whole of $F \times
[-1,3]$. 

Note that for suitable flat/charges ($f_1=0$ in (\ref{symRmatdil}) and
$c_1=0$ in (\ref{symqmatdil})) the matrix dilogarithms have
well-defined finite limits when the cross-ratio modulus $w_0
\rightarrow 0$. From the symmetry relations of the matrix dilogarithms
(see \cite{BB2}, Corollary 5.6), this is true more in general for any
degenerating sequence of $\Ii$-tetrahedra, that is when $w_0$ goes to
$0$, $1$ or $\infty$. Now, we can choose in a continuous way the
flat/charges of $\mathcal{F}_t$ so that they satisfy the above
constraints on the tetrahedra of $\mathcal{F}_t$ that become non
idealizable in $\mathcal{F}$, when $t \rightarrow {\rm id}$. Then
$\textstyle \Hh_N(\mathcal{F}):=\lim_t \Hh_N(\mathcal{F}_t)$
exists. As in Lemma \ref{CYLtensor} (2) we see that
$\Hh_N(\mathcal{F}_t)$ is invertible, with inverse
$\Hh_N(-\mathcal{F}_t)$. Since $\textstyle \Hh_N(\mathcal{F}) \circ
\Hh_N(-\mathcal{F})= \lim_t (\Hh_N(\mathcal{F}_t)\circ
\Hh_N(-\mathcal{F}_t))=_N{\rm id}$ (Proposition \ref{phase-space}), we
deduce that $\Hh_N(\mathcal{F})$ is invertible.  \hfill$\Box$
\begin{figure}[ht]
\begin{center}
\includegraphics[width=5cm]{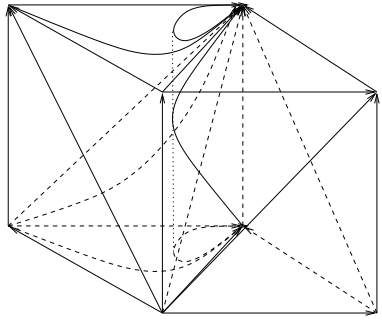}
\caption{\label{tcusp} Pasting opposite vertical sides yields an
instance of $C(T,b)$ for the once-punctured torus $S$, based on an
economic triangulation of $S$ as in Remark \ref{varpar}.}
\end{center}
\end{figure}

\section{Partition functions}\label{CLOSED}

Assume that $W$ is a {\it closed} oriented $3$-manifold, and that $L$
is a link in $W$ with a framed part $L_\Ff$ and an unframed one
$L^0$. Each variant of quantum hyperbolic field theory (see Section
\ref{variant}) leads to the respective {\it partition functions}. 

If $\rho$ is trivial at each meridian of $L^0$, we have the
QHFT {\it partition functions} 
$$ \Hh_N(W,L_\Ff,L^0,\rho,h,k)$$
that specialize to the QHFT$^0$ ones when $L=L^0$:
$$H_N(W,L,\rho,h,k)= \Hh_N(W,\emptyset,L^0,h,k) \ .$$
If $\rho$ is also assumed to be non trivial at each meridian of
$L_\Ff$, we have also
$$\Hh_N^e(W,L_\Ff,L^0,\rho,h,k) \ .$$ These partition functions are
scalars, well-defined up to sign and multiplication by $N$th-roots of
unity. Typical examples of triples $(W,L_\Ff,\rho)$ are given by
hyperbolic {\it cone} manifolds $W$ with framed cone locus $L_\Ff$ and
hyperbolic holonomy $\rho$ on $W \setminus L_\Ff$. The partition
functions can be expressed in terms of manifolds $Y$ with toric
boundary and containing an unframed link $L^0$ in the interior. By
fixing an ordered basis $(m_i,l_i)$ for the integral homology of each
boundary torus, let $W$ be obtained from $Y$ by Dehn filling along the
$m_i$, and $L_\Ff$ be the disjoint union of the cores of the filling
solid tori, framed by the $l_i$. Then the partition functions of
$(W,L_\Ff,L^0,\rho,h,k)$ are in fact invariants of
$(Y,\{(m_i,l_i)\}_i,L^0,\rho,h,k)$.

\subsection{QHFT vs QHFT$^0$ partition functions}\label{ALL}
For $\Bb=(W,L,\rho)$, $L=L^0$, with weights $h=k=0$,
$H_N(W,L,\rho,0,0)$ coincide with the invariants $H_N(W,L,\rho)$
constructed in \cite{BB1, BB2}.  Let us consider more generally
$(W,L_\Ff,L^0,\rho,0,0)$.  Fix also a framing $\Ff_0$ for
$L^0$. Then we can consider the partition function
$\Hh_N(W,L_\Ff\cup L^0_{\Ff_0},\emptyset,\rho,0,0)$.  Let us denote by
$\bar \lambda$ the unframed link obtained by splitting each component
of $L^0$ in the two corresponding parallel boundary components of the
ribbon link $L^0_{\Ff_0}$. We have:
\begin{prop}\label{QHFTQHI} 
$$ \Hh_N(W,L_\Ff\cup L^0_{\Ff_0},\emptyset,\rho,0,0) =_N
\Hh_N(W,L_\Ff, \bar \lambda,\rho,0,0) \ . $$ 
\end{prop}

\begin{figure}[ht]
\begin{center}
\includegraphics[width=3cm]{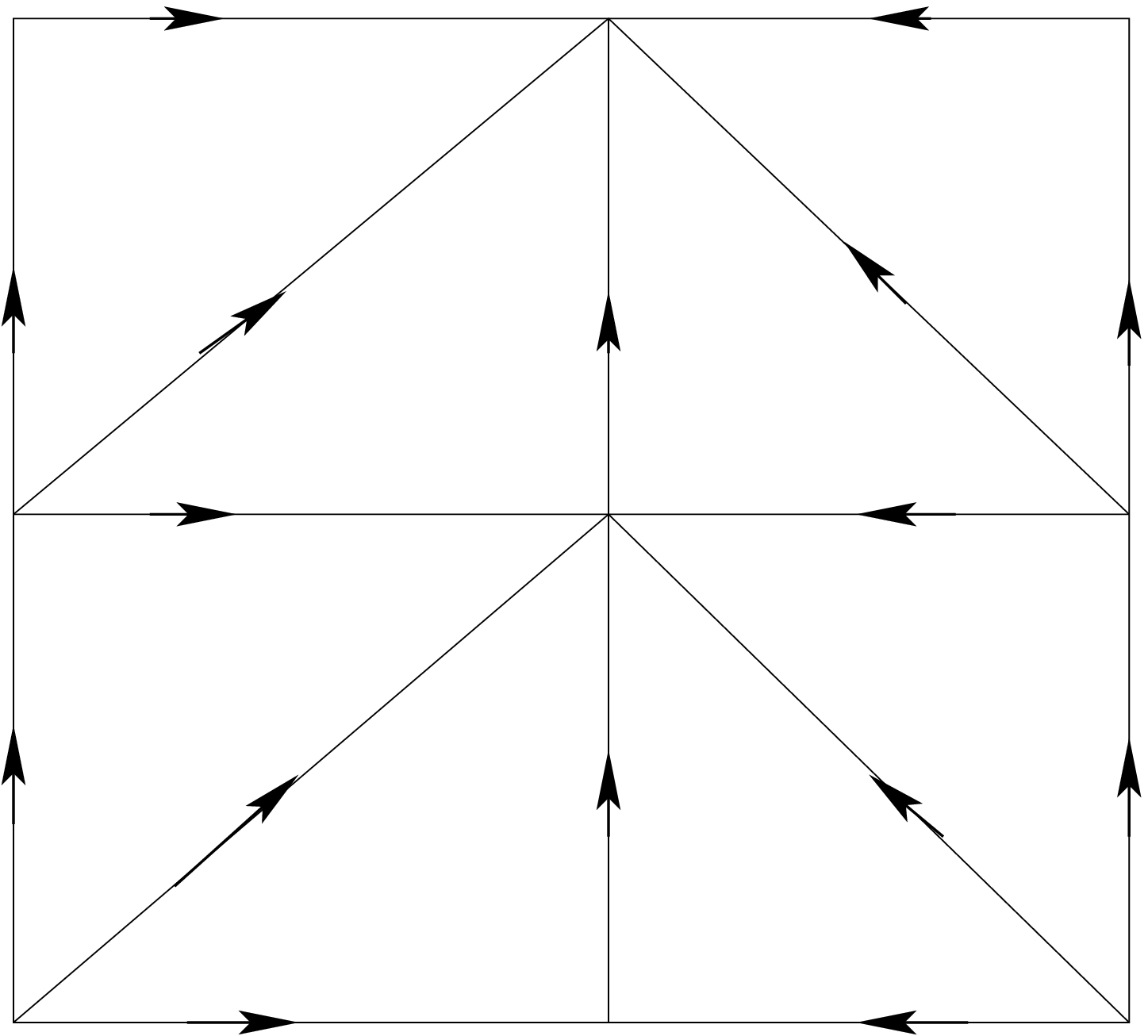}
\caption{\label{carre} A special admissible triangulation of $B$.}
\end{center}
\end{figure}

\noindent {\it Proof.} For simplicity, assume that $L=L^0$.  Fix a
$\Dd$-triangulation of $(\tilde{W},\tilde{L}_\Ff,\rho)$ where each
tunnel component $B$ has a symmetric admissible triangulation as in
Figure \ref{carre} (opposite sides of the quadrilateral are
identified). The tangle $\bar{\lambda}$ cuts open $B$ into symmetric
annuli, left and right to the central vertical line in Figure
\ref{carre}.

Because $\rho$ has trivial holonomy at the meridians of $L_\Ff$,
we can assume that the cocycle takes the same values on symmetric
edges. Identifying the annuli we thus get a $\Dd$-triangulation
for the QHFT$^0$ triple $(W,\bar{\lambda},\rho)$. Since
$H_N(W,\bar{\lambda},\rho)$ is computed from the idealization and
symmetric tetrahedra in the cusps have opposite branching
orientation, the result will follow if we show the existence of
symmetric flat/charges. Then each cusp tensor will be the identity
map.

The existence of flattenings with this property can be shown using
Remark \ref{canonicalflat}, but for charges we need to take
another route. Recall from the end of section \ref{QHG} that
flat/charges form affine spaces over an integral lattice generated
by vectors attached to the edges. For an edge of $B$, such vectors
can be represented as adding $+1$ at one of the left adjacent
corner and $-1$ at the other, and the inverse for the right
adjacent corners. Using these rules and (\ref{2-3comp}) it is
straightforward (though tedious) to check that any given
flat/charge can be turned into one with equal quantum log-branches
on symmetric tetrahedra.\hfill$\Box$

\subsection{Invariants of cusped manifolds and surgery formulas}\label{DEHNF}

Let us recall the QHG pseudo-manifold triangulations $\Tt$
used in \cite{BB2} (see Definition 6.2 and Definition 6.3 in that
paper) to define the quantum hyperbolic invariant trace tensors
$\Hh_N(\Tt)$ for oriented cusped hyperbolic manifolds.\smallskip

Let $M$ be a cusped manifold. Denote $Z$ the pseudo-manifold obtained
by taking the one point compactification of each cusp of $M$. $M$
admits a triangulation by positively embedded hyperbolic ideal
tetrahedra, possibly including some degenerate ones of null volume
(ie. having real cross-ratios). Such a triangulation can be obtained
by subdividing the canonical Epstein-Penner cell decomposition of
$M$. This gives rise to triangulation $(T_0,z_0)$ of $Z$, where $z_0$ is
the cross-ratio function of the abstract edges of $T$, the imaginary
part of every cross-ratio being $\geq 0$. We call it a {\it
  quasi-geometric} ideal triangulation of $Z$. If some quasi-geometric
triangulation admits a global branching, we say that $M$ is {\it
  gentle}. More generally, $M$ is said {\it weakly-gentle} if there is
an $\Ii$-triangulation $(T,b,w)$ of $Z$ such that $(T,z)$,
$z=w^{*_b}$, is obtained via a (possibly empty) finite sequence
$(T_0,z_0) \rightarrow \ldots \rightarrow (T_i,z_i) \rightarrow \ldots
\rightarrow (T,z)$ of positive $2 \rightarrow 3$ transits, where
$(T_0,z_0)$ is a quasi-geometric triangulation of $Z$ as above. Each
transit $(T_i,z_i) \rightarrow (T_{i+1},z_{i+1})$ is defined by
$W_{T_i}(e)=W_{T_{i+1}}(e)$, with all exponents $*_b=1$ (see
(\ref{totw})). Every such a $(T,b,w)$ can be enhanced to flat/charged
$\Ii$-triangulations $\Tt=(T,b,w,f,c)$.
\smallskip

Given $(T_0,z_0)$, it is certainly possible to get an
$\Ii$-triangulation $(T,b,w)$ by performing {\it also} some bubble
moves (hence introducing new interior vertices). The authors do not
know any example of non weakly-gentle cusped manifold, that is, such
that we are forced to do it. Anyway, dealing with bubble moves is a technical
difficulty which will appear also in the proof of Theorem
\ref{defscissors} (1) below. We overcome it as follows. We fix an edge
$a$ of the canonical Epstein-Penner cell decomposition of $M$, and
take $A$ made by two copies of $a$ that intersect at non manifold
points of $Z$; the second copy runs parallel to $a$ within an open
cell of the decomposition.  Hence, $A$ is a circle covered by two
arcs.  We need to enlarge the notions introduced in Definition
\ref{DISFLATCHARGE}. We say that {\it $(T,H)$ is a distinguished
  triangulation of $(Z,a)$} if $H$ is a subcomplex of the $1$-skeleton
of $T$ isotopic to $A$, that contains all the regular vertices of $T$,
and such that one arc of $A$ is covered by an edge $l$ of $H$. We say
that {\it $c$ is a global charge on $(T,H)$} if
\begin{equation}\label{condc2}C_{\Tt}(e)=\left
\lbrace\begin{array}{l} 4\\0\\ 2\end{array}\right.\begin{array}{l}
{\rm if}\ e=l\\{\rm if}\ e \subset H\setminus l\\ {\rm if}\ e \subset
T \setminus H \end{array}
\end{equation}
If $H=\emptyset$ this reduces to the usual notion of global charge on
a closed triangulated pseudo manifold whose non manifold points have
toric links. By using bubble moves and the existence of such usual
global charges, it is easily seen that $(T,H)$ supports global charges
as in (\ref{condc2}) (see the proof of Theorem 6.8 in \cite{BB2} for
the details).

We say that $\Tt=(T,H,b,w,f,c)$ is a flat/charged $\Ii$-triangulation
of $(Z,a)$ if $(T,H,b,c)$ is a branched, charged and distinguished
triangulation of $(Z,a)$, and $(T,z)$, $z=w^{*_b}$, is obtained from a
quasi geometric $(T_0,z_0)$ via a finite sequence $(T_0,z_0)
\rightarrow \ldots \rightarrow (T,z)$ of transits supported by
positive $2 \rightarrow 3$ moves {\it and bubble moves}. By setting $a
= \emptyset$ and $H=\emptyset$, this definition incorporates that for
the weakly-gentle case.
\smallskip

In \cite{BB2} it is shown that flat/charged $\Ii$-triangulations $\Tt$
of $(Z,a)$ (with arbitrary weights) do exist and that, for every odd
$N\geq 1$,
$$ H_N(M,a)=\Hh_N(\Tt)$$ is a well defined invariant of $(M,a)$,
providing the weights of flat/charges to be $0$. To simplify the
exposition, below we continue with this normalization.  When $M$ is
weakly-gentle we get invariants $H_N(M)$. In fact, as a by product of
the following discussion, we will realize that $H_N(M,a)$ does not
depend on the choice of $a$, so that $H_N(M)$ is always well defined
(see Corollary \ref{DF-converge}).

\medskip

Let us recall now a few facts related to {\it hyperbolic Dehn filling}
(see \cite{Th}, \cite{BP}, \cite{Pe-Po}). A quasi geometric
triangulation $(T_0,z_0)$ as above corresponds to the complete
structure of $M$. It can be deformed in a complex variety of dimension
equal to the number of cusps. If $z'$ is close enough to $z_0$,
$(T_0,z')$ is a triangulation by (possibly negative - see
\cite{Pe-Po}) embedded hyperbolic ideal tetrahedra in a non-complete
hyperbolic structure, say $M'$, close to $M$.  In some case the
completion of $M'$ gives rise to a compact closed hyperbolic manifold
$W$, topologically obtained by Dehn filling of the (truncated) cusps
of $M$. The core of each attached solid torus is a ``short'' simple
closed geodesic $L_j$ of $W$, so that we have the (geodesic) link
$\textstyle L=\coprod_j L_j$. Moreover, there are sequences
$(W^n,L^n)$ obtained in this way such that the length of $L^n$ goes to
$0$ when $n\to +\infty$. Hence $(W^n,L^n)$ converges to the cusped
manifold $M$ (in a neat geometric sense). From now on we will consider
small deformations $z'$ leading to such closed completions.
\smallskip

As well as $(T_0,z_0)$ gives rise to a triangulation $\Tt=(T,H,b,w,f,c)$
of $(Z,a)$, $z'$ close to $z_0$ gives rise to another
flat/charged $\Ii$-triangulation $\Tt'=(T,H,b,w',f',c)$, where $w'$ is
close to $w$ and the log-branch associated to the global flattening
$f'$ corresponds to a continuous deformation of the one for $f$.

\begin{lem}\label{deformedlogb}{\rm (See \cite{N}, p. 469) Let $z'$ be
a small deformation of $z_0$ producing $(W,L)$, and $m_j$ be a
meridian of each link component $L_j$. Then there exist flattenings
$f''$ for the deformed triangulation $(T,b,w')$ such that the weight
$\gamma(f'')$ associated to the collection of log-branches of
$\Tt''=(T,H,b,w',f'',c)$ satisfies $\gamma(f'')(m_j)=0$ for all $j$.}
\end{lem}
\noindent {\it Proof.} The logarithm of the derivative of the
holonomy of $m_j$ is $0$ at the complete structure, but after
deformation it represents a full
$2\pi$-rotation about $L_j$ (see eg. \cite{Th} or \cite{BP}).
Hence $\gamma(f')(m_j) \in \mz\sqrt{-1}\pi$. The result then follows from
the comments after Theorem \ref{theoweight}.\hfill$\Box$
\smallskip

An explicit construction of flattenings $f''$ as in the lemma shall be
recalled during the proof of theorem \ref{defscissors}. Hence, we
dispose of two flat/charged $\Ii$-triangulations of $(Z,a)$, $\Tt'=
(T,H,b,w',f',c)$ and $\Tt''=(T,H,b,w',f'',c)$, relative to a small
deformation $z'$ of the complete structure $z_0$ as above, leading to
respective trace tensors $\Hh_N(\Tt')$ and $\Hh_N(\Tt'')$. We can now
state the main results of this Section. 

\begin{teo}\label{defscissors}
{\rm[Cusped manifold surgery formula]} Let $(W,L)$ be obtained by
completion of a small deformation $z'$ of $z_0$, and $\Tt'$, $\Tt''$
be associated triangulations. Denote by $\rho$ the hyperbolic holonomy
of $W$. Then we have $H_N(W,L,\rho) =_N
\Hh_N(\Tt'')$. Moreover, associated to each cusp $C_j$ of $M$ there is
an explicitely known map $\Lambda_N^j(\Tt''): \{ N-{\rm states\ of\ T}
\}\rightarrow \mc$ such that the following {\rm surgery formula}
holds:
$$H_N(W,L,\rho) =_N \sum_s \prod_{\Delta \subset T}
    \Rr_N(\Delta,b,w',f',c)_s\ \prod_j\Lambda_N^j(\Tt'')(s)$$ where
    $s$ runs over the $N$-states of $T$ and
    $\Rr_N(\Delta,b,w',f',c)_s$ is the matrix dilogarithm entry
    determined by $s$, for the tetrahedron with the continuously
    deformed structure.
\end{teo}
\begin{cor}\label{DF-converge} 
If $\{(W_n,L_n,\rho_n)\}$ is a sequence of closed hyperbolic Dehn
fillings converging to the cusp manifold $M$, then for every arc $a$ we
have $\textstyle \lim_n H_N (W_n,L_n,\rho_n)=_N H_N(M,a)$. Hence
$H_N(M)=H_N(M,a)$ is always a well defined invariant of $M$ (beyond
the weakly-gentle case).
\end{cor}

\begin{remarks}\label{CS-surg}
{\rm (1) Theorem \ref{defscissors} is the analog for $N>1$ of Theorems
  14.7 and Theorem 14.5 in \cite{N}, which describe surgery formulas
  for the volume, Vol($W$), and Chern-simons invariant, CS($W$), of
  $W$:
$$\sqrt{-1}({\rm Vol}(W) + \sqrt{-1} {\rm CS}(W)) = \sum_{\Delta
  \subset T} \mathcal{R}(\Delta,b,w',f')-\frac{\pi\sqrt{-1}}{2}
\sum_j \lambda(L_j)$$ where $\Rr$ is given by (\ref{symRmatdil})
and $\lambda(L_j)$ is the complex length of $L_j$, that is, the logarithm of
the dilation factor of its holonomy, which is a loxodromic
transformation of $\mh^3$. The technical complications due to the bubble
moves disappear for $N=1$.
\smallskip

(2) If $M$ is gentle and has a geometric branched ideal
triangulation $(T,b,w)$ without degenerate tetrahedra, then for each
$3$-simplex the flattenings of $\Tt'$ for a sufficiently small
deformation are just $-*_b$ times integral charges. It follows from
the proof of Theorem \ref{defscissors} that the scissors congruence
class $\mathfrak{c}_{\Ii_{fc}}(W,L,\rho)$ of \cite{BB2}, section 7,
coincides with Neumann's deformed scissors congruence class
$\hat{\beta}(M')$ in \cite{N}, Theorem 14.7 (see also Remark 6.12 and
Conjecture 7.9 in \cite{BB2}, where the undeformed $\hat{\beta}(M)$ is
denoted $\mathfrak{c}_{\Ii_{fc}}(M)$).
\smallskip

(3) In general there are small deformations $z$ of $z_0$ leading to
complete manifolds that are still cusped, that is only some cusps of
$M$ have supported a hyperbolic Dehn filling. There are also sequences
of such cusped manifolds $M^n$, with (short) geodesic links $L^n$,
converging to $M$. Similarly to the fusion of QHFT with QHFT$^0$
(see Section \ref{variant}) we can define quantum hyperbolic
invariants $H_N(M^n,L^n)$ for which the natural extensions of Theorem
\ref{defscissors} and Corollary \ref{DF-converge} hold.}
\end{remarks}
Let us consider now $(W,L_\Ff,L^0,\rho)$. Let $L_j$ be a component
of $L_\Ff$, $\bar \lambda_j = L_j \cup L_j'$, $L'_j$ being
the longitude of $L_j$ specifying the framing. Let $U=U(L_j)$
be a tubular neighbourhooh of $L_j$ in $W$, and $l\subset \partial U$
be a non separating simple closed curve. Let $W(l)$ be obtained from
$W$ by the Dehn filling of $W\setminus {\rm Int}(U)$ along $l$. Denote by $l^*$
the core of the attached solid torus.
\begin{teo} \label{surcomp} {\rm
[Closed manifold surgery formula]} Assume that $\rho(l)={\rm id} \in
  PSL(2,\mc)$ and the weight $k$ satisfies $k([l])=0$. Denote: $\rho'$
  the natural extension of $\rho_{|W \setminus U}$ to $W(l)$;
  $\widetilde L_\Ff = L_\Ff \setminus L_j$; $k'$ the restriction of
  $k$ to $W(l)$. We have
\begin{equation}\label{eqsurg1}
\Hh_N(W,L_\Ff,L^0,\rho,0,k) =_N \ \Hh_N(W(l),\widetilde L_\Ff, L^0\cup l^* \cup
\bar{\lambda}_j,\rho',0,k') \ .
\end{equation}
If moreover $\rho$ is not trivial at the meridians of $L_\Ff$ and $l$
is a longitude of $\partial U$, then

\begin{equation}\label{eqsurg2}
\Hh_N^e(W,L_\Ff,L^0,\rho,0,k) =_N \Hh_N^e(W(l),\widetilde L_\Ff, L^0\cup l^*
\cup L'_j,\rho',0,k') \ .
\end{equation}
\end{teo}
Let us assume now that $L^0$ is made by $r$ parallel copies of $L_j$
along the ribbon $L_\Ff$ that encodes the framing. So denote $L^0$ by
$\lambda_r$; with this notation, $\bar \lambda = \lambda_2$. Assume furthemore
that $l=m$ is a meridian of $L_j$, so that $l^*=L_j$. By appling inductively
both (\ref{eqsurg1}) and Proposition \ref{QHFTQHI} to this situation we get
\begin{cor} For every $r\geq 2$ we have
$$H_N(W,\lambda_r,\rho)=_N H_N(W,\lambda_2,\rho) \ . $$
\end{cor}

\begin{remarks}\label{remsurg}{\rm  
\smallskip

(1) Though disjoint and complementary by hypothesis, formula
    (\ref{eqsurg2}) is formally the same as that of Proposition
    \ref{QHFTQHI}, when replacing $l$ by $m$.
\smallskip

(2) Assume (for simplicity) that $L=L_\Ff$. When $l$ is a longitude of
$L_j$, $l^*$ inherits a natural framing in $W(l)$. Hence we get a
triple $(W(l),\hat L_{\hat \Ff},\rho)$. It follows from the very
definition of the QHFT tensors that $\Hh_N(W,L_\Ff,\rho,0,0)=
\Hh_N(W(l),\hat L_{\hat \Ff},\rho,0,0)$ and the same with $\Hh_N^e$
(when defined) replacing $\Hh_N$.

\smallskip

(3) We have seen that both $\Hh_N$ and $\Hh_N^e$ partitions functions
display interesting features of QHG. A main advantage of the $\Hh_N$
ones is the possibility to set in a {\it same ``holomorphic family''}
the QHG tensors associated to characters that are both trivial and non
trivial at link meridians.  Consider for example a hyperbolic knot $L$
in $S^3$, endowed with the canonical framing $\Ff$. Kashaev's volume
conjecture concerns the asymptotic behaviour of $H_N(S^3,L,\rho_{{\rm
    triv}})$ when $N\to +\infty$.  A reasonable variant of it is in
terms of the partition functions $\Hh_N(S^3,L_\Ff,\rho_{{\rm triv}}
,0,0)=_N H_N(S^3,\bar \lambda,\rho_{{\rm triv}}) $.  A family as above
could be useful in order to establish connections with the $\Hh_N$
partition functions of $(S^3,L_\Ff,\rho_{{\rm hyp}})$, where
$\rho_{\rm hyp}$ is the hyperbolic holonomy of the cusped manifold $M=
S^3\setminus L$.}
\end{remarks}

The rest of the section is devoted to the proof of these results.
This goes in several steps.

\smallskip

A main tool is the {\it simplicial blowing up/down} procedure
considered by Neumann in \cite{N}, section 11. We use it just
to get a simplicial version of (topological) Dehn filling.  Let $Z$ be
a pseudo manifold without boundary such that every non-manifold point
has toric link. Let $v$ be a non-manifold point. Consider a closed
cone neighborhood $N(v)$ of $v$, and a non separating simple closed
curve $C$ on the torus  $\partial N(v)$. The topological Dehn
filling of $Z$ at $v$ along $C$ is the pseudo manifold $Z'$ obtained by
gluing a $2$-handle to $Z\setminus {\rm Int}(N(v))$ along $C$, and then
collapsing to one point the resulting boundary component. 
\smallskip

\begin{figure}[ht]
\begin{center}
\includegraphics[width=6cm]{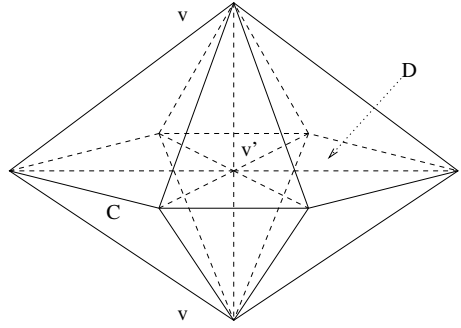}
\caption{\label{doublecone} An instance of double cone on a disc.}
\end{center}
\end{figure}

Now, let $T$ be a pseudo manifold triangulation of $Z$. Consider the
{\it abstract star} Star$^0$($v$) of $v$ in $T$. The boundary of
Star$^0$($v$) is the {\it abstract} link Link($v$) which is
homeomorphic to $\partial N(v)$. Assume that the curve $C$ is realized
as a {\it simplicial} curve on Link($v$). Then the cone from $v$ over
$C$ in Star$^0$($v$) is a triangulated disk $D^0$. The interior of
Star$^0$($v$) embeds onto the interior of the actual {\it star} of $v$
in $T$, Star($v$), which is made of the union of the $3$-simplices
having $v$ as a vertex. In this way $D^0$ maps onto a triangulated
singular disk $D$ in $Z$, that has embedded interior and singular
boundary immersed in the boundary of Star($v$). Cut open $T$ along
${\rm Int}(D)$ and glue the double cone $CD$ of $D$ (this is a
triangulated singular $3$-ball, see Figure \ref{doublecone}) so that
the top and the bottom get identified with the two copies of $D$
resulting from slicing. This gives a triangulation $T'$ of the
pseudo-manifold $Z'$ obtained by Dehn filling along $C$. It has the
property that every (abstract) tetrahedron of $T$ persists in
$T'$. Referring to the topological description, the interior of the
co-core of the $2$-handle attached to $Z\setminus {\rm Int}(N(v))$ is
isotopic to the interior of the union $H'$ of two edges, each joining
$v$ to the new vertex $v'$ at the ``center" of $CD$. In fact $H'$ is
the core of the solid torus added by the Dehn filling.

\begin{remark}\label{fewsimp}{\rm 
Note that in general there are very few simplicial curves on a given
Link($v$). Hence to get such a simplicial description of an arbitrary
Dehn filling, we will usually have to modify a given
triangulation. For the peculiar QHG pseudo-manifold triangulations
considered in this section, retriangulating will be possible by using
QHG isomorphisms, hence without altering the trace tensors. In fact,
any two triangulations of $\partial N(v)$ are connected by a finite
sequence of $2$-dimensional $1 \leftrightarrow 1$ ``flip'' moves (see
Figure \ref{eTflip}), and $1 \leftrightarrow 3$ moves obtained by
replacing a $2$-simplex with the cone of its boundary to a
point. Since $N(v)$ is homeomorphic to Link($v$), any such a sequence
is the boundary trace of a sequence of $2\leftrightarrow 3$ moves and
bubble moves in Star($v$). (Note, in particular, that $H$ passes
through the new interior vertices). By using the arguments of Theorem
6.8 in \cite{BB2}, we will always be able to choose that sequence so
that it lifts to a sequence of QHG transits.}
\end{remark}

Let us consider now a distinguished triangulation $(T,H)$ of $(Z,a)$
as above, and let $Z'$, triangulated by $T'$, be the result of a
simplicial Dehn filling of $Z$ along a curve $C$. Denote $H''$ the graph
union of the knot $H'$ (the core of the solid torus) and the image of
$H$ in $T'$. We define the notion of global charge on $(T',H'')$
by formally replacing $H$ by $H''$ in (\ref{condc2}).

\begin{lem}\label{extc} Let $c$ be a global charge on $(T,H)$
such that the charge weight of the curve $C$ is $0$. Then $c$ extends
to a global charge $c'$ on $(T',H')$.
\end{lem}
\noindent {\it Proof.} The complex $CD$ is made of pairs of
adjacent $3$-simplices, respectively above and below the disk $D$.
For a $3$-simplex of the top layer with charges $c_0$ and $c_1$ at
the edges in $D$ (ordered by using an orientation of $D$, say), we
will put the charges $-c_0$ and $2-c_1$ at these edges for the
symmetric $3$-simplex in the bottom layer. Then the other charges
are $c_2=1-c_0-c_1$ and $-c_2$, respectively. We have the charge
sum $C_{T'}(e)=2$ at each interior edge of $D$, and $C_{T'}(e) =
C_{T}(e)$ at the edges $e$ of $\partial D$. For the top edges $e'$
of $CD$ we can also choose the charges so that $C_{T'}(e')$ equals
$C_T(e)$, where $e$ is the copy of $e'$ in $D \subset M$. Indeed,
there are $n$ degrees of freedom in doing this, where $n$ is the
number of $1$-simplices in the curve $C$ used for blowing down. Then
we check that $C_{T'}(e')=2$ at the bottom edges. In particular,
the subcomplex $H$ survives in $T'$.

\begin{figure}[ht]
\begin{center}
\includegraphics[width=6cm]{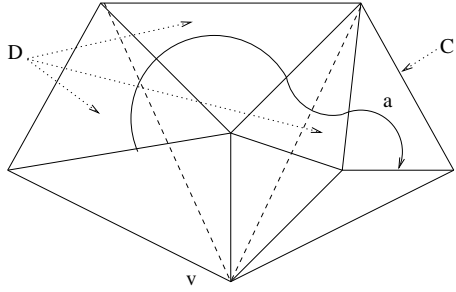}
\caption{\label{normalpath} A normal path $a$ in Link($v$) running
parallel to a blowing-up curve $C$. Four $3$-simplices of $T'
\setminus CD$ glued to three $2$-simplices of $D$ are shown in the
picture.}
\end{center}
\end{figure}

Note that $C_{T'}(e_0)=-C_{T'}(e_1)$ at the edges $e_0$ and $e_1$
of $H'$. We have to check that $C_{T'}(e_0)=0$, so that
(\ref{condc}) is satisfied on $(T',H''\setminus l)$. In fact,
$C_{T'}(e_0)$ is $n$ minus the sum of the $2n$ charges at the
bottom edges of $CD$, which is also the sum of charges in
$T'\setminus CD$ at these edges, minus $n$. We can form $n$ pairs
of such charges corresponding to the $3$-simplices of the ideal
triangulation $T$ of $M$ having a $2$-simplex in $D$. Replacing
for each of them the pair with $1$ minus the last charge, we get
that $C_{T'}(e_0)$ is equal to $\gamma(a)$, with $\gamma$ defined
in section \ref{QHG} and $a$ is a normal path in Link($v$) that
runs parallel to $C$ on one side (see Figure \ref{normalpath}).
Because the weight of $C$ is zero, we deduce $C_{T'}(e_0)=0$.

If $H \ne \emptyset$, we have to show that it can be deleted from
$H''$. As the two components $l$ and $H \setminus l$ are isotopic and
satisfy $C_{T'}(l)=4$ and $C_{T'}(e)=0$ for each edge $e \in H
\setminus l$, we can retriangulate the surgered pseudo-manifold $Z'$
so as to delete them, by using a sequence of charge transits starting
from $(T',H'')$ and terminating with a negative bubble move (see
\cite{BB1}, Proposition 4.27 and \cite{BB2}, proof of Theorem
6.8). Retriangulating $Z'$ backward, we eventually find a sequence of
charge transits terminating at $(T',H',c')$. The result
follows. \hfill$\Box$
\medskip

\noindent {\it Proof of Theorem \ref{defscissors}.}  To simplify
assume that $M$ has only one cusp. Take $\Tt''=(T,H,b,w',f'',c)$. For
the first claim we assume that the meridian $m$ of $L$ is a simplicial
path in Link($v$), where the vertex $v$ of $T$ corresponds to the
filled cusp. This is possible due to Remark \ref{fewsimp}.

Lemma \ref{extc} implies that $c'$ extends to $(T',H')$ after the Dehn
filling along $m$. Extend the branching $b$ by letting the new vertex
$v'$ be a pit of the double cone $CD$ we splice in $T$. By using Lemma
\ref{deformedlogb}, arguments similar to that of Lemma \ref{extc} show
that we can give the same log-branches on the $3$-simplices of $CD$,
in a pair above and below the disk $D$ (see \cite{N}, p. 454). Hence
we get a distinguished flat/charged $\Ii$-triangulation for
$(W,L,\rho)$. The weight $h \in H^1(W;\mz/2\mz)$ is clearly $0$
because of the epimorphism $ H_1(M;\mz/2\mz) \rightarrow
H_1(W;\mz/2\mz)$ induced by inclusion.  As in the first claim of Lemma
\ref{CYLtensor} (2) we see that the (unnormalized) trace tensor for
$CD$ is $N$ times the identity map from the linear space attached to
the top copy of $D$ to that for the bottom one. Combining this with
the normalization of trace tensors in (\ref{deftt}) gives
$H_N(W,L,\rho) = \Hh_N(\Tt'')$.

\begin{figure}[ht]
\begin{center}
\includegraphics[width=6cm]{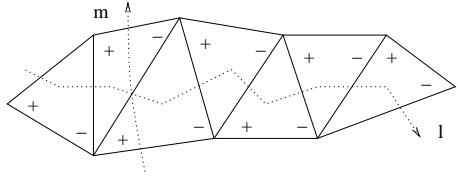}
\caption{\label{surgcurve} Flat/charge corrections for a Dehn
filling.}
\end{center}
\end{figure}

By Lemma \ref{deformedlogb} we know that $\gamma(f'')-\gamma(f') \in
H^1(\partial M;\mz\pi\sqrt{-1})$ is non zero only at the class of $m$,
where it is $-2\sqrt{-1}\pi$. Hence the collection of values of
$f''-f'$ determines a path $l$ normal to the cusp triangulation
induced by $T$, that intersects $m$ once and whose homology class is
Poincar\'e dual to $(\gamma(f'')-\gamma(f'))/2\sqrt{-1}\pi$. Denote
$$\Delta(l_j)=(\Delta^1,\ldots,\Delta^{\mid \Delta(l_j) \mid})$$ the
sequence of flat/charged $\Ii$-tetrahedra (possibly with repetitions)
determined by the $2$-simplices met by $l$. Each time $l$ goes through
a $2$-simplex it selects one of its vertices, whence a cross-ratio
modulus, say $z_i$, of the tetrahedron $\Delta^i$ corresponding to the
$2$-simplex. The values of $f''$ on $\Delta^i$ are obtained from those
of $f'$ by adding or substracting $1$ at the edges corresponding to
the other two vertices, as indicated in Fig. \ref{surgcurve}. For any
fixed tetrahedron $\Delta$ of $T$ all three flattenings may be
eventually altered, and/or differ from those of $f'$ by adding or
substrating $n_i \in \mz$ with $n_i\ne-1$ or $1$, exactly when
$\Delta=\Delta^i=\Delta^j$ for some $i\ne j$. Now, recall from
(\ref{deftt}) that
$$H_N(M,a) = \Hh_N(\Tt') = \sum_s \prod_{\Delta \subset T}
\Rr_N(\Delta,b,w,f,c)_s.$$ Put $\zeta=\exp(2\pi\sqrt{-1}/N)$. For any
$x\in \mc \setminus \{\zeta^j,j=1,\ldots,N-1\}$ the function $g$
defined in the Appendix satisfies (\cite{BB2}, Lemma
8.2): $$g(x\zeta^k)=g(x)\prod_{j=1}^k
\frac{(1-x^N)^{1/N}}{1-x\zeta^j}.$$ Using this formula, it is easily
checked that given a flat/charged $\Ii$-tetrahedron
$(\Delta,b,w',f',c)$ with $f'=(f_0',f_1',f_2')$ and positive branching
orientation, if $f''=(f_0'+n,f_1',f_2'-n)$ we have
\begin{equation}\label{defqR} \Rr_N(\Delta,b,w',f'',c)_s
=\Rr_N(\Delta,b,w',f',c)_s\
\prod_{j=1}^n\frac{(w_1')^{-1}}{1-w_0'\zeta^{i-k+j}}
\end{equation}
up to multiplication by $N$th roots of unity, where $i$ and $k$ are as
in (\ref{symqmatdil}). For each $2$-simplex met by $l$ we can apply
(\ref{defqR}) to the corresponding tetrahedron, or the similar formula
(deduced from Corollary 5.6 of \cite{BB2}) for any other
branching. This defines the function $\Lambda_N(\Tt'')$, so that we get
$$\begin{array}{lll} \Hh_N(\Tt'') & = &
  \sum_s \prod_{\Delta \subset T} \Rr_N(\Delta,b,w',f'',c)_s\\ & = &
  \sum_s \prod_{\Delta \subset T}
  \Rr_N(\Delta,b,w',f',c)_s\ \Lambda_N(\Tt'').
\end{array}$$
The conclusion follows from the equality $H_N(W,L,\rho) =
\Hh_N(\Tt'')$.\hfill$\Box$

{\it Proof of Theorem \ref{surcomp}.} Again for simplicity, assume
that $L$ is a {\it knot} (one component).  We apply the very same
arguments as for the first claim of Theorem \ref{defscissors}. In
particular, Lemma \ref{extc} applies verbatim. Since $\rho(l)$ is
trivial, for an arbitrary flattening the weight along $l$ (computed in
the flattened $\Ii$-cusps) lies in $\mz\pi\sqrt{-1}$. Hence we can
again use Theorem \ref{theoweight} to deduce the existence of
flattenings with zero weight along $l$. Then we give the same
log-branches on the $3$-simplices of the singular $3$-ball $CD$, in a
pair above and below the disk $D$ (see \cite{N}, p. 454). Note that if
we use $\Dd$ triangulations leading to $\Hh_N$ partition functions,
then both parallel components $L,L'$ that make $\bar \lambda$ survive
in the Hamiltonian subcomplex. If we can deal with $\Hh_N^e$-ones,
only the framing longitude $L'$ survives.  Hence we eventually get a
distinguished flat/charged $\Ii$-triangulation for $(W(l),l^* \cup
\bar{\lambda},\rho)$, or $(W(l),l^* \cup
\lambda',\rho)$ respectively. \hfill$\Box$

\subsection{Manifolds that fiber over $S^1$ - Examples}\label{FIBERED} 
\begin{figure}[ht]
\begin{center}
\includegraphics[width=6cm]{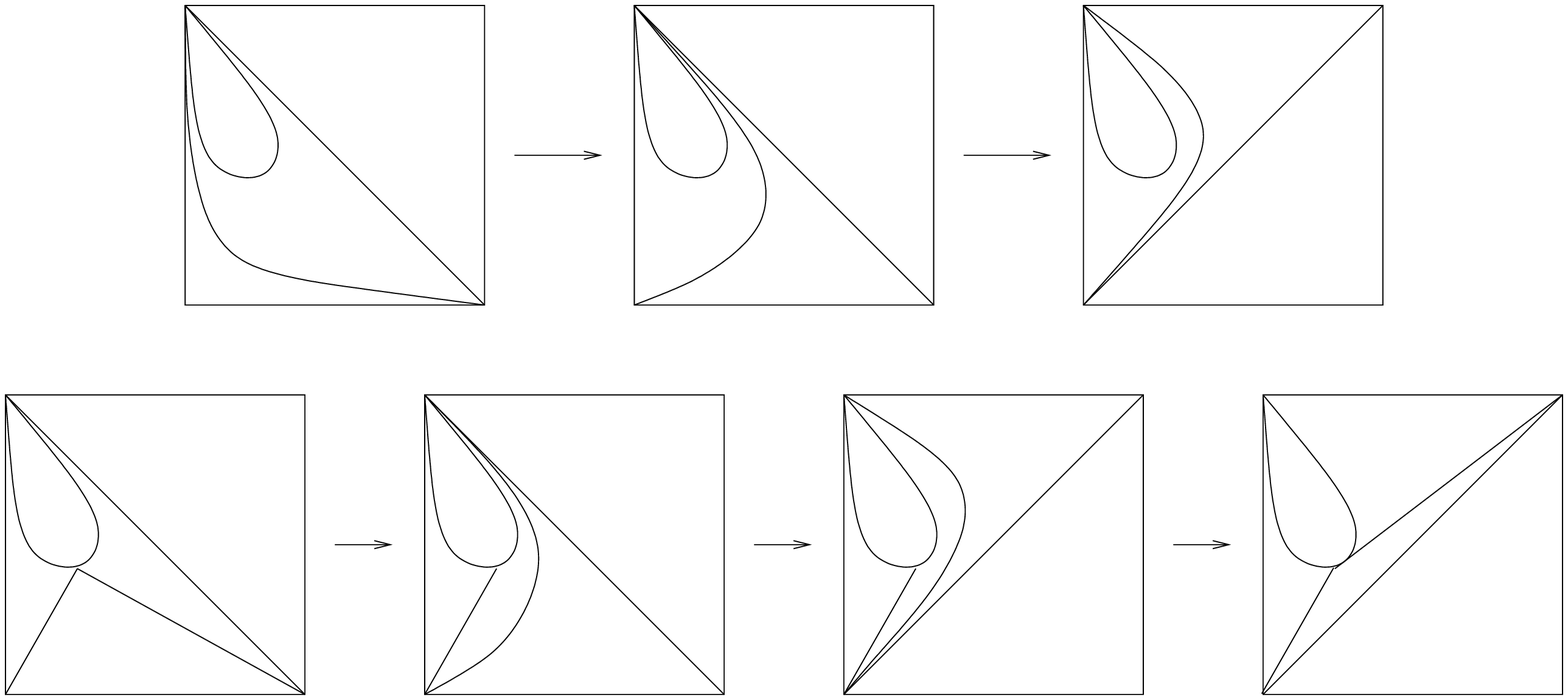}
\caption{\label{liftflip} Lifts to $e$-triangulations (economic ones -
  see Remark \ref{varpar} - at the first row) of flip moves on the
  corresponding ideal triangulations near marked corners. The
  tetrahedron associated to the first flip (first and third ones for
  the second row) degenerates for a sequence $s$ with
  ($-$)-Log-$\Ii$-parameters.}
\end{center}
\end{figure}

Lemma
\ref{CYLtensor} (3) gives a practical recipe to compute the QHFT
partition functions of mapping tori. A specific class of
distinguished flat/charged $\Ii$-triangulations of
$(\tilde{W}_{[\psi]},\tilde{L}_{[\psi]},\rho)$ is obtained by
composing one for the trivial mapping cylinder $F \times [-1,1]$, say
$\Tt_{\rm triv}$, with the monodromy action on the $e$-triangulation
$(T,b)$ of $F \times \{1\}$, and then gluing the two boundary
components. The monodromy action can always be decomposed as a
sequence of flip moves: a single flip on the ideal triangulation $T'$
associated to $T$ defines a flip on $(T,b)$ if it is not adjacent to a
marked corner, and it lifts to sequences as in Figure \ref{liftflip}
otherwise. We view these sequences as the result of gluing
tetrahedra. Hence the monodromy action determines a branched
triangulated pseudo-manifold $T_s$. This can be completed with global
charges, and, as for any $\rho$ we are free to choose the cocycle in
$\Tt_{\rm triv}$, we can also complete $T_s$ to a flattened
$\Ii$-triangulation.  Equivalently we can define a sequence
\begin{equation}\label{mon1}
s: (T,b,\Ll) \rightarrow \ldots \rightarrow
\psi(T,b,\Ll)\end{equation} of $e$-triangulations with
($+$)-Log-$\Ii$-parameters compatible with $\rho$. Note that the
edges of the associated pattern $\Tt_{s}$ of flat/charged
$\Ii$-tetrahedra are disjoint from the Hamiltonian link $H$. Since
$\Hh_N(\Tt_{\rm triv})=_N {\rm id}$, we deduce that
$\Hh_N(\tilde{W}_{[\psi]},\tilde{L}_{[\psi]},\rho)=_N
\Hh_N(\Tt_{s})$. 

Using a similar construction we now prove the
relationship with the quantum hyperbolic invariants $H_N$ of
fibered cusped manifolds \cite{BB2}. Recall that $W\setminus L$ is
homeomorphic to ${\rm Int}(\tilde{W}_{[\psi]})$. Denote by $l$ the
number of components of $L$.
\begin{prop} \label{teofibered} 
If ${\rm Int}(\tilde{W}_{[\psi]})$ supports a (necessarily unique)
complete hyperbolic structure with holonomy $\rho_c$, then
$\Hh_N(\tilde{W}_{[\psi]},\tilde{L}_{[\psi]},\rho_c,0,0)=_N
N^{2l}\ H_N(W\setminus L)$.
\end{prop}

\noindent {\it Proof.} Let $S\setminus V$ be the fiber of $W\setminus
L \rightarrow S^1$, equipped with an ideal triangulation $T'$. First
we show the existence of sequences $s:T'\rightarrow \ldots \rightarrow
\psi(T')$ of flip moves decomposing the monodromy action, such that
the associated pseudo-manifolds $T'_s$ are topological ideal
triangulations of $W\setminus L$, which moreover have maximal
volume. 

The first condition follows from the fact that the monodromy is
homotopically aperiodic (ie. pseudo-Anosov), so that $T'_s$ is
genuinely three-dimensional. When $S\setminus V$ is a once-punctured
torus, the second condition is a consequence of a result of Lackenby
\cite{La}, showing that the {\it monodromy ideal triangulation} of
Floyd and Hatcher \cite{FH} is isotopic to the canonical
Epstein-Penner cellulation. More in general, since no edge of $T'_s$
is homotopically trivial, the results of \cite{F1} imply that we can
straighten the tetrahedra to oriented geodesic ones, possibly with
overlappings, so that the algebraic sum of volumes is Vol($W\setminus
L$). This is known to be maximal \cite{F2}.

As in Section \ref{DEHNF} we can complete $T_s'$ to a flat/charged
$\Ii$-triangulation $\Tt'$. Hence the invariants $H_N(W\setminus L)$
can be computed as trace tensors $\Hh_N(\Tt')$. We note that in the
case when there are several fibrations of $W\setminus L$, or $T_s'$ is
not canonical, the invariance follows from Theorem 6.8 (2) in
\cite{BB2}, which shows that any two flat/charged $\Ii$-triangulations
of $W\setminus L$ with maximal volume are QHG-isomorphic.

Let us denote $\Tt_{s}'$ the result of cutting $\Tt'$ along the
fiber. The two boundary copies are marked pleated hyperbolic surfaces
$(T',b',\mathcal{L}) \rightarrow S\setminus V$ and
$(T',b',\mathcal{L}) \rightarrow \psi(S\setminus V)$, with shear-bend
coordinates (ie. $(-)$-Log-$\Ii$-parameters) $\Ll$ that determine
completely the log-branches of $\Tt_s'$. Recall from Section
\ref{+-versus} the families of flat/charged $\Ii$-triangulations
$\mathcal{F}_t$ and $\mathcal{F}$, and let $C \in \mathcal{F}$, $C_t
\in \mathcal{F}_t$ have boundary structures $(T,b,\mathcal{L})$,
associated to $(T',b',\mathcal{L})$, and $(T,b,\mathcal{L}_t)$,
respectively, at $F\times \{1\}$. We have: $$\begin{array}{lll}
\Hh_N(\Tt') & = & {\rm Tr}(\Hh_N(\Tt_{s}'))\\ & = & N^{-2l}\ {\rm
Tr}\left( \Hh_N(\Tt_{s}') \otimes {\rm id}^{\otimes 2l} \right)\\ & =
& N^{-2l}\ {\rm Tr}\left(\Hh_N(C) \circ (\Hh_N(\Tt_{s}') \otimes {\rm
id}^{\otimes 2l}) \circ \Hh_N(\psi(C))^{-1}\right)\\ & = & N^{-2l}\
\lim_{t \rightarrow {\rm id}} {\rm Tr}\left(\Hh_N(C_t) \circ
(\Hh_N(\Tt_{s,t}')\otimes {\rm id}^{\otimes 2l})\circ
\Hh_N(\psi(C_t))^{-1}\right)\\ & = & N^{-2l}\ \lim_{t \rightarrow {\rm
id}} {\rm Tr}(\Hh_N(C_t * \Tt_{s,t}'*(-\psi(C_t)))).
\end{array}$$
Here we use the invertibility of $\Hh_N(C)$ (Proposition
\ref{tunnel+-}) and the equality $\Hh_N(C) =_N \Hh_N(\psi(C))$. We
define $\Tt_{s,t}'$ as the continuous deformation of $\Tt_{s}'$
obtained from the sequence
$$s_t: (T,b,\Ll_t) \rightarrow \ldots \rightarrow \psi(T,b,\Ll_t)$$
similarly as in (\ref{mon1}) (see also Figure \ref{liftflip}). In the
last equality, $C_t * \Tt_{s,t}'*(-\psi(C_t))$ is for any fixed $t$ a
distinguished flat/charged $\Ii$-triangulation of the mapping cylinder
of $\psi$. Hence $\Hh_N(C_t * \Tt_{s,t}'*(-\psi(C_t)))$ does not
depend on $t$ up to conjugacy, and we conclude with Lemma
\ref{CYLtensor} (3).\hfill$\Box$

\medskip

By following the above computation backwards, we see more in general
that for {\it any} $PSL(2,\mc)$-character $\rho$ that can be realized
by ($-$)-Log-$\Ii$-parameters on some ideal triangulation of the fiber
$S\setminus V$, we have
\begin{equation}\label{mon2}
\Hh_N(\tilde{W}_{[\psi]},\tilde{L}_{[\psi]},\rho,h,k) = N^{2l}\ {\rm
Tr}( \Hh_N(\mathcal{T}_s'))
\end{equation}
where $\mathcal{T}_s'$ is a pattern of flat/charged $\Ii$-tetrahedra
associated to a sequence similar to $(\ref{mon1})$, but with
$e$-triangulations equipped with $(-)$-Log-$\Ii$-parameters compatible
with $\rho$.
\begin{remark}\label{QTeich}{\rm The formula $H_N(W\setminus
L) = {\rm Tr}(\Hh_N(\Tt_{s}'))$ expresses the quantum hyperbolic
invariants of fibered cusped manifolds as amplitudes between two
markings of the fiber, identified with a pleated hyperbolic
surface. For a similar construction based on representations of
quantum Teichmuller spaces, see \cite{BoL}.}
\end{remark}
{\bf Example: the figure-eight knot complement}\label{EXAMPLES}
Here we compute the QHFT partition functions of $(S^3,K_0)$, where
$K_0$ is the $0$-framed figure-eight knot in $S^3$. Recall that
$S^3\setminus K$ is fibered over $S^1$, with fiber the once-punctured
torus $\Sigma_{1,1}$; the $0$-framing of $K$ is induced by the
fibration. For simplicity, below we consider only characters $\rho$ of
{\it injective} representations $\pi_1(S^3\setminus K)\rightarrow
PSL(2,\mc)$. The restriction to $\Sigma_{1,1}$ of such representations
can be realized by $(-)$-Log-$\Ii$-parameters on any ideal
triangulation of $\Sigma_{1,1}$, so that, by (\ref{mon2}), we can
determine the corresponding subspace (still denoted $Def(S^3,K)$) of
the phase space of Definition \ref{defps} by using the monodromy ideal
triangulation.
\smallskip

\begin{figure}[ht]
\begin{center}
\includegraphics[width=10cm]{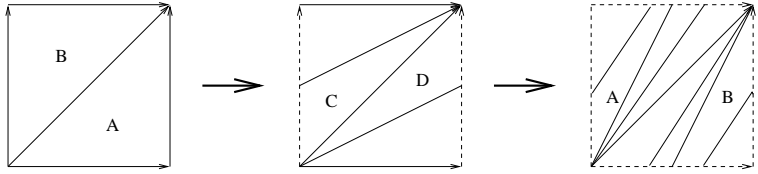}
\caption{\label{Tdyn} The composition of flip transformations
  for the monodromy of $S^3\setminus K$.}
\end{center}
\end{figure}

 The monodromy $\Phi:\Sigma_{1,1}\rightarrow \Sigma_{1,1}$ of
 $S^3\setminus K$ is isotopic to the hyperbolic element
$$\left(  \begin{array}{cc} 2 & 1 \\ 1 & 1  \end{array} \right) =
\left( \begin{array}{cc} 1 & 1 \\ 0 & 1 \end{array}\right) \left(
  \begin{array}{cc} 1 & 0 \\ 1 & 1 \end{array} \right) \quad \in
SL(2,\mz).$$ This description of $[\Phi] \in {\rm Mod}(\Sigma_{1,1})$
can be understood in terms of the Diagram of $PSL(2,\mz)$ (see
eg. \cite{FH}) via the action of $\Phi$ on topological ideal
triangulations of $\Sigma_{1,1}$, which can be represented by two flip
moves. See Figure \ref{Tdyn}, where the left picture is a lift to
$\mr^2 \setminus \mz^2$ of such a triangulation.

\begin{figure}[ht]
\begin{center}
\includegraphics[width=8cm]{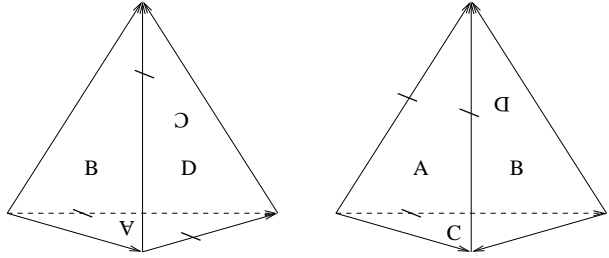}
\caption{\label{Teight} The face and edge identifications for the canonical
  geodesic ideal triangulation of $S^3\setminus K$.}
\end{center}
\end{figure}

The monodromy ideal triangulation $T$ of $S^3 \setminus K$ is obtained
by realizing each flip move via the gluing of an ideal tetrahedron,
first on a fixed triangulation of $\Sigma_{1,1}$, then on the
resulting one. The remaining four free faces are identified under
$\Phi$. It is not difficult to see that $T$ is isotopic to the
canonical geodesic ideal triangulation of $S^3 \setminus K$ with its
complete hyperbolic structure. The gluing pattern of the tetrahedra in
$T$ is shown in Figure \ref{Teight}.

It is well-known (see \cite{Th}) that the deformation space of smooth
(non necessarily complete) hyperbolic structures on $S^3 \setminus K$
is isomorphic to the algebraic set $Def_{{\rm hyp}}(S^3 \setminus K)
\subset \mathbb{H}^2 \times \mathbb{H}^2$ of points $(w_2,z_0)$ such
that
$$w_2 \in \mathbb{H}^2 \setminus \left\lbrace \frac{1}{2} +
\frac{t}{2}i \mid t \geq \sqrt{15} \right\rbrace,\quad z_0 =
\frac{1}{2} + \left(\frac{1}{4} +
\frac{1}{w_2(w_2-1)}\right)^{\frac{1}{2}}$$ where $w_2$ is the
cross-ratio modulus of the edge $e_2$ in the tetrahedron $\Delta^+$
with positive branching orientation (back edge in the left tetrahedron
of Figure \ref{Teight}), and similarly for $z_0$ in $\Delta^-$. The
space $Def_{{\rm hyp}}(S^3 \setminus K)$ is a subspace of 
\begin{equation}\label{subs}
\mathcal{C}=\{
(w_2,z_0) \in (\mc \setminus \{0,1\})^2 \mid
w_1w_2^2z_0^{-2}z_1^{-1}=1\},
\end{equation} which is isomorphic to the whole set of solutions of the edge
compatibility relations for cross-ratio moduli (see Definition
\ref{non-rough}). By an easy computation we find that the edge
compatibility relations for log-branches and charges are:
\begin{equation}\label{eqfc}
(S)\left\lbrace\begin{array}{l}f_1^-+2f_0^-+2f_0^++f_1^+=
(\arg(w_1)+2\arg(w_2)-\arg(z_1)-2\arg(z_0))-2*_+\\
c_1^-+2c_0^-+2c_0^++c_1^+=0\end{array}\right.
\end{equation}
where $*_+$ is the sign of the imaginary part of the $w_i$, and the
flattenings $f_i^+$ and $f_i^-$ correspond to the cross-ratio moduli
$w_i$ and $z_i$ of $\Delta^+$ and $\Delta^-$, respectively. Hence, by
(\ref{qlb}) we see that the QHFT phase space $Def(S^3,K)$ is the
covering of $\mathcal{C}$ given by
$$Def(S^3,K)=\{((w_2;f_2^+-c_2^+,f_0^+-c_0^+),(z_0;f_0^-+c_0^-,f_1^-+c_1^-))
\in \hat{\mc}\times\hat{\mc} \mid (w_2,z_0) \in \mathcal{C},\ (S){\rm\
is\ satisfied}\ \}.$$ Let us point out some geometrically meaningful
subspaces (compare with \cite{N}, Section 15). The dilation factors of
the standard meridian $m$ and longitude $l$ of $K$ are
$$\mu(m) = z_2w_2^{-1},\quad \mu(l)=z_0^2z_2^{-2}.$$ The complete
hyperbolic structure of $S^3\setminus K$ is obtained by solving
$\mu(m)=1$ and $\mu(l)=1$ for points of $Def_{{\rm hyp}}(S^3 \setminus
K)$, which gives $z_j=(w_j)^* = \exp(i\pi/3)$. We already noticed that
the weights $k$ for this solution take values in $2\mz$ (see \cite{N},
Proposition 5.2). Corresponding flattenings satisfy
\begin{equation}\label{eqcompfc}f_2^--f_2^+ = k(m),\quad2f_0^--2f_2^-=k(l)
\end{equation}
and (\ref{eqfc}) simplifies to
\begin{equation}\label{simpfs}
f_1^-+2f_0^-+2f_0^++f_1^+=0.
\end{equation}
Together with (\ref{eqcompfc}) this gives
\begin{equation}\label{standardf} f_1^+ = \frac{k(l)}
{2}-1-2f_0^+,\quad \left\lbrace \begin{array}{l} f_0^-=k(m)-f_0^+\\
    f_1^- = 2k(m)-\frac{k(l)}{2}+1+2f_0^+\end{array}\right.
\end{equation}
Similarly, integral charges are obtained by letting $c_i^+ = -f_i^+$
on $\Delta^+$ and $c_i^- = +f_i^-$ on $\Delta^-$. It can be checked
that $h=0$. Of course, (\ref{standardf}) is still true for points of
$Def(S^3,K)$ sufficiently near the complete solution, if now $k(m)$
and $k(l)$ are replaced with
$\tilde{k}(m)=k(m)-\log(\mu(m))/\sqrt{-1}\pi$ and
$\tilde{k}(l)=k(l)-\log(\mu(l))/\sqrt{-1}\pi$. Hence, in the vicinity
of the complete structure we have the strata
$$\{((w_2;2f_0^+-\tilde{k}(l),2f_0^+),(z_0;2\tilde{k}(m)-2f_0^+,
2(2\tilde{k}(m)+1+2f_0^+)-\tilde{k}(l)),\ f_0^+ \in \mz\},$$
parametrized by the lifts $k(m)$ and $k(l)$ of $\mu(m)$ and $\mu(l)$.
\smallskip

For points of $\mathcal{C}$ corresponding to hyperbolic structures
whose metric completion is obtained from $S^3 \setminus K$ by $(p,q)$
Dehn filling, we have $p\mu(m) +q\mu(l) = 2\pi \sqrt{-1}$. Hence the
flattenings $f=f''$ of Lemma \ref{deformedlogb} are to satisfy
$$p(f_2^--f_2^+)+q(2f_0^--2f_2^-)=-2$$ Let us fix $r$, $s \in \mz$
such that $ps-qr=1$. Solving simultaneously (\ref{simpfs}) and the
last equation gives
\begin{equation}\label{surgf}
f_1^+ = r-1-2f_0^+,\quad \left\lbrace \begin{array}{l}
    f_0^-=-2s-f_0^+\\ f_1^- = -r+4s+1+2f_0^+\end{array}\right.
\end{equation}
The parity condition $f_0^++f_1^++f_0^-+f_1^- \in 2\mz$ (for the class
$h$ to be $0$) is automatically satisfied.
\smallskip

Now we can compute the QHFT invariants of $(S^3,K_0)$ as functions on
$Def(S^3,K)$. Given a $N$-state $s$ of $T$, put $\alpha=s(\hat{2})$,
$\beta=s(\hat{0})$, $\gamma=s(\hat{3})$ and $\delta=s(\hat{1})$, where
$\hat{i}$ denotes the face of $\Delta ^+$ opposite to the $i$th
vertex. By \ref{mon2}) we have:
$$N^{-2}\Hh_N(S^3,K_0,\rho) = \sum_{\alpha,\beta,\gamma,\delta=0}^{N-1}
  \Rr_N(\Delta^+,b^+,w,f^+,c^+)_{\alpha,\beta}^{\gamma,\delta}\
  \Rr_N(\Delta^-,b^-,z,f^-,c^-)_{\beta,\alpha}^{\delta,\gamma}$$
$$\begin{array}{lll} & = & \sum_{\alpha,\beta,\gamma,\delta=0}^{N-1}
    (w_0'^{-c_1^+}w_1'^{c_0^+})^{\frac{N-1}{2}} \
    \frac{g(w_0')}{g(1)}\ \zeta^{\gamma\beta+(m+1)\gamma^2}\
    \omega(w_0',w_1'^{-1}\vert \alpha-\gamma) \ \delta(\alpha+\beta
    -\delta)\\ & & \hspace{4cm} \times
    (z_0'^{-c_1^-}z_1'^{c_0^-})^{\frac{N-1}{2}}\
    \frac{[z_0']g(1)}{g(z_0')}\ \zeta^{-\alpha\delta-(m+1)\delta^2}\
    \frac{\delta(\alpha+\beta
    -\gamma)}{\omega(z_0'/\zeta,z_1'^{-1}\vert \beta-\delta)} \\ \\ &
    = &
    (w_0'^{-c_1^+}w_1'^{c_0^+}z_0'^{-c_1^-}z_1'^{c_0^-})^{\frac{N-1}{2}}
    \frac{[z_0']g(w_0')}{g(z_0')} \sum_{\alpha,\beta=0}^{N-1}
    \zeta^{\beta^2-\alpha^2}\ \omega(w_0',w_1'^{-1}\vert N-\beta)\
    \omega(z_0'/\zeta,z_1'^{-1}\vert N-\alpha)^{-1}.
\end{array}$$
By the proof of Proposition 8.6 in \cite{BB2} we have
$$\frac{[z_0']g(w_0')}{g(z_0')} =
\frac{g((z_0')^*)^*g(w_0')}{|g(1)|^2}$$
and
$$\omega(z_0'/\zeta,z_1'^{-1}\vert N-\alpha)^{-1} =
\omega((z_0')^*,(z_1'^{-1})^*\vert \alpha)^*$$
where $z^*$ is the complex conjugate of $z$. Thus, setting
$$S(w_0',w_1') = \sum_{\beta=0}^{N-1} \zeta^{\beta^2}
\omega(w_0',w_1'^{-1}\vert \beta) = 1+\ \sum_{\beta=1}^{N-1}
\zeta^{\beta^2} \prod_{k=1}^{\beta} \frac{w_1'^{-1}}{1-w_0'\zeta^k}$$
we get
$$\Hh_N(S^3,K_0,\rho) =
N^{2}(w_0'^{-c_1^+}w_1'^{c_0^+}z_0'^{-c_1^-}z_1'^{c_0^-})^{\frac{N-1}{2}}
\frac{g((z_0')^*)^*g(w_0')}{|g(1)|^2}S(w_0',w_1')
\ S((z_0')^*,(z_1')^*)^*.$$ Using (\ref{standardf}) with
$k(m)=k(l)=f_0^+=0$ and the global charge with $c_0^+=c_0^-=0$, we see
that for the complete hyperbolic structure $\rho_{\rm comp}$ on
$S^3\setminus K$ we have
$$z_0' = (w_0')^*= \exp(i\pi/3N),\quad z_1' = (w_1')^*=
\exp(-5i\pi/3N).$$ Hence $$\Hh_N(S^3,K_0,\rho_{\rm comp})
=N^{2}\frac{|g(e^{i\pi/3N})|^2}{|g(1)|^2} \ |
S(e^{i\pi/3N},e^{-5i\pi/3N})|^2.$$ Let us finally consider hyperbolic
$(p,q)$ Dehn filling of $S^3\setminus K$. Denote $S^3(K_{(p,q)})$ the
surgered manifold, $L$ the core of the surgery, and $\rho_{(p,q)}$ its
hyperbolic holonomy. Because of (\ref{surgf}) the difference $f''-f$
is given on the edges $e_0$, $e_1$ and $e_2$ of $\Delta^+$
(resp. $\Delta^-$) by $0$, $r$ and $-r$ (resp. $-2s$, $4s$ and
$0$). Put $N=2m+1$. From Theorem \ref{defscissors} we deduce
$$\begin{array}{l}\Hh_N(S^3(K_{(p,q)}),L,\rho_{(p,q)}) = N^{2}
(w_0'^{-c_1^+}w_1'^{c_0^+}z_0'^{-c_1^-}z_1'^{c_0^-})^{\frac{N-1}{2}}\
\frac{g((z_0')^*)^*g(w_0')}{|g(1)|^2}\\
\hspace{0.8cm}\sum_{\alpha,\beta=0}^{N-1}
  \zeta^{\beta^2-\alpha^2}\omega(w_0',w_1'^{-1}\vert N-\beta)\
  \omega((z_0')^*,(z_1'^{-1})^*\vert \alpha)^*\
  \zeta^{r(N-\beta)(m+1)}\prod_{j=1}^{N-2s}
  \frac{(z_1'^{-1})^*\zeta^{-4s\alpha(m+1)}}{1-(z_0')^*\zeta^{j+\alpha}}\
  .
\end{array}$$
\begin{remark}{\rm Recall the space $\mathcal{C}$ in (\ref{subs}). 
As already mentionned after Definition \ref{non-rough}, we have a
holonomy map $hol:\mathcal{C}\rightarrow X$ to the character variety
$X=X(\pi_1(S^3\setminus K))$. (See \cite{Ril} or \cite{GAM} for a
complete description of the latter). The map $hol$ is generically
$2:1$, and is onto the geometric component of $X$ \cite{Cha}. We can
express the above partition functions in terms of standard generators
of $X$ by the following observation. Considering $S^3\setminus K$ as
the mapping torus of the monodromy $\Phi$, the edges $e_0$, $e_1$ of
$\Delta^+$ are identified with a longitude $l$ and meridian $m$ of the
punctured torus $\Sigma_{1,1}$, and in $\Delta^-$ we have ($e_0'$ is
opposite to $e_0$):
$$e_0=\Phi(l.m),\quad e_1=\Phi(l),\quad e_0'= \Phi(m).$$ As above,
assume that $\rho$ has non trivial holonomy at $m$, $l$ and
$l.m$. Take a flat/charged $\Ii$-triangulation of $(S^3,K_0,\rho)$ as
in Section \ref{FIBERED}, with $PSL(2,\mc)$-valued cocycle $z$. Denote
$z_l$ the value at $l$, and so on. Note that $z_{\Phi(m)}=Az_mA^{-1}$,
where $A=z(S^1)$, the cocycle value on the standard meridian of the
knot $K$. Then the cross ratio moduli of $\Delta^+$ and $\Delta^-$ are
given by
$$w_0=[0:z_l(0):z_lz_m(0):z_{\Phi(m)}(0)],\quad z_0=
[0:z_{\Phi(l)}z_{\Phi(m)}(0):z_{\Phi(l)}z_{\Phi(m)}z_{\Phi(l)}(0):z_mz_l(0)].$$
(We use the branching to remove the twofold ambiguity of $hol$, as it
allows to specify an equivariant association of a fixed point for each
peripheral subgroup of $\rho(\pi_1(S^3\setminus K))$.)}
\end{remark}
\begin{remark}
{\rm Formulas for Cheeger-Chern-Simons invariants
$\Hh_1(S^3,K_0,\rho)$ of arbitrary $PSL(2,\mc)$-characters of
$S^3\setminus K$ come exactly in the same way (See Remark \ref{OKCCS}
and Remark \ref{scissors}). In the peculiar situation of the complete
hyperbolic structure and its hyperbolic Dehn fillings, they coincide
with those of \cite{N}, Section 15.}
\end{remark}


\begin{thebibliography}{99}

\bibitem{Ab} W. Abikoff, \emph{The Real Analytic Theory of
    Teichm\"uller Space}, Lect. Notes Math. 820, Springer (1980)

\bibitem{At} M.F. Atiyah, {\it Topological quantum field theories},
  Publ. Math. IHES 68 (1988) 175--186

\bibitem{B} S. Baseilhac, \emph{Hyperbolic field theories, geometric
  quantization and the dynamics of surface diffeomorphisms}, in
  preparation

\bibitem{BB1} S. Baseilhac, R. Benedetti, \emph{Quantum hyperbolic
    invariants of 3-manifolds with $PSL(2,\C)$-characters}, Topology
  43 (2004) 1373--1423

\bibitem{BB2} S. Baseilhac, R. Benedetti, \emph{Classical and quantum
    dilogarithmic invariants of flat $PSL(2,\C)$-bundles over
    $3$-manifolds}, Geom. Topol. 9 (2005) 493--570

\bibitem{BP} R. Benedetti, C. Petronio, \emph{Lectures on
Hyperbolic Geometry}, Universitext, Springer Verlag, Berlin Heidelberg (1992)

\bibitem{BBo} R. Benedetti, F. Bonsante, \emph{Canonical Wick
    rotations in 3-dimensional gravity}, ArXiv:math.DG/0508485 (v3), 25
    October 2006, monograph to appear on Memoirs of AMS

\bibitem{Bo} F. Bonahon, \emph{Shearing hyperbolic surfaces, bending
    pleated surfaces and Thurston's symplectic form}, Ann. Fac. Sc.
  Toulouse Math. 5 (1996) 233--297

\bibitem{BoL} F. Bonahon, X. Liu, \emph{Representations of the
        quantum Teichmuller space and invariants of surface
        diffeomorphisms} arXiv:math.GT/0407086


\bibitem{CEG} R.D. Canary, D.B.A. Epstein, P. Green, \emph{Notes on
    notes of Thurston}, Analytic and Geometrical aspects of Hyperbolic
  space, D.B.A. Epstein ed., Cambridge University Press, LMS Lecture
  Notes Series, 111 (1987) 3--92

\bibitem{Choi} Y.-E. Choi, \emph{Positively oriented ideal
triangulations on hyperbolic $3$-manifolds}, Topology 43 (2004) 1345--1371



\bibitem{Cha} A. Champanerkar, \emph{A-polynomial and Bloch invariants
of hyperbolic $3$-manifolds}, PhD Thesis, Columbia University (2003)

\bibitem{DZ} J.L. Dupont, C. Zickert, \emph{A dilogarithmic formula
for the Cheeger-Chern-Simons class}, Geom. Topol. 10 (2006) 1347--1372

\bibitem{Ep} D.B.A. Epstein, A. Marden, \emph{Convex hulls in
    hyperbolic space, a theorem of Sullivan and measured pleated
    surfaces}, LMS Lecture Notes Series 111 (1987) 113--253

\bibitem{F1} S. Francaviglia, \emph{Algebraic and geometric solutions
of hyperbolicity equations}, Topol. Appl. 145 (1-3) (2004)
91--118

\bibitem{F2} S. Francaviglia, \emph{Hyperbolic volume of
representations of fundamental groups of cusped $3$-manifolds},
Int. Math. Res. Not. 9 (2004) 425--459

\bibitem{FH} W. Floyd, A. Hatcher, \emph{Incompressible surfaces in
  punctured torus bundles}, Topol. Appl. 13 (1982) 263--282

\bibitem{GAM} F. Gonz\'alez-Acu\~na, J-M. Montesinos-Amilibia,
\emph{On the character variety of group representations in $SL(2,\mc)$
and $PSL(2,\mc)$}, Math. Zeit. 214 (1993) 627--652

\bibitem{HP} M. Heusener, J. Porti, \emph{The variety of characters in
$PSL(2,\C)$}, arXiv:math.GT/0302075

\bibitem{Iv} N.V. Ivanov, \emph{Mapping class groups}, Handbook of
Geometric Topology, Elsevier Science B.V. (2002) 523--633


\bibitem{K2} R. M. Kashaev, \emph{The hyperbolic volume of knots
from the quantum dilogarithm}, Lett. Math. Phys. 39 (1997) 269--275


\bibitem{Kmod} R.M. Kashaev, \emph{Coordinates for the moduli space of
flat $PSL(2,\mr)$-connections}, Math. Res. Lett. 12 (2005) 23--36


\bibitem{La} M. Lackenby, \emph{The canonical decomposition of
    once-punctured torus bundles}, Comment. Math. Helv. 78 (2003)
    363--384

\bibitem{N0} W. D. Neumann, \emph{Combinatorics of triangulations
and the Chern-Simons invariant for hyperbolic $3$-manifolds},
Topology'90 (Columbus, OH, 1990), De Gruyter, Berlin (1992) 243--271

\bibitem{N} W.D. Neumann, {\it Extended Bloch group and the
Cheeger-Chern-Simons class}, Geom. Topol. 8 (2004) 413--474


\bibitem{Pe} R.C. Penner, \emph{The decorated Teichmuller space of
    punctured surfaces}, Comm. Math. Phys. 113 (1987) 299--339

\bibitem{Pe-Po} C. Petronio, J. Porti, \emph{Negatively oriented ideal
triangulations and a proof of Thurston's hyperbolic Dehn filling
theorem}, Expositiones Mathematicae, vol. 18 (2000) 1--35

\bibitem{Ril} R. Riley, \emph{Nonabelian representations of $2$-bridge
knot groups}, Quart. J. Math. Oxford (2), 35 (1984) 191--208

\bibitem{Th} W.P. Thurston, \emph{The geometry and topology of
    $3$-manifolds}, Princeton University lecture notes (1979)

\bibitem{T} V. Turaev, \emph{Quantum Invariants of $3$-Manifolds},
Studies in Mathematics 18, De Gruyter, Berlin (1994)

\bibitem{Tu2} V. Turaev, \emph{Homotopy field theory in dimension 3
    and crossed group-categories}, arXiv:math.GT/0005291

\end{thebibliography}
\end{document}